\documentclass{amsart}
\allowdisplaybreaks

\usepackage{amssymb
,amsthm
,amsmath
,amscd
,mathtools	
}
\usepackage[all]{xy}
\usepackage{graphicx}
\usepackage[top=30truemm,bottom=30truemm,left=25truemm,right=25truemm]{geometry}

\newcommand{\A}{\mathbb{A}}
\newcommand{\Z}{\mathbb{Z}}
\newcommand{\R}{\mathbb{R}}
\newcommand{\C}{\mathbb{C}}
\newcommand{\W}{\mathbb{W}}
\newcommand{\X}{\mathbb{X}}

\newcommand{\CC}{\mathcal{C}}
\newcommand{\TT}{\mathcal{T}}
\newcommand{\Sc}{\mathcal{S}}
\newcommand{\VV}{\mathcal{V}}

\newcommand{\cent}{\mathrm{Cent}}
\newcommand{\disc}{\mathrm{disc}}
\newcommand{\id}{\mathrm{id}}
\newcommand{\mg}{\mathrm{mg}}
\newcommand{\reg}{\mathrm{reg}}
\newcommand{\Int}{\mathrm{Int}}
\newcommand{\inv}{\mathrm{inv}}
\newcommand{\im}{\mathrm{Im}}
\newcommand{\HC}{\mathrm{HC}}
\newcommand{\Irr}{\mathrm{Irr}}
\newcommand{\temp}{\mathrm{temp}}
\newcommand{\Ind}{\mathrm{Ind}}
\newcommand{\tr}{\mathrm{tr}}
\newcommand{\pr}{\mathrm{pr}}
\newcommand{\sgn}{\mathrm{sgn}}
\newcommand{\Hom}{\mathrm{Hom}}
\newcommand{\Gal}{\mathrm{Gal}}
\newcommand{\SL}{\mathrm{SL}}
\newcommand{\GL}{\mathrm{GL}}
\newcommand{\Sp}{\mathrm{Sp}}
\newcommand{\Mp}{\mathrm{Mp}}
\newcommand{\Lie}{\mathrm{Lie}}
\newcommand{\M}{\mathrm{M}}
\newcommand{\U}{\mathrm{U}}

\newcommand{\ep}{\varepsilon}
\newcommand{\lam}{\lambda}
\newcommand{\WD}{{\it WD}}

\renewcommand{\1}{{\bf 1}}
\newcommand{\I}{\sqrt{-1}}
\newcommand{\pair}[1]{\langle #1 \rangle}
\newcommand{\half}[1]{\frac{#1}{2}}
\newcommand{\cl}[1]{\widetilde{#1}}
\newcommand{\iif}{&\quad&\text{if }}
\newcommand{\other}{&\quad&\text{otherwise}}
\newcommand{\resp}{resp.~}
\newcommand{\bs}{\backslash}

\newcommand{\g}{\mathfrak{g}}
\newcommand{\kk}{\mathfrak{k}}
\renewcommand{\tt}{\mathfrak{t}}
\newcommand{\w}{\mathfrak{w}}

\newtheorem{thm}{Theorem}[section]
\newtheorem{lem}[thm]{Lemma}
\newtheorem{prop}[thm]{Proposition}
\newtheorem{cor}[thm]{Corollary}
\newtheorem{defi}[thm]{Definition}
\newtheorem{rem}[thm]{Remark}
\newtheorem{conj}[thm]{Conjecture}
\newtheorem{example}[thm]{Example}
\newtheorem{problem}[thm]{Problem}

\makeatletter
\def\iddots{\mathinner{\mkern1mu\raise\p@
	\hbox{.}\mkern2mu\raise4\p@\hbox{.}\mkern2mu
	\raise7\p@\vbox{\kern7\p@\hbox{.}}\mkern1mu}}
\def\adots{\mathinner{\mkern2mu\raise\p@\hbox{.}
 \mkern2mu\raise4\p@\hbox{.}\mkern1mu
 \raise7\p@\vbox{\kern7\p@\hbox{.}}\mkern1mu}}
\makeatother

\title{On the non-vanishing of theta liftings of tempered representations of $\U(p,q)$}
\author{Hiraku Atobe}
\date{}
\address{Department of mathematics, Kyoto University, Kitashirakawa-Oiwake-cho, Sakyo-ku, Kyoto, 606-8502, Japan}
\email{atobe@math.kyoto-u.ac.jp}

\begin{document}
\maketitle

\begin{abstract}
In this paper, we give an explicit determination of 
the non-vanishing of the theta liftings for unitary dual pairs $(\U(p,q), \U(r,s))$.
Assuming the local Gan--Gross--Prasad conjecture, 
we determine when theta lifts of tempered representations are nonzero 
in terms of the local Langlands correspondence.
\end{abstract}

\tableofcontents

\section{Introduction}\label{intro}
The theta lifting is an important tool in automorphic, real, and $p$-adic representation theories.
In \S \ref{global}, we explain a global motivation of this paper.
The local archimedean introduction is in \S \ref{arch}.
In \S \ref{nonarch}, we recall the non-archimedean result in \cite{AG}.
Finally, in \S \ref{vs}, 
we compare the proof of the archimedean result with 
the one the non-archimedean result.
If a reader were not familiar with the automorphic or $p$-adic representation theory, 
he or she could see only \S \ref{arch}.

\subsection{Global motivation}\label{global}
In the theory of modular forms and automorphic representations, 
one of the most important problem is 
how to construct (non-trivial) cusp forms or cuspidal representations.
\par

Let $k$ be a number field. 
We denote by $\A_k$ the adele ring of $k$.
Let $G$ and $H$ be connected reductive groups over $k$.
In this subsection, 
for an irreducible cuspidal automorphic representation $\pi$ of $G(\A_k)$, 
we shall call $\Pi$ a ``lift'' of $\pi$ to $H(\A_k)$ if
$\Pi$ is an automorphic representation (or ``packet'') of $H(\A_k)$ such that
an arithmetic property (e.g., $L$-function) of $\Pi$ is given explicitly by the one of $\pi$.
There is a general problem.

\begin{problem}\label{prob}
Let $\pi \mapsto \Pi$ be a ``lift'' from 
automorphic representations of $G(\A_k)$ to the ones of $H(\A_k)$.
\begin{enumerate}
\item
Determine when $\Pi$ is nonzero.
\item
Determine when $\Pi$ is irreducible and cuspidal.
\item
Determine the local components $\Pi_v$ of $\Pi$ for each place $v$ of $k$.
\end{enumerate}
\end{problem}

Now we consider this problem for global theta liftings, which
contain many classical liftings, e.g., the Shimura correspondence and the Saito--Kurokawa lifting.
\par

In this paper, we only consider theta liftings for unitary dual pairs.
Let $K/k$ be a quadratic extension of number fields 
and $\A_K$ be the adele ring of $K$.
We denote by $W_n$ an $n$-dimensional hermitian space over $K$.
Also, we fix an anisotropic skew-hermitian space $V_0$ over $K$ of dimension $m_0$, 
and denote by $V_l$ the $m$-dimensional skew-hermitian space obtained from $V_0$ 
by addition of $l$ hyperbolic planes, 
where $m=m_0+2l$.
Let $G=\U(W_n)$ and $H=\U(V_l)$ be 
the isometry groups of $W_n$ and $V_l$, respectively. 
Then $\W = W_n \otimes_K V_l$ can be regarded 
as a $2mn$-dimensional symplectic space over $k$, 
and there exists a canonical map
\[
\alpha_{V,W} \colon \U(W_n)(\A_k) \times \U(V_l)(\A_k) \rightarrow \Sp(\W)(\A_k).
\]
\par

Let $\Mp(\W)_\A$ is the metaplectic cover of $\Sp(\W)(\A_k)$.
Take a complete polarization $\W = \X + \mathbb{Y}$.
Fixing a non-trivial additive character $\psi$ of $\A_k/k$, 
we obtain the Weil representation $\omega_\psi$ of $\Mp(\W)_\A$ associated to $\psi$
which is realized on the space $\Sc(\X(\A_k))$ of Schwarts--Bruhat functions on $\X(\A_k)$.
Let $\Theta \colon \Sc(\X(\A_k)) \rightarrow \C$ be the functional 
$\phi \mapsto \Theta(\phi) = \sum_{x \in \X(k)}\phi(x)$.
There exists a unique canonical splitting $\Sp(\W)(k) \rightarrow \Mp(\W)_\A$
which satisfies that 
$\Theta(\omega_\psi(\gamma)\phi) = \Theta(\phi)$
for any $\gamma \in \Sp(\W)(k)$ and $\phi \in \Sc(\X(\A_k))$.
\par

Fix Hecke characters $\chi_W$ and $\chi_V$ of $\A_K^\times/K^\times$
such that $\chi_W|\A_k^\times = \eta_{K/k}^n$ and $\chi_V|\A_k^\times = \eta_{K/k}^m$, 
where $\eta_{K/k}$ is the quadratic Hecke character of $\A_k^\times/k^\times$ associated to $K/k$.
Kudla gave an explicit splitting
\[
\cl\alpha_{\chi_V, \chi_W} \colon \U(W_n)(\A_k) \times \U(V_l)(\A_k) \rightarrow \Mp(\W)_\A.
\]
Let $\omega_{\psi, V,W} = \omega_\psi \circ \cl\alpha_{\chi_V, \chi_W}$
be the pullback of the Weil representation of $\Mp(\W)_\A$ 
to the product group $\U(W_n)(\A_k) \times \U(V_l)(\A_k)$.
For each $\phi \in \Sc(\X(\A_k))$, we consider the theta function
\[
\Theta(g,h; \phi) = \Theta(\omega_{\psi, V, W}(g,h)\phi)
\]
for $g \in \U(W_n)(\A_k)$ and $h \in \U(V_l)(\A_k)$, 
which is an automorphic form on $\U(W_n)(\A_k) \times \U(V_l)(\A_k)$.
If $\pi$ is an irreducible cuspidal representation of $\U(W_n)(\A_k)$, 
the global theta lift $\theta_{\psi,V_l,W_n}(\pi)$ of $\pi$ to $\U(V_l)(\A_k)$ 
is the automorphic representation of $\U(V_l)(\A_k)$ spanned by the functions
\[
\theta(h; f, \phi) = \int_{\U(W_n)(k) \bs \U(W_n)(\A_k)} \theta(g,h; \phi) \overline{f(g)} dg, 
\]
where $f \in \pi$, $\phi \in \Sc(\X(\A_k))$ and $dg$ is the Tamagawa measure.
\par

We may consider Problem \ref{prob} for $\theta_{\psi,V_l,W_n}(\pi)$.
The cuspidality issue can be solved by the so-called tower property.
\begin{thm}[{Tower property \cite{R}}]
When $\theta_{\psi,V_l,W_n}(\pi)$ is nonzero, 
$\theta_{\psi,V_{l+1},W_n}(\pi)$ is also nonzero.
The theta lift $\theta_{\psi,V_l,W_n}(\pi)$ is cuspidal (possibly zero)
if and only if $\theta_{\psi,V_{j},W_n}(\pi) = 0$ for all $j < l$.
\end{thm}
\par

The irreducibility and the local components are determined by Kudla--Rallis.
\begin{thm}[{\cite[Corollary 7.1.3]{KR}}]
Assume that $\theta_{\psi,V_l,W_n}(\pi)$ is nonzero and cuspidal.
Then $\theta_{\psi,V_l,W_n}(\pi)$ is irreducible and 
$\theta_{\psi,V_l,W_n}(\pi) \cong \otimes_v \theta_{\psi,V_l,W_n}(\pi_v)$, 
where $\theta_{\psi,V_l,W_n}(\pi_v)$ is the ``local theta lift'' of 
the local component $\pi_v$ of $\pi$.
\end{thm}
\par

Therefore, the remaining issue is the non-vanishing problem.
There is a local-global criterion for the non-vanishing 
of the global theta lifting $\theta_{\psi, V_l, W_n}(\pi)$
established by Yamana \cite{Y} and Gan--Qiu--Takeda \cite{GQT}. 
It is roughly stated as follows:
When $\theta_{\psi, V_j, W_n}(\pi)=0$ for any $j < l$ so that 
$\theta_{\psi, V_l, W_n}(\pi)$ is cuspidal, 
\begin{itemize}
\item
the global theta lifting $\theta_{\psi, V_l, W_n}(\pi)$ is nonzero
\end{itemize}
``if'' and only if
\begin{itemize}
\item
the local theta lifting $\theta_{\psi, V_l, W_n}(\pi_v)$ is nonzero for all places $v$ of $k$; 
\item
the ``standard $L$-function $L(s, \pi)$'' of $\pi$ 
is non-vanishing or has a pole at a distinguished point $s_0$.
\end{itemize}
Since a local archimedean property is not proven, 
the ``if'' part is not completely established.
For more precisions, see \cite[Theorem 1.3]{GQT} and a remark after this theorem.
\par

In any case, Problem \ref{prob} for the global theta lifting is reduced to 
a local analogous problem.
The purpose of this paper is to give a criterion for the non-vanishing of local theta liftings.

\subsection{Archimedean case}\label{arch}
The theory of local theta correspondence was initiated by Roger Howe.
Since then, it has been a major theme in the representation theory.
\par

In this paper, we only consider the theta liftings for unitary dual pairs $(\U(p, q), \U(r, s))$.
More precisely, let $W_{p,q}$ and $V_{r,s}$
be a hermitian space and a skew hermitian space over $\C$ of
signature $(p,q)$ and $(r,s)$ and of dimension $n=p+q$ and $m=r+s$, respectively.
Then $\W = W_{p,q} \otimes_\C V_{r,s}$ is a symplectic space over $\R$
of dimension $2mn$.
The isometry group of $W_{p,q}$ (\resp $V_{r,s}$) is isomorphic to
$\U(p,q)$ (\resp $\U(r,s)$).
Fix characters $\chi_W$ and $\chi_V$ of $\C^\times$ such that 
$\chi_W|\R^\times = \sgn^{n}$ and $\chi_V|\R^\times = \sgn^{m}$, 
and a non-trivial additive character $\psi$ of $\R$.
Kudla \cite{Ku2} gave an explicit splitting
\[
\cl\alpha_{\chi_V, \chi_W} \colon \U(p,q) \times \U(r,s) \rightarrow \Mp(\W)
\]
of the natural map
\[
\alpha_{V, W} \colon \U(p,q) \times \U(r,s) \rightarrow \Sp(\W), 
\]
where $\Mp(\W)$ is the metaplectic cover of $\Sp(\W)$. 
We assume that $\chi_W$ and $\chi_V$ depend 
only on $n \bmod 2$ and $m \bmod 2$, respectively.
Let $\omega_\psi$ be the Weil representation of $\Mp(\W)$ associated to $\psi$.
Then we obtain the Weil representation 
$\omega_{\psi,V,W} = \omega_\psi \circ \cl\alpha_{\chi_V, \chi_W}$ of $\U(p,q) \times \U(r,s)$.
For an irreducible (unitary) representation $\pi$ of $\U(p,q)$, 
the maximal $\pi$-isotypic quotient of $\omega_{\psi, V,W}$ has the form
\[
\pi \boxtimes \Theta_{r,s}(\pi)
\]
for some representation $\Theta_{r,s}(\pi)$ of $\U(r,s)$ (known as the big theta lift of $\pi$). 
The most important result is the Howe duality correspondence stated as follows:

\begin{thm}[Howe duality correspondence {\cite[Theorem 2.1]{Ho}}]
If $\Theta_{r,s}(\pi)$ is nonzero, 
then it has a unique irreducible quotient $\theta_{r,s}(\pi)$.
\end{thm}
We shall interpret $\theta_{r,s}(\pi)$ to be zero if so is $\Theta_{r,s}(\pi)$.
We call $\theta_{r,s}(\pi)$ the small theta lift 
(or simply, the local theta lift) of $\pi$ to $\U(r,s)$. 
After the above theorem, it is natural to consider the following two basic problems:
\begin{problem}\label{local}
\begin{enumerate}
\item
Determine precisely when $\Theta_{r,s}(\pi)$ is nonzero.
\item
Determine $\theta_{r,s}(\pi)$ precisely if $\Theta_{r,s}(\pi)$ is nonzero.
\end{enumerate}
\end{problem}

This problem is solved in the following cases:
\begin{itemize}
\item
in the case when $q=0$, i.e., $\U(p,q) = \U(p,0)$ is compact
by Kashiwara--Vergne \cite{KV} (see also \cite[\S 6]{A}); 
\item
in the (almost) equal rank case by Paul \cite{P1}, \cite{P3}; 
\item
in special cases when $\pi$ is discrete series by J.-S. Li \cite{Li} and Paul \cite{P2}.
\end{itemize}
\par

To formulate an answer to Problem \ref{local}, 
it is necessary to have some sort of classifications of irreducible representations 
of the groups $\U(p,q)$ and $\U(r,s)$. 
In this paper, we shall use the local Langlands correspondence (LLC) as a classification, 
and address Problem \ref{local} (1) for tempered representations $\pi$.
\par

The local Langlands correspondence for real reductive connected groups 
is well understood by the work of many mathematicians, including 
Adams, Arthur, Barbasch, Johnson, Langlands, M{\oe}gline, Shelstad, and Vogan.
For discrete series representations, 
it is essentially the same parametrization as using Harish-Chandra parameters
(see Theorem \ref{LLC} (4) below).
In this introduction, we explain the main result (Theorem \ref{main})
only for discrete series representations in terms of their Harish-Chandra parameters.
\par

Let $\pi$ be a discrete series representation of $\U(p,q)$.
Its Harish-Chandra parameter $\lam = \HC(\pi)$ is of the form
\[
\lam = (\lam_1, \dots, \lam_p; \lam'_1, \dots, \lam'_q), 
\]
where $\lam_i, \lam_j' \in \Z + \half{n-1}$, 
$\lam_1 > \dots > \lam_p$, $\lam'_1 > \dots > \lam'_q$, and
$\{\lam_1, \dots, \lam_p\} \cap \{\lam'_1, \dots, \lam'_q\} = \emptyset$.
To state our main result, we define some terminologies.

\begin{defi}[{Definition \ref{def}}]\label{intro.def}
Fix $\kappa \in \{1,2\}$ and 
suppose that $\chi_V$ is of the form 
\[
\chi_V(ae^{\I\theta}) = e^{\nu\I\theta}
\]
for $a > 0$ and $\theta \in \R/2\pi\Z$ with $\nu \in \Z$ 
such that $\nu \equiv n+\kappa \bmod 2$.
Let $\pi$ be a discrete series representation of $\U(p,q)$ with
Harish-Chandra parameter $\lam = (\lam_1, \dots, \lam_p; \lam'_1, \dots, \lam'_q)$.

\begin{enumerate}
\item
Write
\[
\lam-\left(\half{\nu}, \dots, \half{\nu}\right) = 
\left(a_1, \dots, a_x, 
\underbrace{\half{k-1}, \half{k-3}, \dots, \half{-k+1}}_k, 
b_1, \dots b_y; 
c_1, \dots, c_z, 
d_1, \dots, d_w\right)
\]
or
\[
\lam-\left(\half{\nu}, \dots, \half{\nu}\right) = 
\left(a_1, \dots, a_x, 
b_1, \dots b_y; 
c_1, \dots, c_z, 
\underbrace{\half{k-1}, \half{k-3}, \dots, \half{-k+1}}_k,
d_1, \dots, d_w\right)
\]
with $a_x,-b_1, c_z, -d_1 \geq \half{k+1}$.
If $k = 0$ and $a_i, b_i, c_i, d_i \in \Z \setminus \{0\}$ for all $i$, we set $k_\lam = -1$.
Otherwise, we set $k_\lam$ to be the maximal choice of $k$.

\item
When $k$ is chosen to be maximal in (1), 
we put $r_\lam = x+w$ and $s_\lam = y+z$.

\item
Define $X_\lam \subset \half{1}\Z \times \{\pm1\}$ by 
\[
X_\lam = \left\{ \left(\lam_1-\half{\nu}, +1\right), \dots, \left(\lam_p-\half{\nu}, +1\right) \right\}
\cup
\left\{ \left(\lam'_1-\half{\nu}, -1\right), \dots, \left(\lam'_q-\half{\nu}, -1\right) \right\}.
\]

\item
We define a sequence 
$X_\lam=X_\lam^{(0)} \supset X_\lam^{(1)} \supset \dots 
\supset X_\lam^{(n)} \supset \cdots$ as follows:
Let $\{\beta_1, \dots, \beta_{u_j}\}$ be the image of $X_\lam^{(j)}$ 
under the projection $\half{1}\Z \times \{\pm1\} \rightarrow \half{1}\Z$
such that $\beta_1 > \dots > \beta_{u_j}$.
Set $S$ to be the set of $i \in \{2, \dots, u_j\}$ such that
one of the following holds:
\begin{itemize}
\item
$\beta_{i-1} \in \{a_1, \dots, a_x\}$ and $\beta_i \in \{c_1, \dots, c_z\}$;
\item
$\beta_{i-1} \in \{b_1, \dots, b_y\}$ and $\beta_i \in \{d_1, \dots, d_w\}$.
\end{itemize}
Here, we assume that $k$ is chosen to be maximal in (1).
Then we define a subset $X_\lam^{(j+1)}$ of $X_\lam^{(j)}$ by
\[
X_\lam^{(j+1)} 
= X_\lam^{(j)} \setminus \left(\bigcup_{i \in S}\{(\beta_{i-1}, +1), (\beta_i, -1)\}\right).
\]
Finally, we set $X_\lam^{(\infty)} = X_\lam^{(n)} = X_\lam^{(n+1)}$.

\item
For an integer $T$ and $\epsilon \in \{\pm1\}$, we define a set $\CC^\epsilon_\lam(T)$ by 
\[
\CC^\epsilon_\lam(T) = \left\{(\alpha, \epsilon) \in X_\lam^{(\infty)} |\ 
0 \leq \epsilon\alpha + \half{k_\lam-1} < T
\right\}.
\]
In particular, if $T \leq 0$, then $\CC^\epsilon_\lam(T) = \emptyset$.

\end{enumerate}
\end{defi}

The main result for discrete series representations is stated as follows:
\begin{thm}[{Theorem \ref{main}}]\label{intro.main}
Let $\pi$ be a discrete series representation of $\U(p,q)$ 
with Harish-Chandra parameter $\lam$.
Set $k=k_\lam$, $r=r_\lam$ and $s=s_\lam$.
\begin{enumerate}
\item
Suppose that $k = -1$.
Then for integers $l$ and $t \geq 0$, 
the theta lift $\Theta_{r+2t+l+1, s+l}(\pi)$ is nonzero if and only if
\[
l \geq 0
\quad\text{and} \quad
\#\CC^\epsilon_\lam(t+l) \leq l
\quad\text{for each $\epsilon \in \{\pm1\}$}.
\]

\item
Suppose that $k \geq 0$.
Then for integers $l$ and $t \geq 1$, 
the theta lift $\Theta_{r+2t+l, s+l}(\pi)$ is nonzero if and only if
\[
l \geq k
\quad\text{and} \quad
\#\CC^\epsilon_\lam(t+l) \leq l
\quad\text{for each $\epsilon \in \{\pm1\}$}.
\]
Moreover, for an integer $l$, 
the theta lift $\Theta_{r+l, s+l}(\pi)$ is nonzero if and only if $l \geq 0$.
\end{enumerate}
\end{thm}

In fact, we use the local Gan--Gross--Prasad conjecture (Conjecture \ref{GGP} below)
to prove our main result (Theorem \ref{main}).
This conjecture gives a conjectural answer to restriction problems 
in terms of the local Langlands correspondence.
For more precisions, see \S \ref{sec.GGP} below.
To prove our main result for discrete series representations (Theorem \ref{intro.main}), 
we use the local Gan--Gross--Prasad conjecture only for discrete series representations. 
This case has been established by He \cite{He} (in terms of Harish-Chandra parameters), 
so that the statement in Theorem \ref{intro.main} holds unconditionally.
For tempered representations, 
only a weaker version of the local Gan--Gross--Prasad conjecture
is proven by Beuzart-Plessis \cite{BP}.
\par

The proof of the main theorem consists of three steps:
First, we show the sufficient condition for the non-vanishing of theta lifts when $t \geq 1$.
This is proven by induction by using the local Gan--Gross--Prasad conjecture 
and seesaw identities (\S \ref{nonvanish}).
The initial steps of the induction are the equal rank case and the almost equal rank case, 
which are established by Paul (\cite{P1}, \cite{P3}).
Next, we show the necessary condition for the non-vanishing of theta lifts when $t \geq 1$.
This is proven by using a non-vanishing result of theta lifts of 
one dimensional representations and seesaw identities (\S \ref{vanish}).
Finally, by using the conservation relation (Theorem \ref{CR} below), 
we obtain the $t=0$ case.
\par

When $\chi_V$ is of the form $\chi_V(ae^{\I\theta}) = e^{\nu\I\theta}$
for $a > 0$ and $\theta \in \R/2\pi\Z$ with $\nu \in \Z$ such that $\nu \equiv m \bmod 2$, 
it is known that $\Theta_{r,s}(\pi) \not= 0$ 
if and only if $\Theta_{s,r}(\pi^\vee \otimes {\det}^{\nu}) \not= 0$
(see Proposition \ref{vee} below).
By this property together with our main result, 
we can determine the non-vanishing of all theta lifts of tempered representations.
\par

We shall give some examples.
First, we consider the case when $\U(p,q)$ is compact, i.e, $p=0$ or $q=0$.
\begin{example}[Compact case]
Suppose that $p=n$ and $q=0$.
Let $\pi$ an irreducible (continuous, discrete series) representation $\pi$ of $\U(n) = \U(n,0)$
with Harish-Chandra parameter $\lam$ satisfying
\[
\lam-\left(\half{\nu}, \dots, \half{\nu}\right) = 
\left(a_1, \dots, a_x, 
\underbrace{\half{k-1}, \half{k-3}, \dots, \half{-k+1}}_k, 
b_1, \dots b_y\right), 
\]
where $k = \max\{0, k_\lam\}$, $a_x, -b_1 \geq (k_\lam+1)/2$, $x+y+k=n$, 
and $\nu \in \Z$ such that 
$\chi_V(ae^{\I\theta}) = e^{\nu\I\theta}$ for $a > 0$ and $\theta \in \R/2\pi\Z$.
Note that $r_\lam=x$ and $s_\lam=y$.
Then 
\[
X_\lam = X_\lam^{(\infty)} = \left\{
(a_1, +1), \dots, (a_x, +1), \left(\half{k-1}, +1\right), \dots, \left(\half{-k+1}, +1\right),
(b_1, +1), \dots. (b_y,+1)
\right\}.
\]
Hence $\CC^-_\lam(T) = \emptyset$ for any $T$.
Moreover, for $t \geq 1$ and $l \geq k$, we see that $\#\CC^+_\lam(t+l) \leq l$ 
if and only if $l-k\geq x$, or $l-k<x$ and
\[
a_{x-l+k} \geq t+l-\half{k_\lam-1}.
\]
Similarly, $\pi^\vee \otimes {\det}^\nu$ has the Harish-Chandra parameter $\lam^\vee$ satisfying
\[
\lam^\vee-\left(\half{\nu}, \dots, \half{\nu}\right) = 
\left(-b_y, \dots, -b_1, 
\underbrace{\half{k-1}, \half{k-3}, \dots, \half{-k+1}}_k, 
-a_x, \dots -a_1 \right).
\]
Hence $\CC^-_{\lam^\vee}(T) = \emptyset$ for any $T$. 
Moreover, for $t \geq 1$ and $l \geq k$, we see that $\#\CC^+_\lam(t+l) \leq l$ 
if and only if $l-k \geq y$, or $l-k<y$ and 
\[
b_{1+l-k} \leq -\left(t+l-\half{k_\lam-1}\right).
\]
Therefore, by Theorem \ref{intro.main}, 
$\Theta_{r,s}(\pi)$ is nonzero if and only if one of the following holds:
\begin{itemize}
\item
$|(r-s) - (x-y)| \leq 1$, and $r \geq x$ and $s \geq y$;
\item
$(r-s)-(x-y) > 1$, and $s \geq n-x$, and if $n-x \leq s \leq n-1$, then
\[
a_{n-s} \geq \half{m-n+1};
\]
\item
$(r-s)-(x-y) < -1$, and $r \geq n-y$, and if $n-y \leq r \leq n$, then
\[
b_{1+y+r-n} \leq -\half{m-n+1}.
\]
\end{itemize}
Here, we put $m=r+s$.
It is easy to check that this condition is compatible with \cite[Proposition 6.6]{A}.
\end{example}

Second, we consider the case when $\pi$ is holomorphic discrete series.
\begin{example}[Holomorphic case]
Suppose that $p = q$ so that $n = 2p$. 
Let $k \geq n$ be an even integer.
Consider the holomorphic discrete series representation $\pi$ of $\U(p,p)$ 
of weight $k$.
Its Harish-Chandra parameter $\lam$ is given by
\[
\lam = 
\left( \half{k-1}, \half{k-3}, \dots, \half{k-n+1}; \half{-k+n-1}, \dots, \half{-k+3}, \half{-k+1}\right).
\]
For simplicity, we assume that $\chi_V$ is trivial, i.e., $\nu=0$.
Then we have
$k_\lam = 0$, $r_\lam = n$, $s_\lam = 0$ and
\[
X_\lam = X_\lam^{(\infty)} = \left\{
\left(\half{k-1}, +1\right), \dots, \left(\half{k-n+1}, +1\right), 
\left(\half{-k+n-1}, -1\right), \dots, \left(\half{-k+1}, -1\right)
\right\}.
\]
Hence for integer $T$, we have
\[
\#\CC_\lam^+(T) = \#\CC_\lam^-(T) = \left\{
\begin{aligned}
&0 \iif T \leq \half{k-n},\\
&T-\half{k-n} \iif \half{k-n} \leq T \leq \half{k},\\
&p \iif T \geq \half{k}.
\end{aligned}
\right.
\]
On the other hand, the Harish-Chandra parameter $\lam^\vee$ of $\pi^\vee$ is given by
\[
\lam^\vee = 
\left( \half{-k+n-1}, \dots, \half{-k+3}, \half{-k+1}; \half{k-1}, \half{k-3}, \dots, \half{k-n+1} \right)
\]
so that $\CC_{\lam^\vee}^\epsilon(T) = \emptyset$ 
for any $T \in \Z$ and $\epsilon \in \{\pm1\}$.
Therefore, by Theorem \ref{intro.main}, 
$\Theta_{r,s}(\pi) \not= 0$ with $r+s \equiv 0 \bmod 2$ 
if and only if one of the following holds:
\begin{itemize}
\item
$r-s \leq n$ and $r \geq n$;
\item
$n \leq r-s \leq k$ and $s \geq 0$;
\item
$r-s > k$ and $s \geq p$.
\end{itemize}
The following summarizes when $\Theta_{r,s}(\pi) \not=0$
with $n=2p=4$ and $k=8$.
Here, a black plot $(r,s)$ means that $\Theta_{r,s}(\pi) \not= 0$, 
and a white plot $(r,s)$ means that $\Theta_{r,s}(\pi) = 0$.

\begin{center}
\includegraphics{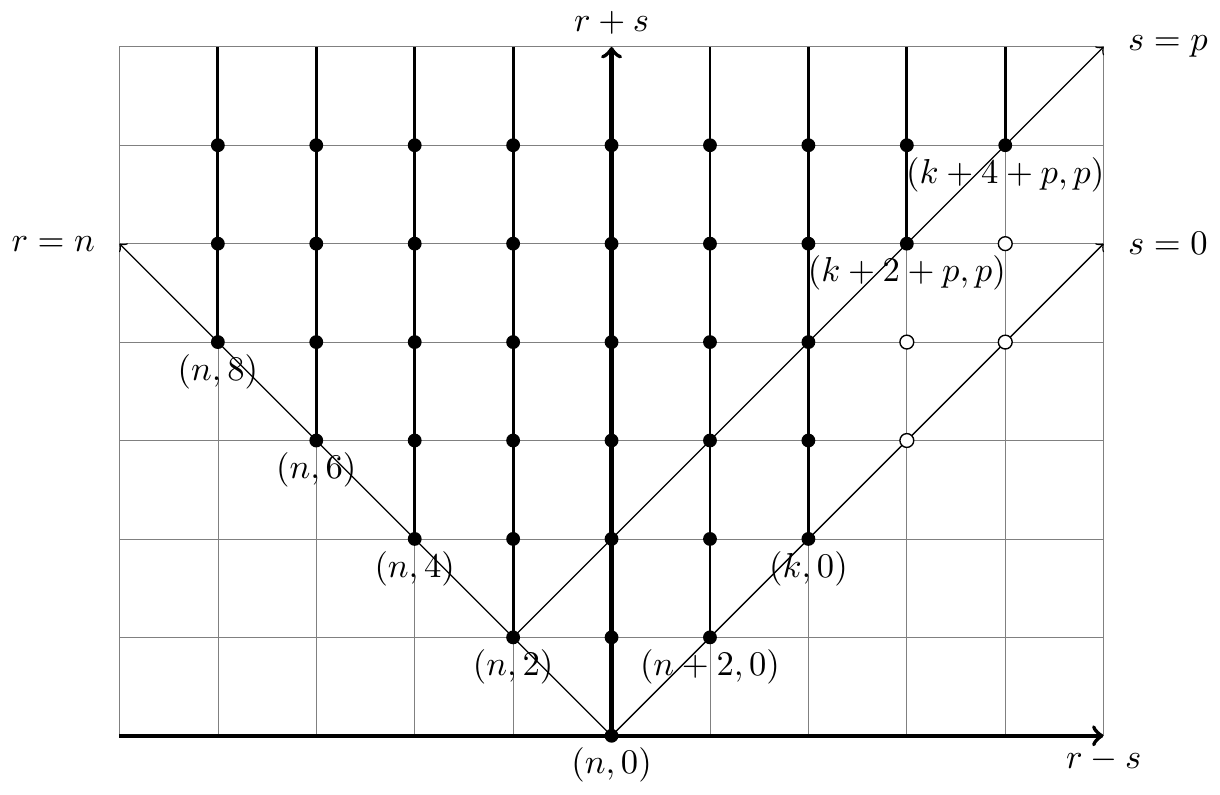}
\end{center}

\end{example}

Third, we consider the case when $\pi$ is generic (large) discrete series.
\begin{example}[Generic case]
Suppose that $p=q+1$ so that $n=2p-1$.
Consider a generic (large) discrete series representation $\pi$.
Its Harish-Chandra parameter 
\[
\lam = (a_1, \dots, a_p; b_1, \dots, b_q)
\]
satisfies that
\[
a_1 > b_1 > a_2 > b_2 > \dots > a_{q} > b_q > a_{p}.
\]
Hence
\begin{align*}
X_\lam = \left\{ \left(a_1-\half{\nu}, +1 \right), \left(b_1-\half{\nu},-1\right), \dots, 
\left(a_q-\half{\nu}, +1\right), \left(b_q-\half{\nu},-1\right), \left(a_p-\half{\nu}, +1\right) \right\}.
\end{align*}
In particular, we have
\begin{align*}
k_\lam &= \left\{
\begin{aligned}
&1 \iif \nu \in \{2a_1, \dots, 2a_p\} \cup \{2b_1, \dots, 2b_q\},\\
&0 \iif \nu \not\equiv 2a_i \equiv 2b_j,\\
&-1 \other,
\end{aligned}
\right.
\\
r_\lam &= \#\{i \in \{1, \dots, p\}\ |\ 2a_p > \nu\} + \#\{j \in \{1, \dots, q\}\ |\ 2b_j < \nu\}, \\
s_\lam &= \#\{i \in \{1, \dots, p\}\ |\ 2a_p < \nu\} + \#\{j \in \{1, \dots, q\}\ |\ 2b_j > \nu\}.
\end{align*}
Now we further assume that $\nu = 2a_{i_0}$ for some $1 \leq i_0 \leq p$. 
Then 
$k_\lam = 1$, $r_\lam = s_\lam = q$, and
\[
X_\lam^{(1)} = X_\lam^{(\infty)} \subset 
\left\{
(0,+1), \left(b_{i_0}-\half{\nu}, -1\right), \left(a_p-\half{\nu}, +1\right)
\right\}, 
\]
which is equal if $i_0 < p$.
Hence $\#\CC_\lam^\epsilon(T) \leq 1$ for any $T \in \Z$ and $\epsilon \in \{\pm1\}$.
Similarly, the Harish-Chandra parameter $\lam^\vee$ of $\pi^\vee \otimes {\det}^\nu$ 
satisfies that
\begin{align*}
X_{\lam^\vee} &= \left\{ \left(-a_p+\half{\nu}, +1 \right), 
\left(-b_q+\half{\nu},-1\right), \left(-a_q+\half{\nu}, +1 \right)
\dots, 
\left(-b_1+\half{\nu},-1\right), \left(-a_1+\half{\nu}, +1\right) \right\},\\
X_{\lam^\vee}^{(\infty)} &\subset \left\{
(0,+1), \left(-b_{i_0-1}+\half{\nu}, -1\right), \left(-a_1+\half{\nu}, +1\right)
\right\}, 
\end{align*}
which is equal if $i_0 > 1$.
Hence $\#\CC_{\lam^\vee}^\epsilon(T) \leq 1$ for any $T \in \Z$ and $\epsilon \in \{\pm1\}$.
Therefore, by Theorem \ref{intro.main}, 
$\Theta_{r,s}(\pi) \not= 0$ with $r+s \equiv p+q+1 \bmod 2$ 
if and only if $(r,s) = (q,q)$ or $\min\{r-q, s-q\} \geq 1$.
The following summarizes when $\Theta_{r,s}(\pi) \not=0$.
Here, a black plot $(r,s)$ means that $\Theta_{r,s}(\pi) \not= 0$, 
and a white plot $(r,s)$ means that $\Theta_{r,s}(\pi) = 0$.

\begin{center}
\includegraphics{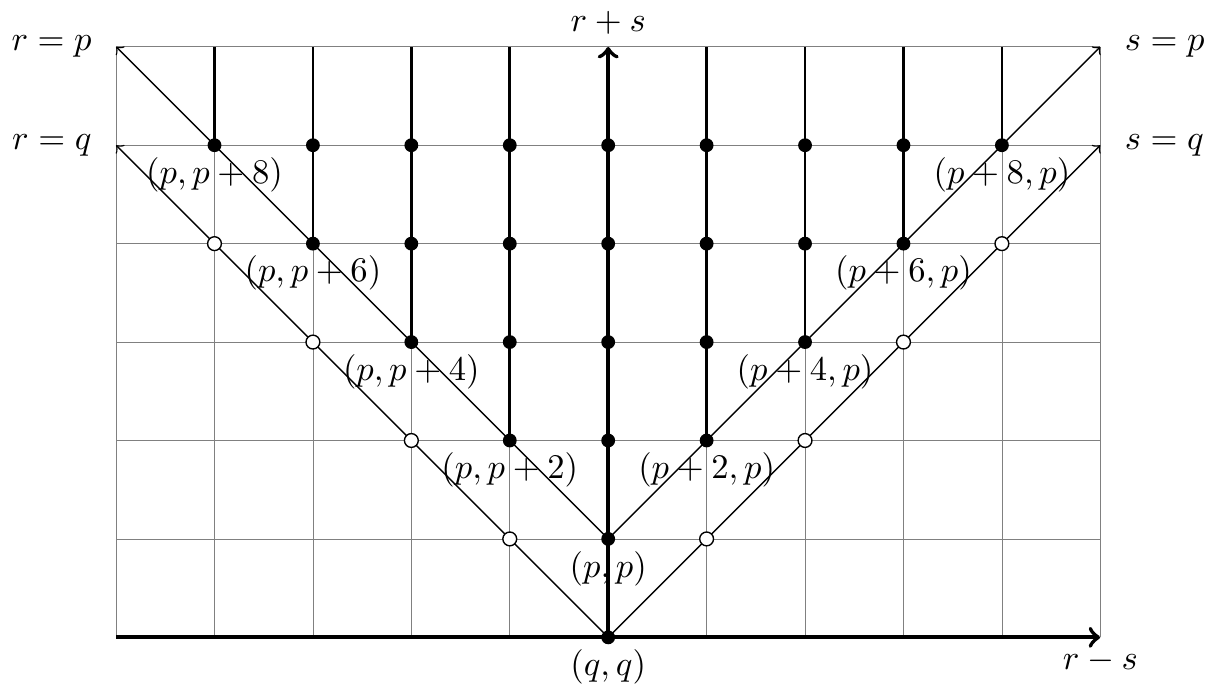}
\end{center}
\end{example}

Finally, we give a more specific example.
\begin{example}
Suppose that $\chi_V$ is trivial, i.e., $\nu=0$.
Let $\pi$ be an irreducible discrete series representation of $\U(4,5)$
with Harish-Chandra parameter
\[
\lam = (6,5,4,-8; 3,1,0,-3,-7).
\]
Then $k_\lam = 1$, $r_\lam = 5$, $s_\lam = 3$, and
\begin{align*}
X_\lam &= \{ (6, +1), (5, +1), (4, +1), (3,-1), (1, -1), (0,-1), (-3, -1), (-7, -1), (-8, +1)\}, \\
X_\lam^{(\infty)} &= \{ (6, +1), \phantom{(5, +1), (4, +1), (3,-1), (1, -1),} 
(0,-1), (-3, -1), (-7, -1), (-8, +1)\}.
\end{align*}
Hence
\[
\CC^+_\lam(T) = \left\{
\begin{aligned}
&\emptyset \iif T \leq 6,\\
&\{(6,+1)\} \iif T>6,
\end{aligned}
\right.
\quad
\CC^-_\lam(T) = \left\{
\begin{aligned}
&\{(0,-1)\}	\iif 0< T \leq 3,\\
&\{(0,-1), (-3,-1)\} \iif 3 < T \leq 7,\\
&\{(0,-1), (-3,-1), (-7,-1)\} \iif T>7.
\end{aligned}
\right.
\]
Therefore, by Theorem \ref{intro.main}, 
$\Theta_{5+2t+l,3+l}(\pi)$ is nonzero if and only if one of the following holds:
\begin{itemize}
\item
$t=0$ and $l \geq 0$;
\item
$t \geq 1$, $l \geq 1$ and $t+1\leq 3$;
\item
$t \geq 1$, $l \geq 2$ and $t+2\leq 7$;
\item
$t \geq 1$ and $l \geq 3$.
\end{itemize}
Similarly, the Harish-Chandra parameter of $\pi^\vee$ is given by
\[
\lam^\vee = (8, -4,-5,-6; 7,3,0,-1,-3).
\]
We have $k_{\lam^\vee} = 1$, $r_{\lam^\vee}=3$, $s_{\lam^\vee} = 5$, and
\begin{align*}
X_{\lam^\vee} &= \{(8,+1), (7,-1), (3,-1), (0,-1), (-1,-1), (-3,-1), (-4,+1), (-5,+1),(-6,+1)\},\\
X_{\lam^\vee}^{(\infty)} &= \{\phantom{(8,+1), (7,-1),}  (3,-1), (0,-1), (-1,-1), (-3,-1), (-4,+1), (-5,+1),(-6,+1)\}.
\end{align*}
Hence $\CC^+_{\lam^\vee}(T) = \emptyset$ for any $T$, and
\[
\CC^-_{\lam^\vee}(T) = \left\{
\begin{aligned}
&\{(0,-1)\} \iif 0 < T \leq 1,\\
&\{(0,-1), (-1,-1)\} \iif 1 < T \leq 3,\\
&\{(0,-1), (-1,-1), (-3,-1)\} \iif T > 3.
\end{aligned}
\right.
\]
Therefore, by Theorem \ref{intro.main}, 
$\Theta_{5+l,3+2t+l}(\pi)$ is nonzero if and only if one of the following holds:
\begin{itemize}
\item
$t=0$ and $l \geq 0$;
\item
$t \geq 0$, $l \geq 1$ and $t+1 \leq 1$;
\item
$t \geq 0$, $l \geq 2$ and $t+2 \leq 3$;
\item
$t \geq 0$ and $l \geq 3$.
\end{itemize}
The following summarizes when $\Theta_{r,s}(\pi) \not=0$.
Here, a black plot $(r,s)$ means that $\Theta_{r,s}(\pi) \not= 0$, 
and a white plot $(r,s)$ means that $\Theta_{r,s}(\pi) = 0$.

\begin{center}
\includegraphics{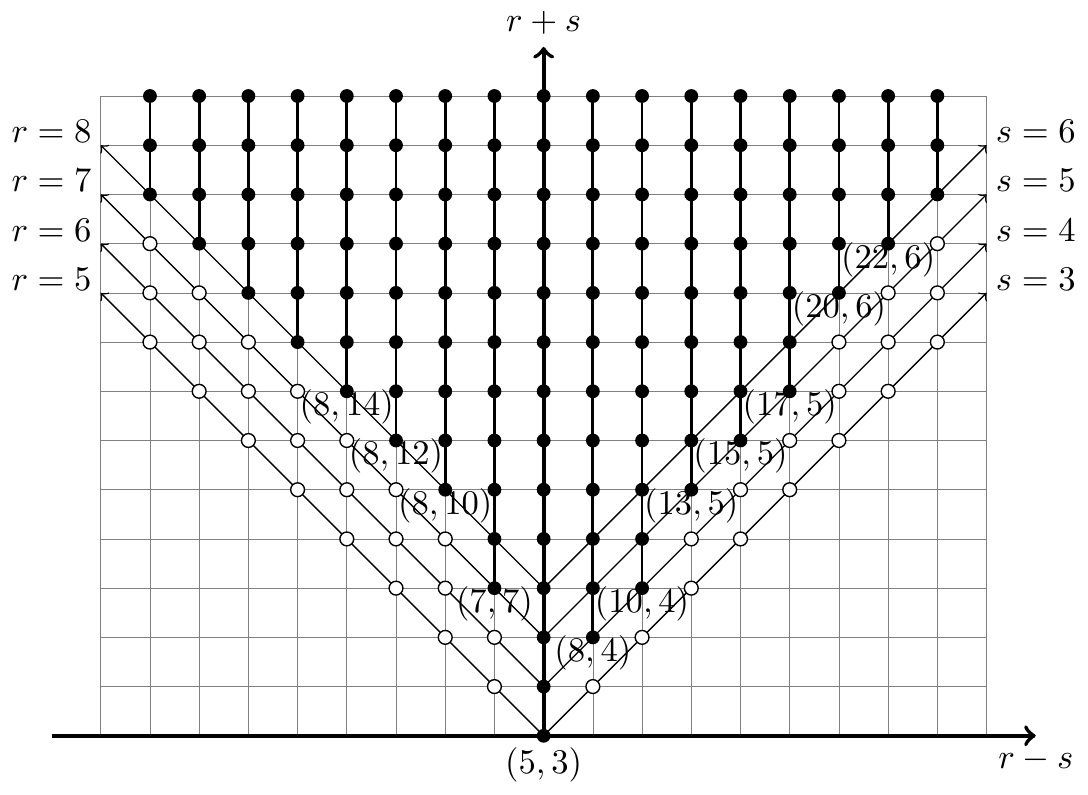}
\end{center}

\end{example}

\subsection{Non-archimedean case}\label{nonarch}
In this subsection, we recall a result in the non-archimedean case in \cite{AG}.
Let $E/F$ be a quadratic extension of non-archimedean local fields of characteristic zero.
We denote by $W_n$ an $n$-dimensional hermitian space over $E$.
Also, we fix an anisotropic skew-hermitian space $V_0$ over $E$ of dimension $m_0$, 
and denote by $V_l$ the $m$-dimensional skew-hermitian space obtained from $V_0$ 
by addition of $l$ hyperbolic planes, where $m=m_0+2l$.
As in \S \ref{global}, fixing a non-trivial additive character $\psi$ of $F$ and
characters $\chi_W$ and $\chi_V$ of $E^\times$ such that
$\chi_W|F^\times = \eta_{E/F}^n$ and $\chi_V|F^\times = \eta_{E/F}^m$, 
we obtain the Weil representation $\omega_{\psi, V, W}$ of $\U(W_n) \times \U(V_l)$.
Here, $\eta_{E/F}$ is the quadratic character of $F^\times$ associated to $E/F$.
\par

For an irreducible smooth representation $\pi$ of $\U(W_n)$, 
the big theta lift $\Theta_{\psi, V_l, W_n}(\pi)$ is defined similarly as in \S \ref{arch}.
The Howe duality correspondence is proven by Waldspurger \cite{W} 
when the residual characteristic of $F$ is not two, 
and by Gan--Takeda \cite{GT1}, \cite{GT2} in general.
Hence we can define the small theta lift (or simply, the local theta lift) 
$\theta_{\psi, V_l, W_n}(\pi)$ of $\pi$ similarly as in \S \ref{arch}, 
and we may consider Problem \ref{local} in the non-archimedean case.
\par

In \cite{AG}, we addressed both Problem \ref{local} (1) and (2) 
in terms of the local Langlands correspondence established by Mok \cite{Mo} 
and Kaletha--M{\'i}nguez--Shin--White \cite{KMSW}.
We recall only the result concerning Problem \ref{local} (1).
\par

The local Langlands correspondence classifies 
the irreducible tempered representations $\pi$ of $\U(W_n)$
by their $L$-parameters $\lam = (\phi,\eta)$, 
where 
\[
\phi \colon \WD_E \rightarrow \GL_n(\C)
\]
is a conjugate self-dual representation of the Weil--Deligne group 
$\WD_E = W_E \times \SL_2(\C)$ of sign $(-1)^{n-1}$ with bounded image, and
\[
\eta \in \Irr(A_\phi)
\]
is an irreducible character of the component group $A_\phi$ associated to $\phi$, 
which is an elementary two abelian group.
More precisely, if we decompose 
\[
\phi = m_1\phi_1 \oplus \dots \oplus m_u\phi_u \oplus (\phi' \oplus {}^c\phi'^\vee), 
\]
where
\begin{itemize}
\item
$\phi_i$ is an irreducible conjugate self-dual representation of sign $(-1)^{n-1}$; 
\item
$m_i > 0$ is the multiplicity of $\phi_i$ in $\phi$; 
\item
$\phi'$ is a sum of irreducible representations of $\WD_E$ which are
not conjugate self-dual of sign $(-1)^{n-1}$; 
\item
${}^c\phi'^\vee$ is the conjugate dual of $\phi'$, 
\end{itemize}
then 
\[
A_\phi = (\Z/2\Z) e_{\phi_1} \oplus \dots \oplus (\Z/2\Z) e_{\phi_u}.
\]
Namely, $A_\phi$ is a free $(\Z/2\Z)$-module of rank $u$ 
equipped with a canonical basis $\{e_{\phi_1}, \dots, e_{\phi_u}\}$
associated to $\{\phi_1, \dots, \phi_u\}$.
We call the element $z_\phi = m_1e_{\phi_1} + \dots + m_ue_{\phi_u}$ in $A_\phi$
the central element of $A_\phi$.
\par

The tower property (Proposition \ref{tower} below) also holds in the non-archimedean case.
It asserts that if $\Theta_{\psi, V_l, W_n}(\pi)$ is nonzero, 
then so is $\Theta_{\psi, V_{l+1}, W_n}(\pi)$.
Fix $\kappa \in \{1,2\}$ and $\delta \in E^\times$ such that $\tr_{E/F}(\delta) = 0$. 
There are exactly two anisotropic skew-hermitian spaces $V_0^+$ and $V_0^-$ such that
$\dim(V_0^\pm) \equiv n+\kappa \bmod 2$.
For $\epsilon \in \{\pm1\}$, 
when $V_0 = V_0^\epsilon$, we also write $V_l = V_l^\epsilon$.
This is characterized by 
\[
\chi_V(\delta^{-m}(-1)^{(m-1)m/2}\det(V_l)) = \epsilon.
\]
We call $\VV_\epsilon = \{V_l^\epsilon\ |\ l \geq 0\}$ the $\epsilon$-Witt tower, 
and
\[
m_\epsilon(\pi) = \min\{\dim(V_l^\epsilon)\ |\ \Theta_{\psi, V_l^\epsilon, W_n}(\pi) \not = 0\}
\]
the first occurrence index of $\pi$ in the Witt tower $\VV_\epsilon$.
\par

The following theorem is one of the main results in \cite{AG}.
Here, we denote the unique $k$-dimensional irreducible representation of $\SL_2(\C)$ by $S_k$.
\begin{thm}[{\cite[Theorem 4.1]{AG}}]\label{AGmain}
Fix $\kappa \in \{1,2\}$ and $\delta \in E^\times$ such that $\tr_{E/F}(\delta) = 0$. 
Let $\pi$ be an irreducible tempered representation of $\U(W_n)$ 
with $L$-parameter $\lam = (\phi,\eta)$.
\begin{enumerate}
\item
Consider the set $\TT$ containing $\kappa-2$ 
and all integers $k >0$ with $k \equiv \kappa \bmod 2$
satisfying the following conditions:
\begin{description}
\item [(chain condition)] 
$\phi$ contains $\chi_V S_{\kappa} + \chi_V S_{\kappa+2} + \dots + \chi_V S_k$;
\item [(odd-ness condition)] 
the multiplicity of $\chi_V S_{r}$ in $\phi$ is odd for $r=\kappa, \kappa+2, \dots, k-2$;
\item [(initial condition)] 
if $\kappa=2$, then $\eta(e_{V, 2})=-1$;
\item [(alternating condition)] 
$\eta(e_{V, r})=-\eta(e_{V, r+2})$ for $r=\kappa, \kappa+2, \dots, k-2$.
\end{description}
Here, $e_{V, r}$ is the element in $A_\phi$ corresponding to $\chi_V S_r$. 
Let
\[  
k_\lam = \max \, \TT.  
\]
Then
\[  
\min\{m_+(\pi), m_-(\pi)\} =  n - k_\lam
\quad \text{and} \quad   
\max\{m_+(\pi), m_-(\pi)\} = n + 2 + k_\lam. 
\]

\item
If $k_\lam = -1$, then $m_+(\pi)=m_-(\pi)$.
Suppose that $k_\lam \geq 0$.
Then $\phi$ contains $\chi_V$ if $\kappa=1$.
Moreover, $\min\{m_+(\pi), m_-(\pi)\} = m_\alpha(\pi)$ if and only if 
\[
\alpha = \left\{
\begin{aligned}
&\eta(z_\phi+e_{V, 1}) \iif \kappa=1,\\
&\eta(z_\phi) \cdot \ep(\phi \otimes \chi_V^{-1}, \psi_2^E) \iif \kappa=2, 
\end{aligned}
\right.
\]
where $\ep(\phi \otimes \chi_V^{-1}, \psi_2^E) = \ep(1/2, \phi \otimes \chi_V^{-1}, \psi_2^E)$
is the root number of $\phi \otimes \chi_V^{-1}$ with respect to 
the additive character $\psi_2^E$ of $E$ defined by $\psi_2^E(x) = \psi(\tr_{E/F}(\delta x))$.
\end{enumerate}
\end{thm}

It seems that this theorem closely resembles our main result (Theorem \ref{main} below).

\subsection{Archimedean case v.s. non-archimedean case}\label{vs}
In this subsection, we explain some differences 
between the archimedean case and the non-archimedean case.
\par

First of all, one of clear differences is the parametrization of Witt towers.
Fixing $\kappa \in \{1,2\}$, 
in the non-archimedean case, there exist exactly two Witt towers 
\[
\VV_\pm = \{V_l^\pm\ |\ l \geq 0,\ \dim(V_l^\pm) \equiv n + \kappa \bmod2\}, 
\]
whereas, in the archimedean case, for each integer $d$ with $d \equiv \kappa \bmod 2$, 
there is a Witt tower
\[
\VV_d = \{V_{r,s}\ |\ r-s = d,\ r+s \equiv n + \kappa \bmod 2\}.
\]
The first occurrence indices are defined similarly in \S \ref{nonarch} for each Witt towers.
The conservation relation proven by Sun--Zhu \cite{SZ2} asserts that 
$m_+(\pi) + m_-(\pi) = 2n+2$ for any irreducible smooth representation $\pi$ of $\U(W_n)$
in the non-archimedean case, 
whereas, 
it is more complicated in the archimedean case (see Theorem \ref{CR} below).
\par

To study the theta correspondence, 
we often need to know a relation between theta lifts and induced representations.
To show such a relation, in non-archimedean case, Kudla's filtration \cite{Ku1} is useful.
This is a finite explicit filtration of the Jacquet module of the Weil representation.
Its archimedean analogue is the induction principle (see e.g., \cite[Theorem 4.5.5]{P1}).
As a matter of fact, however, the induction principle is 
just an analogue of a corollary of Kudla's filtration, 
so that it is less useful than Kudla's filtration itself.
\par

Both of the proofs of Theorem \ref{AGmain} and our main result (Theorem \ref{main} below)
use the local Gan--Gross--Prasad conjecture (for tempered representations).
In the non-archimedean case, Beuzart-Plessis proved it completely, 
whereas, 
in the archimedean case, he proved only a weaker version of it (\cite{BP}).
For discrete series representations of $\U(p,q)$, He \cite{He} proved 
the local Gan--Gross--Prasad conjecture in terms of Harish--Chandra parameters.
\par

Also, both of the proofs of two results (Theorem \ref{AGmain} and Theorem \ref{main})
use explicit descriptions of theta correspondence 
in the equal rank case and the almost equal rank case.
In the archimedean case, they are Paul's results (\cite{P1}, \cite{P3}).
In the non-archimedean case, they are two of Prasad's conjectures (\cite{Pr}), 
both of which are proven by Gan--Ichino \cite{GI}.

\subsection*{Organization of this paper}
This paper is organized as follows. 
In \S \ref{parameter}, we explain two parametrizations of 
irreducible representations of $\U(p,q)$, 
Harish-Chandra parameters and $L$-parameters, 
and we compare them.
Also we state the local Gan--Gross--Prasad conjecture.
In \S \ref{theta}, we recall some basic results of theta correspondence, 
including the Howe duality correspondence, the induction principle, 
seesaw identities, and Paul's results (\cite{P1}, \cite{P3}).
For the unitary dual pairs, we use Kudla's splitting, 
whereas, Paul uses double covers of unitary groups.
We also compare them in \S \ref{theta}.
In \S \ref{result}, we state the main result (Theorem \ref{main}) and 
its corollary (Corollary \ref{ind}).
Finally, we prove the main result in \S \ref{proof}.
The relation between Harish-Chandra parameters and $L$-parameters
(Theorem \ref{LLC} (4)) might be well-known, but 
there seems to be no proper reference.
For the convenience of the reader, we explain this relation in Appendix \ref{sec.explicit}.

\subsection*{Acknowledgments}
The author is grateful to Atsushi Ichino, Shunsuke Yamana, 
Yoshiki Oshima and Kohei Yahiro for useful discussions. 
We are most thankful to Hisayosi Matumoto for pointing out using Schmid's character identity
to prove Proposition \ref{generic}.
This is a key proposition to relate $L$-parameters with Harish-Chandra parameters.
This work was supported by the Foundation for Research Fellowships of Japan Society 
for the Promotion of Science for Young Scientists (DC1) Grant 26-1322.

\subsection*{Notation}
The symmetric group on $n$ letters is denoted by $S_n$.
For non-negative integers $p$ and $q$, 
we set $n = p+q$ and 
define the unitary group $\U(p,q)$ of signature $(p,q)$ by
\[
\U(p,q) =\left\{
g \in \GL_{n}(\C)\ |\ 
{}^t\overline{g}
\begin{pmatrix}
\1_p & 0 \\ 0 & -\1_q
\end{pmatrix}
g
=
\begin{pmatrix}
\1_p & 0 \\ 0 & -\1_q
\end{pmatrix}
\right\}.
\]
We also denote by $\U_n(\R)$ a unitary group of size $n$, i.e.,
$\U_n(\R)$ isomorphic to $\U(p,q)$ for some non-negative integers $p,q$ such that $n = p+q$.
\par

For a reductive Lie group $G$, we denote by $\Irr_\disc(G)$ (\resp $\Irr_\temp(G)$) 
the set of equivalence classes of discrete series representations 
(\resp tempered representations) of $G$.
In this paper, we only consider unitary representations of $G$, 
and identify them with the $(\g_\C, K)$-modules associated them, 
where $\g_\C = \Lie(G) \otimes_\R \C$ is the complexification of Lie algebra of $G$ 
and $K$ is a fixed maximal compact subgroup of $G$.
\par

Put $\C^1 = \{z \in \C^\times\ |\ z\overline{z} = 1\}$.
For a symplectic space $\W$ over $\R$, 
we denote the $\C^1$-cover of $\Sp(\W)$ by $\Mp(\W)$, 
and the double cover of $\Sp(\W)$ by $\cl{\Sp}(\W)$, which is a closed subgroup of $\Mp(\W)$.
Namely, 
\[
\begin{CD}
1 @>>> \C^1 @>>> \Mp(\W) @>>> \Sp(\W) @>>>1\\
@. @AAA @AAA @| @.\\
1 @>>> \{\pm1\} @>>> \cl{\Sp}(\W) @>>> \Sp(\W) @>>>1.
\end{CD}
\]
One should not confuse $\Mp(\W)$ with $\cl{\Sp}(\W)$.
When we consider representations of covering groups, 
we always assume that they are unitary and genuine.
\par

We denote representations of Lie groups by $\pi$ and $\sigma$, 
and ones of double covers by $\cl\pi$ and $\cl\sigma$.
Also we denote the contragredient representations of $\pi$, $\sigma$, $\cl\pi$ and $\cl\sigma$
by $\pi^\vee$, $\sigma^\vee$, $\cl\pi^\vee$ and $\cl\sigma^\vee$, respectively.
One should not confuse $\cl\pi$ with $\pi^\vee$.

\section{Classifications of irreducible representations of $\U(p,q)$}\label{parameter}
In this section, we give two parametrizations of irreducible representations of $\U(p,q)$, 
and we compare them.
We also state the local Gan--Gross--Prasad conjecture in \S \ref{sec.GGP}.

\subsection{Harish-Chandra parameters}
The Harish-Chandra parameters classify irreducible discrete series representations.
Let $G=\U(p,q)$.
We set $K \cong \U(p) \times \U(q)$ to be the maximal compact subgroup of $G$ consisting of 
the usual block diagonal matrices, 
and $T$ to be the maximal compact torus of $G$ consisting of diagonal matrices.
We denote the Lie algebras of $G$, $K$ and $T$ by
$\g$, $\kk$ and $\tt$, 
and its complexifications by $\g_\C$, $\kk_\C$ and $\tt_\C$, respectively. 
The set $\Delta_c$ of compact roots of $\g_\C$ with respect $\tt_\C$
and the set $\Delta_n$ of non-compact roots are given by
\begin{align*}
\Delta_c &= \{e_i-e_j\ |\ 1 \leq i,j \leq p\} \cup \{f_i-f_j \ |\ 1 \leq i,j \leq q\},\\
\Delta_n &= \{\pm(e_i-f_j)\ |\ 1 \leq i \leq p,\ 1 \leq j \leq q\}, 
\end{align*}
respectively.
Here, $e_i, f_j \in \tt_\C^*$ are defined by
\[
e_i \colon 
\begin{pmatrix}
t_1&&\\
&\ddots&\\
&&t_n
\end{pmatrix}
\mapsto t_i,\quad
f_j \colon 
\begin{pmatrix}
t_1&&\\
&\ddots&\\
&&t_n
\end{pmatrix}
\mapsto t_{p+j}.
\]
Note that $e_i$ and $f_j$ belong to $\I\tt^*$, i.e., 
the images of $e_i(\tt)$ and $f_j(\tt)$ are in $\I\R$.
\par

The Harish-Chandra parameter $\HC(\pi)$ of 
a discrete series representation $\pi$ of $\U(p,q)$ is of the form
\[
\HC(\pi) = (\lam_1, \dots, \lam_p; \lam_1', \dots, \lam_q') \in \I \tt^*, 
\]
where 
\begin{itemize}
\item
$\lam_i, \lam_j' \in \Z+\half{n-1}$;
\item
$\lam_i \not=\lam_j'$ for $1 \leq i \leq p$ and $1 \leq j \leq q$;
\item
$\lam_1 > \dots > \lam_p$ and $\lam_1' > \dots > \lam_q'$.
\end{itemize}
Here, using the basis $\{e_1 , \dots, e_p, f_1, \dots, f_q\}$, 
we identify $\I\tt^*$ with $\R^p \times \R^q$.
Via this identification, we regard $\HC(\pi)$ as an element of 
$(\Z+\half{n-1})^p \times (\Z+\half{n-1})^q$.
Hence we obtain an injection
\[
\HC \colon \Irr_\disc(\U(p,q)) \hookrightarrow 
\left(\Z+\half{n-1}\right)^p \times \left(\Z+\half{n-1}\right)^q.
\]
\par

The infinitesimal character $\tau_\lam$ of $\pi$ is the $W_G$-orbit of $\lam = \HC(\pi)$, 
where $W_G$ is the Weyl group of $G$ relative to $T$.
As well as $\lam = \HC(\pi)$ is regarded as an element of 
$(\Z+\half{n-1})^p \times (\Z+\half{n-1})^q$, 
we regard $\tau_\lam$ as an element of $(\Z+\half{n-1})^n/S_n$.
Note that given $\tau \in (\Z+\half{n-1})^n/S_n$, 
there are exactly $n!/(p!q!)$ discrete series representations of $\U(p,q)$ 
whose infinitesimal characters are equal to $\tau$.

\subsection{$L$-parameters}
The local Langlands correspondence is a parametrization of
irreducible tempered representations of $\U(p,q)$ in terms of $L$-parameters.
\par

For $\alpha \in \half{1}\Z$, we define a unitary character $\chi_{2\alpha}$ of $\C^\times$ by
\[
\chi_{2\alpha}(z) = \overline{z}^{-2\alpha}(z\overline{z})^{\alpha} = (z/\overline{z})^{\alpha}.
\]
Note that $\chi_{2\alpha}(\overline{z}) = \chi_{2\alpha}(z)^{-1} = \chi_{-2\alpha}(z)$.
When $a > 0$ and $\theta \in \R/2\pi\Z$, we have
\[
\chi_{2\alpha}(ae^{\I\theta}) = e^{2\alpha\I\theta}.
\]
Define $\Phi_\disc(\U_n(\R))$ by 
\[
\Phi_\disc(\U_n(\R)) = \left\{ \chi_{2\alpha_1} \oplus \dots \oplus \chi_{2\alpha_{n}}\ |\ 
\alpha_i \in \half{1}\Z, \ 2\alpha_i \equiv n-1,\ \alpha_1 > \dots > \alpha_n \right\}.
\] 
For $\phi = \chi_{2\alpha_1} \oplus \dots \oplus \chi_{2\alpha_{n}} \in \Phi_\disc(\U_n(\R))$, 
we define a component group $A_\phi$ of $\phi$ by
\[
A_\phi = (\Z/2\Z)e_{2\alpha_1} \oplus \dots \oplus (\Z/2\Z)e_{2\alpha_n}.
\]
Namely, $A_\phi$ is a free $(\Z/2\Z)$-module of rank $n$ 
equipped with a canonical basis $\{e_{2\alpha_1}, \dots, e_{2\alpha_n}\}$
associated to $\{\chi_{2\alpha_1}, \dots, \chi_{2\alpha_n}\}$.
\par

More generally, we define $\Phi_\temp(\U_n(\R))$ 
by the set of representations $\phi$ of $\C^\times$ of the form
\[
\phi = (m_1\chi_{2\alpha_1} \oplus \dots \oplus m_u\chi_{2\alpha_u}) \oplus 
(\xi_1 \oplus \dots \oplus \xi_v) \oplus ({}^c\xi_1^{-1} \oplus \dots \oplus {}^c\xi_v^{-1}), 
\]
where
\begin{itemize}
\item
$\alpha_i \in \half{1}\Z$ satisfies $2\alpha_i \equiv n-1$ and $\alpha_1 > \dots > \alpha_u$;
\item
$m_i \geq 1$ is the multiplicity of $\chi_{2\alpha_i}$ in $\phi$;
\item
$m_1+\dots+m_u+2v=n$;
\item
$\xi_i$ is a unitary character of $\C^\times$, 
which is not of the form $\chi_{2\alpha}$ with $2\alpha \equiv n-1$;
\item
${}^c\xi_i^{-1}$ is the unitary character of $\C^\times$ 
defined by ${}^c\xi_i^{-1}(z) = \xi_i(\overline{z}^{-1})$.
\end{itemize}
For such $\phi$, we define a component group $A_\phi$ of $\phi$ by
\[
A_\phi = (\Z/2\Z)e_{2\alpha_1} \oplus \dots \oplus (\Z/2\Z)e_{2\alpha_u}.
\]
We denote the Pontryagin dual of $A_\phi$ by $\widehat{A_\phi}$.
For $\eta \in \widehat{A_\phi}$, define $-\eta \in \widehat{A_\phi}$ by 
$(-\eta)(e_{2\alpha_i}) = -\eta(e_{2\alpha_i})$ for $i=1, \dots, u$.
\par

We define an additive character $\psi_{-2}^\C$ of $\C$ by
\[
\psi_{-2}^\C(z) = \exp(2\pi (\overline{z} - z))
\]
for $z \in \C$.
For a (continuous, completely reducible) representation $\phi$ of $\C^\times$, 
let $\ep(s, \phi, \psi_{-2}^\C)$ be the $\ep$-factor of $\phi$ (see \cite{T}).
It satisfies that
\begin{itemize}
\item
$\ep(s, \phi_1 \oplus \phi_2, \psi_{-2}^\C) 
= \ep(s, \phi_1, \psi_{-2}^\C) \cdot \ep(s, \phi_2, \psi_{-2}^\C)$;
\item
$\ep(1/2, \xi \oplus {}^c\xi^{-1}, \psi_{-2}^\C) = 1$ for any character $\xi$ of $\C^\times$;
\item
$\ep(1/2, \chi_{2\alpha}, \psi_{-2}^\C) = 1$ for $\alpha \in \Z$;
\item
When $\alpha \in \half{1}\Z \setminus \Z$, 
\[
\ep\left(\half{1}, \chi_{2\alpha}, \psi_{-2}^\C\right) = 
\left\{
\begin{aligned}
&-1 \iif \alpha > 0,\\
&+1 \iif \alpha < 0.
\end{aligned}
\right.
\]
\end{itemize}
For the last equation, see e.g., \cite[Proposition 2.1]{GGP2}.
We call the central value $\ep(1/2, \phi, \psi^\C_{-2})$ 
the root number of $\phi$ with respect to $\psi^\C_{-2}$, 
and we denote it simply by $\ep(\phi, \psi^\C_{-2})$.
\par

The local Langlands correspondence is a parametrization of $\Irr_\temp(\U(p,q))$ as follows:
\begin{thm}[\cite{L}, \cite{V3}, \cite{S1}, \cite{S2}, \cite{S3}]\label{LLC}
\begin{enumerate}
\item
There is a canonical surjection
\[
\bigsqcup_{p+q = n}\Irr_\temp(\U(p,q)) \rightarrow \Phi_\temp(\U_{n}(\R)).
\]
For $\phi \in \Phi_\temp(\U_{n}(\R))$, 
we denote by $\Pi_\phi$ the inverse image of $\phi$ under this map, 
and call $\Pi_\phi$ the $L$-packet associated to $\phi$.

\item
There is a bijection 
\[
J \colon \Pi_\phi \rightarrow \widehat{A_\phi}. 
\]
If $\pi \in \Pi_\phi$ corresponds to $\eta \in \widehat{A_\phi}$ via this bijection, 
we write $\pi = \pi(\phi, \eta)$ and call $(\phi,\eta)$ the $L$-parameter of $\pi$.

\item
If $\phi = \chi_{2\alpha_1} \oplus \dots \oplus \chi_{2\alpha_{n}} \in \Phi_\disc(\U_n(\R))$, 
the $L$-packet $\Pi_\phi$ consists of discrete series representations of various $\U(p,q)$
whose infinitesimal characters are $(\alpha_1, \dots, \alpha_n) \in (\Z + \half{n-1})^n/S_n$.

\item
If $\phi = \chi_{2\alpha_1} \oplus \dots \oplus \chi_{2\alpha_{n}} \in \Phi_\disc(\U_n(\R))$
with $\alpha_1 > \dots > \alpha_n$, 
the Harish-Chandra parameter 
\[
\HC(\pi(\phi, \eta)) = (\lam_1, \dots, \lam_p; \lam'_1, \dots, \lam_q')
\] 
of $\pi(\phi,\eta) \in \Pi_\phi$ is given so that
\begin{itemize}
\item
$\{\lam_1, \dots, \lam_p, \lam_1', \dots, \lam_q'\} = \{\alpha_1, \dots, \alpha_n\}$;
\item
$\alpha_i \in \{\lam_1, \dots, \lam_p\}$ if and only if $\eta(e_{2\alpha_i}) = (-1)^{i-1}$.
\end{itemize}
In particular, $\pi(\phi, \eta) \in \Irr_\disc(\U(p,q))$ with
\[
p = \#\{i\ |\ \eta(e_{2\alpha_i}) = (-1)^{i-1}\},
\quad
q = \#\{i\ |\ \eta(e_{2\alpha_i}) = (-1)^{i}\}.
\]

\item
If $\phi = \xi \oplus \phi_0 \oplus {}^c\xi^{-1}$ with a unitary character $\xi$ of $\C^\times$ 
and an element $\phi_0$ in $\Phi_\temp(\U_{n-2}(\R))$, 
for any $\pi(\phi_0,\eta_0) \in \Pi_{\phi_0} \cap \Irr_\temp(\U(p-1,q-1))$, 
the induced representation $\Ind_P^{\U(p,q)}(\xi \otimes \pi(\phi_0,\eta_0))$ 
decomposes as follows:
\[
\Ind_P^{\U(p,q)}(\xi \otimes \pi(\phi_0, \eta_0)) = 
\bigoplus_{\substack{\eta \in \widehat{A_\phi}, \\ \eta|A_{\phi_0} = \eta_0}}\pi(\phi,\eta).
\]
Here, $P$ is a parabolic subgroup of $\U(p,q)$ 
with Levi subgroup $M_P = \C^\times \times \U(p-1,q-1)$.

\item
The contragredient representation of $\pi(\phi,\eta)$ is given by $\pi(\phi^\vee, \eta^\vee)$, 
where $\eta^\vee \colon A_{\phi^\vee} \rightarrow \{\pm1\}$ is defined by 
\[
\eta^\vee(e_{-2\alpha_i}) = \left\{
\begin{aligned}
&\eta(e_{2\alpha_i}) \iif \text{$n$ is odd},\\
&-\eta(e_{2\alpha_i}) \iif \text{$n$ is even}
\end{aligned}
\right.
\]
for any $e_{-2\alpha_i} \in A_{\phi^\vee}$.

\item
If $\pi = \pi(\phi,\eta) \in \Irr_\temp(\U(p,q))$, then $\pi(\phi, -\eta) \in \Irr_\temp(\U(q,p))$
is the same representation as $\pi$ via the canonical identification $\U(p,q) = \U(q,p)$
as subgroups of $\GL_{n}(\C)$.

\item
If $\phi \in \Phi_\temp(\U_n(\R))$, then $\phi \chi_2 = \phi \otimes \chi_2 \in \Phi_\temp(\U_n(\R))$
and there is a canonical identification $A_\phi = A_{\phi\chi_2}$.
If $\pi = \pi(\phi, \eta)$, the corresponding representation $\pi(\phi \chi_2, \eta)$
is the determinant twist $\pi \otimes \det$.
\end{enumerate}
\end{thm}

In fact, the bijection $J \colon \Pi_\phi \rightarrow \widehat{A_\phi}$ in Theorem \ref{LLC} (2) 
is characterized by endoscopic character identities as in \cite{S3}.
Also $J$ depends on a choice of a pair of a quasi-split $\U(p,q)$ (i.e., $|p-q| \leq 1$) 
and a Whittaker datum of $\U(p,q)$.
There are exactly two such pairs.
Through this paper, we take a specific choice of a pair.
(see Appendix \ref{sec.explicit}).
Theorem \ref{LLC} (4) highly depends on this choice, 
and there seems to be no proper reference.
We will discuss this part in Appendix \ref{sec.explicit} below.
\par

In Theorem \ref{LLC} (4), we see that
\[
(-1)^{i-1} = \ep(\phi \otimes \chi_{-2\alpha_i} \otimes \chi_{-1}, \psi^\C_{-2}).
\]
Hence the Harish-Chandra parameter of $\pi$ and the unitary group $\U(p,q)$ which acts on $\pi$
can be determined by the $L$-parameter $\lam = (\phi,\eta)$ of $\pi$ and certain root numbers.

\subsection{Local Gan--Gross--Prasad conjecture}\label{sec.GGP}
To prove the main result, we use the local Gan--Gross--Prasad conjecture \cite{GGP1}, 
which gives an answer to restriction problems.
\par

Suppose that $(\U_n(\R), \U_{n+1}(\R)) = (\U(p,q), \U(p+1,q))$ or $(\U(p,q), \U(p,q+1))$.
Then there is a canonical injection $\U_n(\R) \hookrightarrow \U_{n+1}(\R)$, so that
we have a diagonal map
\[
\Delta \colon \U_n(\R) \rightarrow \U_n(\R) \times \U_{n+1}(\R).
\]
By a result of Sun--Zhu \cite{SZ1}, 
for $\pi_n \in \Irr_\temp(\U_{n}(\R))$ and $\pi_{n+1} \in \Irr_\temp(\U_{n+1}(\R))$, 
we have
\[
\dim_\C\Hom_{\Delta\U_{n}(\R)}(\pi_n \otimes \pi_{n+1}, \C) \leq 1.
\]
We call a pair $(\U_n(\R), \U_{n+1}(\R))$ is relevant if
\[
(\U_n(\R), \U_{n+1}(\R)) = 
\left\{
\begin{aligned}
&(\U(p,q), \U(p+1,q)) \iif \text{$n$ is even},\\
&(\U(p,q), \U(p,q+1)) \iif \text{$n$ is odd}
\end{aligned}
\right.
\]
for some $(p, q)$ such that $p+q = n$.
Note that if $(\U(p,q), \U(p',q'))$ with $p'+q' = p+q+1$ is not relevant, 
then $(\U(q, p), \U(q',p'))$ is relevant.
Gan, Gross and Prasad predicted when the $\Hom$ space 
\[
\Hom_{\Delta\U_{n}(\R)}(\pi_n \otimes \pi_{n+1}, \C) 
\]
is nonzero for relevant pairs.

\begin{conj}[local Gan--Gross--Prasad conjecture (GGP) {\cite[Conjecture 17.3]{GGP1}}]
\label{GGP}
Let $\phi_n \in \Phi_\temp(\U_{n}(\R))$ and $\phi_{n+1} \in \Phi_\temp(\U_{n+1}(\R))$
such that
\begin{align*}
\phi_{n} &= (m_1\chi_{2\alpha_1} \oplus \dots \oplus m_{u}\chi_{2\alpha_{u}}) 
\oplus (\xi_1 \oplus \dots \oplus \xi_v) \oplus ({}^c\xi_1^{-1} \oplus \dots \oplus {}^c\xi_v^{-1}),\\
\phi_{n+1} &= (m'_1\chi_{2\beta_1} \oplus \dots \oplus m'_{u'}\chi_{2\beta_{u'}}) 
\oplus (\xi'_1 \oplus \dots \oplus \xi'_{v'}) 
\oplus ({}^c{\xi'_1}^{-1} \oplus \dots \oplus {}^c{\xi'_{v'}}^{-1}),
\end{align*}
where
\begin{itemize}
\item
$\alpha_i, \beta_j \in \half{1}\Z$ such that 
$2\alpha_i \equiv n-1 \bmod 2$ and $2\beta_j \equiv n \bmod 2$;
\item
$m_i \geq 1$ (\resp $m'_j \geq 1$) is the multiplicity of $\chi_{2\alpha_i}$ (\resp $\chi_{2\beta_j}$)
in $\phi_{n}$ (\resp $\phi_{n+1}$);
\item
$m_1 + \dots + m_u + 2v = n$ and $m'_1 + \dots + m'_{u'} + 2v' = n+1$; 
\item
$\xi_i$ (\resp $\xi'_j$) is a unitary character of $\C^\times$, 
which is not of the form $\chi_{2\alpha}$ with $2\alpha \equiv n-1 \bmod 2$ 
(\resp $\chi_{2\beta}$ with $2\beta \equiv n \bmod 2$).
\end{itemize}
Then there exists a unique pair of representations 
$(\pi_n,\pi_{n+1}) \in \Pi_{\phi_n} \times \Pi_{\phi_{n+1}}$ such that
\begin{itemize}
\item
$(\pi_n,\pi_{n+1})$ is a pair of representations of a relevant pair $(\U_n(\R), \U_{n+1}(\R))$;
\item
$\Hom_{\Delta\U_{n}(\R)}(\pi_n \otimes \pi_{n+1}, \C) \not= 0$.
\end{itemize}
Moreover,
\begin{align*}
J(\pi_{n})(e_{2\alpha_i}) &= \ep(\chi_{2\alpha_i} \otimes \phi_{n+1}, \psi_{-2}^{\C}),\\
J(\pi_{n+1})(e_{2\beta_j}) &= \ep(\phi_{n} \otimes \chi_{2\beta_j}, \psi_{-2}^\C) 
\end{align*}
for $e_{2\alpha_i} \in A_{\phi_{n}}$ and $e_{2\beta_j} \in A_{\phi_{n+1}}$.
\end{conj}

When $\phi_n \in \Phi_\disc(\U_{n}(\R))$ and $\phi_{n+1} \in \Phi_\disc(\U_{n+1}(\R))$, 
Conjecture \ref{GGP} is proven by He \cite{He}.
In general, Beuzart-Plessis \cite{BP} showed a weaker version of Conjecture \ref{GGP}.
\par

We use this conjecture as the following form.
\begin{prop}\label{GGPprop}
Assume the local Gan--Gross--Prasad conjecture (Conjecture \ref{GGP}).
Let $\phi \in \Phi_\temp(\U_n(\R))$ and $\phi' \in \Phi_\temp(\U_{n+1}(\R))$ such that
\begin{align*}
\phi &= \chi_{2\alpha_1} \oplus \dots \oplus \chi_{2\alpha_u} 
\oplus (\xi_1 \oplus \dots \oplus \xi_v) \oplus ({}^c\xi_1^{-1} \oplus \dots \oplus {}^c\xi_v^{-1}),\\
\phi' &= \chi_{2\beta_1} \oplus \dots \oplus \chi_{2\beta_{u'}} 
\oplus (\xi'_1 \oplus \dots \oplus \xi'_{v'}) 
\oplus ({}^c{\xi'_1}^{-1} \oplus \dots \oplus {}^c{\xi'_{v'}}^{-1}),
\end{align*}
where
\begin{itemize}
\item
$\alpha_i, \beta_j \in \half{1}\Z$ such that 
$2\alpha_i \equiv n-1 \bmod 2$ and $2\beta_j \equiv n \bmod 2$;
\item
$\xi_i$ (\resp $\xi'_j$) is a unitary character of $\C^\times$
(which can be of the form $\chi_{2\alpha}$ (\resp $\chi_{2\beta}$));
\item
$u+2v=n$ (\resp $u'+2v'=n+1$).
\end{itemize}
Then for $(\pi,\pi') \in \Pi_\phi \times \Pi_{\phi'}$, 
the following are equivalent:
\begin{itemize}
\item
$(\pi, \pi') \in \Irr_\temp(\U(p,q)) \times \Irr_\temp(\U(p+1,q))$ for some $(p,q)$ and
$\Hom_{\U(p,q)}(\pi', \pi) \not= 0$;
\item
$J(\pi) \in \widehat{A_\phi}$ and $J(\pi') \in \widehat{A_{\phi'}}$ are given by
\begin{align*}
J(\pi)(e_{2\alpha}) &= (-1)^{\#\{j \in \{1, \dots, u'\}\ |\ \beta_j<\alpha\} + n},\\
J(\pi')(e_{2\beta}) &= (-1)^{\#\{i \in \{1, \dots, u\}\ |\ \alpha_i<\beta\} + n}
\end{align*}
for any $\chi_{2\alpha} \subset \phi$ so that $e_{2\alpha} \in A_\phi$ and
any $\chi_{2\beta} \subset \phi'$ so that $e_{2\beta} \in A_{\phi'}$. 
\end{itemize}
\end{prop}
\begin{proof}
Note that $\Hom_{\U(p,q)}(\pi', \pi) \not= 0$ 
if and only if $\Hom_{\Delta\U(p,q)}(\pi' \otimes \pi^\vee, \C) \not= 0$.
First, we assume that $n=p+q$ is even.
Then $(\U(p,q), \U(p+1,q))$ is a relevant pair.
By the local Gan--Gross--Prasad conjecture (Conjecture \ref{GGP}), 
we should see that $\Hom_{\Delta\U(p,q)}(\pi' \otimes \pi^\vee, \C) \not= 0$
if and only if
\begin{align*}
J(\pi^\vee)(e_{-2\alpha}) 
&= \ep(\chi_{-2\alpha} \otimes \phi', \psi^\C_{-2})
=(-1)^{\#\{j \in \{1, \dots, u'\}\ |\ -\alpha+\beta_j > 0\}},\\
J(\pi')(e_{2\beta}) &= \ep(\phi^\vee \otimes \chi_{2\beta}, \psi^\C_{-2})
=(-1)^{\#\{i \in \{1, \dots, u\}\ |\ -\alpha_i+\beta > 0\}}
\end{align*}
for any $e_{-2\alpha} \in A_{\phi^\vee}$ and any $e_{2\beta} \in A_{\phi'}$.
By Theorem \ref{LLC} (6), we see that
\[
J(\pi)(e_{2\alpha}) = -J(\pi^\vee)(e_{-2\alpha}) 
=-(-1)^{\#\{j \in \{1, \dots, u'\}\ |\ \beta_j > \alpha\}}
=(-1)^{\#\{j \in \{1, \dots, u'\}\ |\ \beta_j < \alpha\}}
\]
for any $e_{2\alpha} \in A_{\phi}$
since $u' \equiv n+1 \equiv 1 \bmod 2$.
Hence we obtain the assertion in the case when $n=p+q$ is even.
\par

Next, we assume that $n=p+q$ is odd.
Then $(\U(q,p), \U(q,p+1))$ is a relevant pair.
By Theorem \ref{LLC} (7), we see that
$\Hom_{\U(p,q)}(\pi', \pi) \not= 0$ 
if and only if
$\Hom_{\U(q,p)}(\pi(\phi', -J(\pi')), \pi(\phi, -J(\pi))) \not= 0$. 
By a similar argument to the first case, this is equivalent that
\begin{align*}
-J(\pi)(e_{2\alpha}) &= (-1)^{\#\{j \in \{1, \dots, u'\}\ |\ \beta_j<\alpha\}},\\
-J(\pi')(e_{2\beta}) &= (-1)^{\#\{i \in \{1, \dots, u\}\ |\ \alpha_i<\beta\}}
\end{align*}
for any $e_{2\alpha} \in A_{\phi}$ and any $e_{2\beta} \in A_{\phi'}$.
Hence we obtain the assertion in the case when $n=p+q$ is odd.
This completes the proof.
\end{proof}

\section{Theta liftings}\label{theta}
In this subsection, we review the theory of theta liftings.
First, we recall Kudla's splitting \cite{Ku2} of a unitary dual pair (\S \ref{sec. splitting}).
Then we can consider theta lifts of irreducible unitary representations of unitary groups as in \S \ref{arch}.
On the other hand, in various results, including Paul's ones (\cite{P1}, \cite{P2}, \cite{P3}),
irreducible genuine representations of certain double covers of unitary groups are used for theta lifts.
We compare Kudla's splitting with double covers of unitary groups in \S \ref{sec. double}.
In \S \ref{sec. property} and \S \ref{sec. paul}, we recall
basic properties on theta liftings and Paul's results \cite{P1}, \cite{P3}, respectively.

\subsection{Kudla's splitting}\label{sec. splitting}
Let $W = W_{p,q}$ (\resp $V = V_{r,s}$) 
be a (right) complex vector space of dimension $n=p+q$
(\resp $m=r+s$) equipped with a hermitian form $\pair{\cdot, \cdot}_W$ 
(\resp a skew-hermitian form $\pair{\cdot, \cdot}_V$)
of signature $(p,q)$ (\resp $(r,s)$). 
Namely, 
\begin{itemize}
\item
the pairings $\pair{\cdot, \cdot}_W$ and $\pair{\cdot, \cdot}_V$ satisfy
\begin{align*}
\pair{w_1a, w_2b}_W = a\overline{b}\pair{w_1,w_2}_W, 
&\quad \pair{w_2,w_1}_W = \overline{\pair{w_1,w_2}_W},\\
\pair{v_1a, v_2b}_V = a\overline{b}\pair{v_1,v_2}_V, 
&\quad \pair{v_2,v_1}_V = -\overline{\pair{v_1,v_2}_V}
\end{align*}
for $a,b \in \C$, $w_1, w_2 \in W$ and $v_1,v_2 \in V$; 
\item
there exist $e_1, \dots, e_n \in W$ and $e'_1, \dots, e'_m \in V$ such that
\[
\pair{e_i,e_j}_W = \left\{
\begin{aligned}
&0	\iif i \not= j,\\
&1	\iif i = j \leq p,\\
&-1	\iif i = j > p,
\end{aligned}
\right.
\quad
\pair{e'_i,e'_j}_V = \left\{
\begin{aligned}
&0	\iif i \not= j,\\
&\I	\iif i = j \leq r,\\
&-\I	\iif i = j > r.
\end{aligned}
\right.
\]
\end{itemize}
The isometry group $\U(W)$ of $\pair{\cdot, \cdot}_W$ (\resp $\U(V)$ of $\pair{\cdot, \cdot}_V$), 
which has a left action on $W$ (\resp on $V$), 
is isomorphic to $\U(p,q)$ (\resp $\U(r,s)$).
Let $\W = V \otimes_\C W$ 
be the symplectic space over $\R$ of dimension $2mn$
equipped with the symplectic form
\[
\pair{v_1\otimes w_1, v_2\otimes w_2} = \tr_{\C/\R}(\pair{v_1,v_2}_V \cdot \pair{w_1,w_2}_W)
\]
for $v_1,v_2 \in V$ and $w_1,w_2 \in W$. 
The symplectic group $\Sp(\W)$ acts on $\W$ on the left.
\par

We note that the convention in \cite{Ku2} and \cite{HKS} differs from ours.
They use the following:
\begin{itemize}
\item
$W'$ is a left vector space and
the hermitian form $\pair{\cdot, \cdot}_{W'}$ on $W'$ satisfies 
\[
\pair{aw'_1, bw'_2}_{W'} = a\pair{w_1,w_2}_{W'} \overline{b}
\]
for $a,b \in \C$ and $w'_1,w'_2 \in W'$; 
\item
$V'$ is a right vector space and
the skew-hermitian form $\pair{\cdot, \cdot}_{V'}$ on $V'$ satisfies
\[
\pair{v'_1a, v'_2b}_{V'} = \overline{a}\pair{v'_1,v'_2}_{V'}b
\]
for $a,b \in \C$ and $v'_1,v'_2 \in V'$; 
\item
the symplectic form on $\mathbb{W}' = V' \otimes_\C W'$ is defined by
\[
\pair{v'_1\otimes w'_1, v'_2\otimes w'_2}' = 
\half{1} \cdot \tr_{\C/\R}(\pair{v'_1,v'_2}_{V'} \cdot \overline{\pair{w'_1,w'_2}_{W'}})
\]
for $v'_1,v'_2 \in V'$ and $w'_1,w'_2 \in W'$;
\item
$\U(W')$, $\U(V')$ and $\Sp(\mathbb{W}')$ act on $W'$, $V'$ and $\mathbb{W}'$ 
on the right, left and right, respectively.
\end{itemize}
To use results in \cite{Ku2} and \cite{HKS}, 
we have to compare these conventions.
\par

First, we compare $\U(W)$ with $\U(W')$.
Assume that $W'$ has a basis $\{e_1, \dots, e_n\}$ satisfying the same conditions as above, 
so that $W=W'$.
However, since $W$ is a right $\C$-vector space, whereas, $W'$ is a left $\C$-vector space, 
we obtain expressions
\begin{align*}
W &= \left\{
\begin{pmatrix}
e_1 & \ldots & e_n
\end{pmatrix}
\begin{pmatrix}
a_1 \\ \vdots \\ a_n
\end{pmatrix}
\ |\ a_i \in \C
\right\} \cong \C^n \quad \text{(column vectors)},\\
W' &= \left\{
\begin{pmatrix}
a'_1 & \ldots & a'_n
\end{pmatrix}
\begin{pmatrix}
e_1 \\ \vdots \\ e_n
\end{pmatrix}
\ |\ a_i' \in \C
\right\} \cong \C^n \quad \text{(row vectors)}.
\end{align*}
Via the $\C$-linear isomorphism
\[
W \ni w = 
\begin{pmatrix}
e_1 & \ldots & e_n
\end{pmatrix}
\begin{pmatrix}
a_1 \\ \vdots \\ a_n
\end{pmatrix}
\mapsto 
w'=
\begin{pmatrix}
a_1 & \ldots & a_n
\end{pmatrix}
\begin{pmatrix}
e_1' \\ \vdots \\ e'_n
\end{pmatrix}
\in W', 
\]
we identify $W$ with $W'$.
Then we have $\pair{w_1',w_2'}_{W'} = \pair{w_1, w_2}_W$ for any $w_1, w_2 \in W$.
Also, via the identifications $W \cong \C^n$ (row vectors) and $W' \cong \C^n$ (column vectors), 
we have expressions
\begin{align*}
\U(W) &= \left\{
g \in \GL_n(\C) \ |\ 
{}^tg \begin{pmatrix}
\1_p & 0 \\
0 & -\1_q \\
\end{pmatrix} \overline{g}
= 
\begin{pmatrix}
\1_p & 0 \\
0 & -\1_q \\
\end{pmatrix}
\right\},\\
\U(W') &= \left\{
g' \in \GL_n(\C) \ |\ 
g' \begin{pmatrix}
\1_p & 0 \\
0 & -\1_q \\
\end{pmatrix} {}^t\overline{g'}
= 
\begin{pmatrix}
\1_p & 0 \\
0 & -\1_q \\
\end{pmatrix}
\right\}.
\end{align*}
The above identification $W = W'$ implies the map
\[
\U(W) \ni g \mapsto {}^tg \in \U(W').
\]
\par

Next, we compare $\U(V)$ with $\U(V')$.
Assume that $V'$ has a basis $\{e'_1, \dots, e'_m\}$ so that $V=V'$ as right $\C$-vector spaces.
Set the skew-hermitian pairing $\pair{\cdot, \cdot}_{V'}$ by
\[
\pair{v_1', v_2'}_{V'} = 2 \overline{\pair{v_1, v_2}_V}
\]
for elements $v_1=v_1'$ and $v_2=v_2'$ in $V=V'$.
Note that the signature of $V'$ is $(s,r)$.
We have expressions
\begin{align*}
\U(V) &= \left\{
h \in \GL_m(\C) \ |\ 
{}^th \begin{pmatrix}
\I\1_r & 0 \\
0 & -\I\1_s \\
\end{pmatrix} \overline{h}
= 
\begin{pmatrix}
\I\1_p & 0 \\
0 & -\I\1_q \\
\end{pmatrix}
\right\},\\
\U(V') &= \left\{
h' \in \GL_m(\C) \ |\ 
{}^t\overline{h'} \begin{pmatrix}
-2\I\1_r & 0 \\
0 & 2\I\1_s \\
\end{pmatrix} h'
= 
\begin{pmatrix}
-2\I\1_r & 0 \\
0 & 2\I\1_s \\
\end{pmatrix}
\right\}.
\end{align*}
Hence $\U(V')$ coincides with $\U(V)$ as subgroups of $\GL_m(\C)$.
\par

Since 
\begin{align*}
\pair{v'_1\otimes w'_1, v'_2\otimes w'_2}' 
&= 
\half{1} \cdot \tr_{\C/\R}(\pair{v'_1,v'_2}_{V'} \cdot \overline{\pair{w'_1,w'_2}_{W'}})\\
&=
\tr_{\C/\R}(\overline{\pair{v_1,v_2}_{V} \cdot \pair{w_1,w_2}_{W}})
= \pair{v_1\otimes w_1, v_2\otimes w_2} 
\end{align*}
for $v_1, v_2 \in V$ and $w_1, w_2 \in W$, 
we see that
the two symplectic forms on $\W = \W'$ agree.
\par

Note that both $\U(W)$ and $\Sp(\W)$ act on $W$ and $\W$ on the left, 
respectively, whereas, 
$\U(W')$ and $\Sp(\W')$ act on $W'$ and $\W'$ on the right, respectively.
Hence the canonical embedding 
$\alpha = \alpha_V \colon \U(W) \rightarrow \Sp(\W)$
coincides with the counterpart 
$\alpha' = \alpha_{V'} \colon \U(W') \rightarrow \Sp(\W')$
via the above identifications.
Therefore the results in \cite{Ku2} can be transferred to our convention.
\par

Fix a non-trivial additive character $\psi$ of $\R$.
Let $\Mp(\W)$ be the $\C^1$-cover of $\Sp(\W)$, i.e.,
\[
\begin{CD}
1 @>>> \C^1 @ >>> \Mp(\W) @>>> \Sp(\W) @>>>1.
\end{CD}
\]
Choosing a character $\chi = \chi_V$ of $\C^\times$ such that $\chi|\R^\times = \sgn^{m}$ 
with $m = r+s = \dim(V)$, 
Kudla gave a splitting 
\[
\xymatrix{
& \Mp(\W) \ar@{->}[d]\\
\U(W) \ar@{-->}^{\cl{\alpha}_\chi}[ur] \ar@{->}^{\alpha}[r] & \Sp(\W).
}
\]

\begin{lem}
We identify $\Mp(\mathbb{W}) = \Sp(\mathbb{W}) \times \C^1$ as sets as in \cite{Ku2}.
If we write $\cl\alpha_\chi(g) = (\alpha(g), \beta_\chi(g))$ for $g \in \U(W)$, 
then $\beta_\chi(g)$ satisfies that
\[
\beta_\chi(g)^8 = \chi(\det(g))^4
\]
\end{lem}
\begin{proof}
Let $W_-$ denote the space $W$ with the hermitian form $-\pair{\cdot, \cdot}_W$.
Consider the space $W \oplus W_-$.
Now we have a canonical embedding
\[
\Delta \colon \U(W) \times \U(W_-) \rightarrow \U(W \oplus W_-).
\]
Let $x(\Delta(g,1))$ be Rao's function (see \cite[\S 1]{Ku2}).
Note that $x(\Delta(g,1))$ is an element in $\C^\times/ \R_{>0}$ (\cite[Corollary 1.5]{Ku2}).
By \cite[Lemma 3.5]{Ku2}, it satisfies that
\[
\det(g)^2 = \left( x(\Delta(g,1)) \cdot {\overline{x(\Delta(g,1))}}^{-1} \right)^2.
\]
By \cite[Theorems 3.1, 3.3]{Ku2}, 
we have 
\[
\beta_\chi(g) = \chi(x(\Delta(g,1)))\zeta
\]
for some $8$-th root of unity $\zeta$.
Since $\chi^2|\R^\times = \1$, we have
\[
\beta_\chi(g)^8 = \chi(x(\Delta(g,1)))^8
= \chi\left( x(\Delta(g,1)) \cdot {\overline{x(\Delta(g,1))}}^{-1} \right)^4
= \chi(\det(g))^4.
\]
This completes the proof.
\end{proof}

\subsection{Double cover of $\U(p,q)$}\label{sec. double}
Let $W$, $V$ and $\W$ be as in \S \ref{sec. splitting}.
Then we have a canonical map $\alpha = \alpha_V \colon \U(W) \rightarrow \Sp(\W)$.
Let $\cl{\Sp}(\W)$ be the double cover of $\Sp(\W)$, 
which is a closed subgroup of $\Mp(\W)$, i.e., 
\[
\begin{CD}
1 @>>> \{\pm1\} @>>> \cl{\Sp}(\W) @>>> \Sp(\W) @>>> 1.
\end{CD}
\]
Fix $\nu \in \Z$ such that $\nu \equiv m = r+s \bmod 2$ and
define the ${\det}^{\nu/2}$-cover of $\U(W)$ by 
\[
\cl{\U}(W) = \{(g,z) \in \U(W) \times \C^1 \ |\ z^2 = \det(g)^{\nu}\}.
\]
It has a genuine character 
\[
{\det}^{\nu/2} \colon \cl{\U}(W) \rightarrow \C^\times,\ (g,z) \mapsto z.
\]
Hence the set of genuine tempered representations of $\cl\U(W)$ is give by
\[
\Irr_\temp(\cl\U(W)) = \{\pi \otimes {\det}^{-\nu/2}\ |\ \pi \in \Irr_\temp(\U(W))\}.
\]
\par

As in \cite[\S 1.2]{P1}, we have a homomorphism
\[
\cl{\alpha} = \cl{\alpha}_V \colon \cl{\U}(W) \rightarrow \cl{\Sp}(\W)
\]
such that $\cl{\alpha}(\cl{\U}(W))$ is the inverse image of $\alpha(\U(W))$,
and the diagram
\[
\begin{CD}
\cl{\U}(W) @>\cl{\alpha}>> \cl{\Sp}(\W)\\
@VVV @VVV \\
\U(W) @>\alpha>> \Sp(\W)
\end{CD}
\]
is commutative, where the left arrow is the first projection.
We write the composition of $\cl\alpha \colon \cl\U(W) \rightarrow \cl\Sp(\W)$
with the inclusion map $\cl\Sp(\W) \hookrightarrow \Mp(\W) = \Sp(\W) \times \C^1$ as
\[
(g, z) \mapsto (\alpha(g), \beta(g,z)).
\]
Note that $\beta(g,-z) = -\beta(g,z)$.
We put $\mu_8 = \{\zeta \in \C^\times \ |\ \zeta^8 = 1\}$
to be the set of $8$-th roots of unity in $\C^\times$.

\begin{lem}
For any $(g,z) \in \cl\U(W)$, 
the value $\beta(g,z)$ belongs to $\mu_8$.
\end{lem}
\begin{proof}
We write 
\[
\Mp(\W) = \Sp(\W) \times \C^1,
\quad
\cl{\Sp}(\W) = \Sp(\W) \times \{\pm1\}
\]
as sets.
The multiplication law of $\Mp(\W)$ 
is given by $(h_1, z_1) \cdot (h_2, z_2) = (h_1h_2, z_1z_2c(h_1,h_2))$, 
where $c(h_1, h_2)$ is Rao's $2$-cocycle (see \cite{Ra} and \cite{Ku2}).
By \cite[Theorem 4.1 (5)]{Ra} (cf.~\cite[Theorem]{Ku2}), 
$c(g_1, g_2)$ is a Weil index, which is an $8$-th root of unity.
We write the inclusion map $\cl\Sp(\W) \hookrightarrow \Mp(\W)$ as
\[
(h,\epsilon) \mapsto (h, \gamma(h,\epsilon)).
\]
When $(h_1,\epsilon_1) \cdot (h_2, \epsilon_2) = (h_1h_2, \epsilon_3)$ in $\cl{\Sp}(\W)$, 
we see that
\[
\gamma(h_1h_2, \epsilon_3) = \gamma(h_1,\epsilon_1)\gamma(h_2,\epsilon_2)c(h_1,h_2).
\]
Since $\gamma(h,-\epsilon) = -\gamma(h,\epsilon)$, 
we see that $(h, \epsilon) \mapsto \gamma(h,\epsilon)^8$ factors through 
a group homomorphism of $\Sp(\W)$.
Since $\Sp(\W)$ is semisimple, there exists a dense subset $X$ of $\cl{\Sp}(\W)$
such that $\gamma(h,\epsilon) \in \mu_8$ for any $(h,\epsilon) \in X$.
Since $\cl{\Sp}(\W)$ is closed in $\Mp(\W)$, 
the closure of $X$ in $\Mp(\W)$, 
which is contained in $\Sp(\W) \times \mu_8 \subset \Mp(\W)$,  
coincides with $\cl{\Sp}(\W)$.
This means that
$\gamma(h,\epsilon) \in \mu_8$ for any $(h,\epsilon) \in \cl\Sp(\W)$.
In particular, since the map $(g, z) \mapsto (\alpha(g), \beta(g,z))$ factors through 
$\cl\Sp(\W) \hookrightarrow \Mp(\W)$, 
we conclude that $\beta(g,z)$ is an $8$-th root of unity for any $(g,z) \in \cl\U(W)$.
\end{proof}

Now we compare $\beta_\chi(g)$ with $\beta(g,z)$.
\begin{prop}\label{beta}
Fix $\nu \in \Z$ such that $\nu \equiv m \bmod 2$.
Let $\chi = \chi_\nu$ be the character of $\C^\times$ given by
\[
\chi(ae^{\I\theta}) = e^{\nu\I\theta}
\]
for $a > 0$ and $\theta \in \R/2\pi\Z$.
Suppose that $\cl\U(W)$ is the ${\det}^{\nu/2}$-cover of $\U(W)$, 
i.e., $z^2 = \det(g)^\nu$ when $(g,z) \in \cl\U(W)$.
Then
\[
\beta_\chi(g) = \beta(g,z)z
\]
for any $(g,z) \in \cl\U(W)$.
\end{prop}
\begin{proof}
Since the maps $\U(W) \rightarrow \Mp(\W),\ g \mapsto (\alpha(g),\beta_\chi(g))$
and $\cl\U(W) \rightarrow \Mp(\W),\ (g,z) \mapsto (\alpha(g),\beta(g,z))$
are homomorphism, we see that
\[
\beta_\chi(g_1g_2)\beta_\chi(g_1)^{-1}\beta_\chi(g_2)^{-1}
= c(\alpha(g_1),\alpha(g_2))
= \beta(g_1g_2, z_1z_2)\beta(g_1, z_1)^{-1}\beta(g_2, z_2)^{-1}
\]
for any $(g_1,z_1), (g_2,z_2) \in \cl\U(W)$.
This implies that the map
\[
\eta \colon \cl\U(W) \ni (g,z) \mapsto \frac{\beta_\chi(g)}{\beta(g,z)z} \in \C^1
\]
factors through a group homomorphism $\eta$ of $\U(W)$.
Since
\[
\left(\frac{\beta_\chi(g)}{\beta(g,z)z}\right)^8
=\frac{\chi(\det(g))^4}{\det(g)^{4\nu}}=1
\]
for any $(g,z) \in \cl\U(W)$, we may regard this map as a group homomorphism
\[
\eta \colon \U(W) \rightarrow \mu_8,\ 
g \mapsto \frac{\beta_\chi(g)}{\beta(g,z)z}.
\]
Since $\U(W)$ is a Lie group, 
we can find an open neighborhood $X$ of $1$ such that 
for any $g \in X$, there exists $h \in \U(W)$ such that $g=h^8$.
This means that $X$ is contained in the kernel of $\eta$, 
so that $\eta$ is continuous.
In conclusion, $\eta$ is a (continuous) character of $\U(W)$ of finite order.
Since $\U(W)$ is connected, it must be the trivial character.
Hence $\beta_\chi(g) = \beta(g,z)z$.
\end{proof}

\subsection{Basic properties of theta liftings}\label{sec. property}
Through this paper, 
for each hermitian space $W=W_{p,q}$ and each skew-hermitian space $V=V_{r,s}$
as in \S \ref{sec. splitting}, 
we fix characters $\chi_{W} = \chi_{W_{p,q}}$ and $\chi_V = \chi_{V_{r,s}}$
such that 
\begin{itemize}
\item
$\chi_W|\R^\times = \sgn^n$ and $\chi_V|\R^\times = \sgn^m$
with $n=p+q=\dim(W)$ and $m=r+s=\dim(V)$, respectively;
\item
$\chi_{W}$ and $\chi_{V}$ depend only on $n \bmod 2$ and $m \bmod 2$, respectively, 
\end{itemize}
and a non-trivial additive character $\psi$ of $\R$.
\par

Let $W=W_{p,q}$ and $V=V_{r,s}$, and 
set $\W = V \otimes_\C W$ as in \S \ref{sec. splitting}.
Then the isometry groups $\U(W)$ and $\U(V)$ are isomorphic to $\U(p,q)$ and $\U(r,s)$, 
respectively.
We have a canonical map
\[
\alpha_V \times \alpha_W \colon \U(W) \times \U(V) \rightarrow \Sp(\W).
\]
As in \S \ref{sec. splitting}, we have a splitting
\[
\cl\alpha_{\chi_V} \times \cl\alpha_{\chi_W} \colon \U(W) \times \U(V) \rightarrow \Mp(\W)
\]
of $\alpha_V \times \alpha_W$.
On the other hand, as in \S \ref{sec. double}, 
there are two-fold covers $\cl\U(W)$ and $\cl\U(V)$ of $\U(W)$ and $\U(V)$, respectively, 
and a lifting
\[
\cl\alpha_V \times \cl\alpha_W \colon \cl\U(W) \times \cl\U(V) \rightarrow \cl\Sp(\W)
\]
of $\alpha_V \times \alpha_W$.
\par

For $a \in \R^\times$, we define an additive character $a\psi$ of $\R$ by
\[
(a\psi)(x) = \psi(ax)
\]
for $x \in \R$.
Let $\omega_{a\psi}$ be the Weil representation of $\Mp(\W)$ associated to $a\psi$.
It satisfies that $\omega_{a\psi}(z) = z \cdot \id$ for $z \in \C^1 \subset \Mp(\W)$.
Moreover, $\omega_{a\psi} \cong \omega_{a'\psi}$ if and only if $aa'>0$.
Hence there are exactly two Weil representations of $\Mp(\W)$.
By the restriction, $\omega_{a\psi}$ is regarded as a representation of $\cl\Sp(\W)$ also.
We choose $\psi$ such that
$\omega_\psi$ is the Weil representation of $\cl\Sp(\W)$
which Paul has used in \cite{P1}, \cite{P2}, \cite{P3} (c.f. \cite[Lemma 1.4.5]{P1}).
We consider two representations $\omega_\psi \circ \cl\alpha_{\chi_V}$ of $\U(W)$
and $\omega_\psi \circ \cl\alpha_V$ of $\cl\U(W)$.
For $\pi \in \Irr_\temp(\U(W))$ and $\cl\pi \in \Irr_\temp(\cl\U(W))$, 
the maximal $\pi$-isotypic quotient of $\omega_\psi \circ \cl\alpha_{\chi_V}$
and the maximal $\cl\pi$-isotypic quotient of $\omega_\psi \circ \cl\alpha_{V}$
are of the form
\[
\pi \boxtimes \Theta_{r,s}(\pi)
\quad\text{and}\quad
\cl\pi \boxtimes \Theta_{r,s}(\cl\pi),
\]
where $\Theta_{r,s}(\pi)$ and $\Theta_{r,s}(\cl\pi)$ are (genuine or possibly zero) representations
of $\U(V)$ and $\cl\U(V)$, respectively.
We call $\Theta_{r,s}(\pi)$ (\resp $\Theta_{r,s}(\cl\pi)$) the big theta lift of $\pi$ (\resp $\cl\pi$).

\begin{thm}[Howe duality correspondence {\cite[Theorem 2.1]{Ho}}]
If $\Theta_{r,s}(\pi)$ (\resp $\Theta_{r,s}(\cl\pi)$) is nonzero, 
then it has a unique irreducible quotient $\theta_{r,s}(\pi)$ (\resp $\theta_{r,s}(\cl\pi)$).
\end{thm}

We interpret $\theta_{r,s}(\pi)$ (\resp $\theta_{r,s}(\cl\pi)$) 
to be zero if so is $\Theta_{r,s}(\pi)$ (\resp $\Theta_{r,s}(\cl\pi)$).
We call $\theta_{r,s}(\pi)$ (\resp $\theta_{r,s}(\cl\pi)$) the small theta lift of $\pi$ (\resp $\cl\pi$).
Similarly, for $\sigma \in \Irr_\temp(\U(V))$ (\resp $\cl\sigma \in \Irr_\temp(\cl\U(V))$), 
we can define the big and small theta lifts
$\Theta_{p,q}(\sigma)$ and $\theta_{p,q}(\sigma)$ 
(\resp $\Theta_{p,s}(\cl\sigma)$ and $\theta_{p,q}(\cl\sigma)$).
\par

In this paper, we determine explicitly 
when $\Theta_{r,s}(\pi)$ (\resp $\Theta_{r,s}(\cl\pi)$) is nonzero for 
$\pi \in \Irr_\temp(\U(p,q))$ (\resp $\cl\pi \in \Irr_\temp(\cl\U(p,q))$).
A relation between the non-vanishing of $\Theta_{r,s}(\pi)$ and 
the one of $\Theta_{r,s}(\cl\pi)$ is given as follows:

\begin{prop}\label{chiV}
Suppose that $\chi_V = \chi_{\nu}$ with $\nu \in \Z$ such that $\nu \equiv m \bmod 2$.
Then for any $\pi \in \Irr_\temp(\U(W))$, we have
\[
\Hom_{\U(W)}(\omega_\psi \circ \cl\alpha_{\chi_V}, \pi) \not= 0
\iff
\Hom_{\cl\U(W)}(\omega_\psi \circ \cl\alpha_V, \pi \otimes {\det}^{-\nu/2}) \not= 0.
\]
In particular, $\Theta_{r,s}(\pi) \not= 0$ if and only if $\Theta_{r,s}(\pi \otimes {\det}^{-\nu/2}) \not= 0$.
\end{prop}
\begin{proof}
For $(g,z) \in \cl\U(W)$, by Proposition \ref{beta}, we have
\begin{align*}
(\omega_\psi \circ \cl\alpha_V) \otimes {\det}^{\nu/2}(g,z)
&= \omega_\psi(\alpha_V(g), \beta(g,z))z\\
&= \omega_\psi(\alpha_V(g), \beta(g,z)z)\\
&= \omega_\psi(\alpha_V(g), \beta_{\chi_V}(g))\\
&= \omega_\psi \circ \cl\alpha_{\chi_V}(g).
\end{align*}
Hence $(\omega_\psi \circ \cl\alpha_V) \otimes {\det}^{\nu/2} = \omega_\psi \circ \cl\alpha_{\chi_V}$
as representations of $\U(W)$, so that
\[
\Hom_{\U(W)}((\omega_\psi \circ \cl\alpha) \otimes {\det}^{\nu/2}, \pi) \not= 0
\iff
\Hom_{\U(W)}(\omega_\psi \circ \cl\alpha_{\chi_V}, \pi) \not= 0.
\]
This completes the proof.
\end{proof}

If $\chi_V = \chi_{\nu}$ with $\nu \in \Z$ such that $\nu \equiv m \bmod 2$, 
the genuine character ${\det}^{\nu/2}$ of $\cl\U(W)$ is also denoted by $\chi_V$.
Hence $\chi_V^2$ is the character ${\det}^\nu$ of $\U(W)$.
By Proposition \ref{chiV}, 
$\Theta_{r,s}(\pi) \not= 0$ if and only if $\Theta_{r,s}(\pi \otimes \chi_V^{-1}) \not= 0$.
\par

Now we recall basic properties on the theta correspondence.
Frist we state a proposition 
which is called the tower property or Kudla's persistence principle.
\begin{prop}[Tower property {\cite{Ku1}}]\label{tower}
If $\Theta_{r,s}(\pi)$ is nonzero, 
then $\Theta_{r+l,s+l}(\pi)$ is also nonzero for any $l \geq 0$.
\end{prop}

Next, we state the conservation relation.
Fix $\pi \in \Irr_\temp(\U(p,q))$.
For each integer $d$, we consider a set of theta lifts $\{\Theta_{r,s}(\pi) \ |\ r-s = d\}$.
We call this set the $d$-th Witt tower of theta lifts of $\pi$.
Also we call 
\[
m_d(\pi) = \min\{r+s\ |\ \Theta_{r,s}(\pi) \not= 0,\ r-s=d\}
\]
the first occurrence index of the $d$-th Witt tower of theta lifts of $\pi$.

\begin{thm}[Conservation relation {\cite{SZ2}}]\label{CR}
Fix $\delta \in \{0,1\}$ and
set 
\[
m_\pm(\pi) = 
\min\{m_d(\pi)\ |\ d \equiv \delta \bmod 2,\ (-1)^{\half{d-\delta}} = \pm 1\}. 
\]
Then 
\[
m_+(\pi) + m_-(\pi) = 2n+2
\]
for any $\pi \in \Irr_\temp(\U(p,q))$ with $n=p+q$.
\end{thm}

The following is a consequence of the induction principle (\cite[Theorem 4.5.5]{P1}).

\begin{prop}[Induction principle]\label{IP}
Let $\pi_0 \in \Irr_\temp(\U(p_0,q_0))$.
Suppose that $\Theta_{r_0,s_0}(\pi_0)$ is nonzero for some $(r_0, s_0)$. 
Let $\xi_1, \dots, \xi_v$ be unitary characters of $\C^\times$.
Put $p = p_0 + v$, $q = q_0 + v$, $r = r_0+v$ and $s = s_0+v$.
Then there exists an irreducible subquotient $\pi$ of the induced representation
\[
\Ind_{P}^{\U(p,q)}(\xi_1 \otimes \dots \otimes \xi_v \otimes \pi_0)
\]
such that $\Theta_{r,s}(\pi)$ is nonzero, 
where $P$ is a parabolic subgroup of $\U(p,q)$
with Levi subgroup of the form $(\C^\times)^v \times \U(p_0, q_0)$.
\end{prop}

For a relation between theta lifts and contragredient representations, 
the following is known.

\begin{prop}\label{vee}
Let $\pi \in \Irr_\temp(\U(p,q))$. 
If $\Theta_{r,s}(\pi) \not= 0$, then $\Theta_{s,r}(\pi^\vee \otimes \chi_V^2) \not= 0$.
\end{prop}
\begin{proof}
By Proposition \ref{chiV}, $\Theta_{r,s}(\pi) \not= 0$ if and only if 
$\Theta_{r,s}(\pi \otimes \chi_{V}^{-1}) \not= 0$.
Similarly, $\Theta_{s,r}(\pi^\vee \otimes \chi_V^2) \not= 0$ if and only if
$\Theta_{s,r}(\pi^\vee \otimes \chi_{V}) \not= 0$.
Hence
the proposition follows from \cite[Proposition 2.1]{P1}.
\end{proof}

There is a non-vanishing result of theta lifts of one dimensional representations.
\begin{prop}\label{trivial}
Fix a positive integer $t$ and a half integer $l$ such that
$2l \equiv t \bmod 2$ and $-t/2 < l < t/2$.
Let ${\det}^l$ be a genuine character of the ${\det}^{l}$-cover $\cl\U(p,q)$.
If $\Theta_{r, s}({\det}^l \otimes \chi_{V_{r,s}}) \not=0$ and $|r-s| = t$, then $\min\{r,s\} \geq p+q$.
\end{prop}
\begin{proof}
Note that ${\det}^l \otimes \chi_{V_{r,s}}$ is a character of $\U(p,q)$.
By Proposition \ref{chiV}, $\Theta_{r,s}({\det}^l \otimes \chi_{V_{r,s}}) \not=0$ 
if and only if $\Theta_{r,s}({\det}^l) \not=0$. 
Hence the proposition is a restatement of \cite[Lemma 3.1]{P2}.
\end{proof}
\par

We denote by $\omega_{p,q,r,s}$ the Weil representation of $\cl\U(p,q) \times \cl\U(r,s)$, 
i.e, $\omega_{p,q,r,s} = \omega_{\psi} \circ (\cl\alpha_{V} \times \cl\alpha_{W})$.
The following are called seesaw identities, 
which are key properties to prove the main result (Theorem \ref{main} below).

\begin{prop}[Seesaw identity]\label{seesaw}
\begin{enumerate}
\item
For $\cl\pi \in \Irr_\temp(\cl\U(p,q))$ and $\cl\sigma \in \Irr_\temp(\cl\U(r,s))$, we have
\[
\Hom_{\cl\U(p,q)}(\Theta_{p+p',q+q'}(\cl\sigma), \cl\pi) \cong 
\Hom_{\cl\U(r,s)}(\Theta_{r,s}(\cl\pi) \otimes \omega_{p',q', r,s}, \cl\sigma).
\]
In particular, if there is $\pi' \in \Irr_\temp(\U(p+p',q+q'))$ such that 
$\Hom_{\U(p,q)}(\pi', \pi) \not= 0$ and $\Theta_{r,s}(\pi') \not= 0$, 
then $\Theta_{r,s}(\pi) \not= 0$.
\item
For $\pi \in \Irr_\temp(\U(p,q))$, 
$\cl\sigma_1 \in \Irr_\temp(\cl\U(r_1,s_1))$ and
$\cl\sigma_2 \in \Irr_\temp(\cl\U(s_2,r_2))$
with $r_1+s_1 \equiv r_2+s_2 \bmod 2$, 
we have
\[
\Hom_{\U(p,q)}(\Theta_{p,q}(\cl\sigma_1) \otimes \Theta_{p,q}(\cl\sigma_2), \pi) \cong 
\Hom_{\cl\U(r_1,s_1) \times \cl\U(s_2,r_2)}
(\Theta_{r_1+s_2,s_1+r_2}(\cl\pi), \cl\sigma_1 \otimes \cl\sigma_2), 
\]
where $\cl\pi$ is the genuine representation of the trivial cover 
$\cl\U(p,q) = \U(p,q) \times \{\pm1\}$ defined by $\cl\pi | \U(p,q) = \pi$.
In particular, for a unitary character $\chi$ of $\U(p,q)$, 
if there is $\pi \in \Irr_\temp(\U(p,q))$ such that 
$\Theta_{r_1,s_1}(\pi) \not= 0$ and $\Theta_{r_2,s_2}(\pi \otimes \chi^{-1}) \not= 0$, 
then $\Theta_{r_1+s_2,s_1+r_2}(\cl\chi \cdot \chi_{V_{0,0}}) \not= 0$.
\end{enumerate}

\end{prop}
\begin{proof}
The first assertions of (1) and (2) immediately follow from \cite[Lemma 2.8]{P1}.
\par

We show the last assertion of (1).
Set $\cl\pi = \pi \otimes \chi_V^{-1}$ and $\cl\pi' = \pi' \otimes \chi_V^{-1}$.
If $\Hom_{\U(p,q)}(\pi', \pi) \not= 0$ and $\Theta_{r,s}(\pi') \not= 0$
then $\Hom_{\cl\U(p,q)}(\cl\pi', \cl\pi) \not= 0$ and $\Theta_{r,s}(\cl\pi') \not= 0$.
If we put $\cl\sigma = \theta_{r,s}(\cl\pi')$, then
$\cl\pi'$ is a quotient of $\Theta_{p+p',q+q'}(\cl\sigma)$. 
Hence we have
\[
0 \not= \Hom_{\cl\U(p,q)}(\cl\pi', \cl\pi)
\subset \Hom_{\cl\U(p,q)}(\Theta_{p+p',q+q'}(\cl\sigma), \cl\pi)
\cong \Hom_{\cl\U(r,s)}(\Theta_{r,s}(\cl\pi) \otimes \omega_{p',q', r,s}, \cl\sigma).
\]
In particular, we have $\Theta_{r,s}(\cl\pi) \not= 0$ so that $\Theta_{r,s}(\pi) \not= 0$.
\par

We show the last assertion of (2).
We denote $\chi_V = \chi_{V_{r_1,s_1}} = \chi_{V_{r_2, s_2}}$ and 
$\chi_{V_{0,0}} = \chi_{V_{r_1+s_2, s_1+r_2}}$.
By Proposition \ref{vee}, we have $\Theta_{s_2,r_2}(\pi^\vee \otimes \chi\chi_V^2) \not= 0$.
Set $\cl\pi_1=\pi \otimes \chi_{V}^{-1}$ and $\cl\pi_2 = \pi^\vee \otimes \chi\chi_{V}$. 
Then $\cl\sigma_1 = \theta_{r_1,s_1}(\cl\pi_1) \not= 0$, 
$\cl\sigma_2 = \theta_{s_2,r_2}(\cl\pi_2) \not=0$, 
and there exist surjections
\[
\Theta_{p,q}(\cl\sigma_1) \otimes \Theta_{p,q}(\cl\sigma_2)
\twoheadrightarrow
\cl\pi_1 \otimes \cl\pi_2 
\twoheadrightarrow
\chi
\]
as representations of $\U(p,q)$.
Since 
\[
0\not=\Hom_{\U(p,q)}(\Theta_{p,q}(\cl\sigma_1) \otimes \Theta_{p,q}(\cl\sigma_2), \chi) 
\cong
\Hom_{\cl\U(r_1,s_1) \times\cl\U(s_2,r_2)}
(\Theta_{r_1+s_2,s_1+r_2}(\cl\chi), \cl\sigma_1 \otimes \cl\sigma_2), 
\]
we have $\Theta_{r_1+s_2, s_1+r_2}(\cl\chi) \not= 0$.
Hence $\Theta_{r_1+s_2, s_1+r_2}(\cl\chi \cdot \chi_{V_{0,0}}) \not= 0$.
\end{proof}

We write Proposition \ref{seesaw} (1) as
\[
\xymatrix{
\U(p+p', q+q') \ar@{-}[d] \ar@{-}[dr] & \U(r, s) \times \U(r, s) \ar@{-}[d]\\
\U(p,q) \times \U(p',q') \ar@{-}[ur] & \U(r, s)
}
\]
and (2) as
\[
\xymatrix{
\U(p, q) \times \U(p,q) \ar@{-}[d] \ar@{-}[dr] & \U(r_1+s_2, s_1+r_2) \ar@{-}[d]\\
\U(p,q) \ar@{-}[ur] & \U(r_1, s_1) \times \U(s_2,r_2),
}
\]
respectively.
Note that in Proposition \ref{seesaw} (2), 
if $\chi = {\det}^a$ for some integer $a$, 
then $\cl\chi$ is the genuine character ${\det}^a$ of the ${\det}^a$-cover $\cl\U(p,q)$.

\subsection{Equal rank case and almost equal rank case}\label{sec. paul}
To prove our main result, 
we use non-trivial results established by Paul in \cite{P1} and \cite{P3}.
These are results on theta lifts in the equal rank case and the almost equal rank case.
In these results, Paul considered the theta lifts from 
irreducible genuine representations of double covers of unitary groups.
In this subsection, we recall Paul's results, and translate them into results of
the theta lifts from irreducible representations of $\U(p,q)$.
\par

Recall that 
\[
\cl\U(p,q) \cong \{(g,z) \in \U(p,q) \times \C^1\ |\ z^2 = \det(g)^{\nu}\}, 
\]
where $\nu \in \Z$ satisfies that $\nu \equiv m \bmod 2$
(so that $\cl\U(p,q)$ depends not only on $(p,q)$ but also on $m \bmod 2$).
It has a genuine character ${\det}^{\nu/2} \colon (g,z) \mapsto z$.
Hence
\[
\Irr_\disc(\cl\U(p,q)) = \{\pi \otimes {\det}^{-\nu/2}\ |\ \pi \in \Irr_\disc(\U(p,q))\}.
\]
Since $\pi \in \Irr_\disc(\U(p,q))$ is characterized by its Harish-Chandra parameter $\HC(\pi)$, 
which is an element of $(\Z+\half{n-1})^p \times (\Z+\half{n-1})^q$ with $n=p+q$, 
the representation $\cl\pi = \pi \otimes {\det}^{-\nu/2}$ is characterized by
its Harish-Chandra parameter
\[
\HC(\cl\pi) = \HC(\pi) - \left(\half{\nu}, \dots, \half{\nu}; \half{\nu}, \dots, \half{\nu}\right)
\in \left(\Z+\half{n+m-1}\right)^p \times \left(\Z+\half{n+m-1}\right)^q.
\]
\par

First, we recall the result in the equal rank case (\cite{P1}).
\begin{thm}[{\cite[Theorems 0.1, 6.1 (a)]{P1}}]\label{P1}
Let $\cl\U(p,q)$ be the ${\det}^{(p+q)/2}$-cover of $\U(p,q)$, 
and $\cl\pi$ be an irreducible tempered genuine representation of $\cl\U(p,q)$.
Then there exists a unique pair $(r,s)$ such that
$r+s = p+q$ and $\Theta_{r,s}(\cl\pi) \not= 0$.
Moreover, if $\cl\pi$ is a direct summand of induced representation 
\[
\Ind_{\cl{P}}^{\cl\U(p,q)}(\xi_1 \otimes \dots \otimes \xi_v \otimes \cl\pi_0), 
\]
where
\begin{itemize}
\item
$\cl{P}$ is a parabolic subgroup of $\cl\U(p,q)$
with Levi subgroup of the form $(\C^\times)^v \times \cl\U(p_0, q_0)$
with $p = p_0+v$ and $q=q_0+v$;
\item
$\cl\pi_0$ is an irreducible genuine discrete series such that
\[
\HC(\cl\pi_0) = 
\left(a_1, \dots, a_x, b_1, \dots b_y; c_1, \dots, c_z, d_1, \dots, d_w\right) 
\in \left(\Z + \half{1}\right)^{p_0} \times \left(\Z + \half{1}\right)^{q_0}
\]
with $a_1 > \dots > a_x > 0 > b_1 > \dots > b_y$ and 
$c_1 > \dots > c_z > 0 > d_1 > \dots > d_w$; 
\item
$\xi_1, \dots, \xi_v$ are unitary characters of $\C^\times$, 
\end{itemize}
then $r=x+w+v$ and $s = y+z+v$.
\end{thm}

We translate this theorem in terms of $L$-parameters.
Fix $m \bmod 2$.
Let $\lam = (\phi, \eta)$ be a pair of 
$\phi \in \Phi_\temp(\U_n(\R))$ and $\eta \in \widehat{A_\phi}$.
Write
\begin{align*}
\phi\chi_V^{-1} =& 
\chi_{2\alpha_1} + \dots + \chi_{2\alpha_{u}} 
+ (\xi_1 + \dots + \xi_v) + ({}^c\xi_1^{-1} + \dots + {}^c\xi_v^{-1}), 
\end{align*}
where
\begin{itemize}
\item
$\alpha_i \in \half{1}\Z$ such that 
$2\alpha_i \equiv n+m-1 \bmod 2$ and $\alpha_1 > \dots > \alpha_{u}$;
\item
$\xi_i$ is a unitary character of $\C^\times$ (which can be of the form $\chi_{2\alpha}$);
\item
$u+ 2v=n$.
\end{itemize}
Then 
\[
A_\phi \supset (\Z/2\Z) e_{V, 2\alpha_1} \oplus \dots \oplus (\Z/2\Z) e_{V, 2\alpha_{u}}, 
\]
where $\{e_{V, 2\alpha_1}, \dots, e_{V, 2\alpha_{u}}\}$
is the canonical basis associated to $\{\chi_V\chi_{2\alpha_1}, \dots, \chi_V\chi_{2\alpha_{u}}\}$.
\par

The following theorem is a translation of Theorem \ref{P1}.
\begin{thm}\label{PE}
Suppose that $m=n=p+q$.
Let $\lam = (\phi, \eta)$ be as above, and set $\pi = \pi(\phi,\eta)$.
Then there exists a unique pair $(r,s)$ such that
$r+s = p+q$ and $\Theta_{r,s}(\pi) \not= 0$.
Moreover $(r,s)$ is given by 
\[
\left\{
\begin{aligned}
r &= \#\{i \in \{1, \dots, u\}\ |\ (-1)^{i-1}\eta(e_{V,2\alpha_i})\alpha_i > 0\}+v,\\
s &= \#\{i \in \{1, \dots, u\}\ |\ (-1)^{i-1}\eta(e_{V,2\alpha_i})\alpha_i < 0\}+v.
\end{aligned}
\right.
\]
\end{thm}

Similarly, we can translate the result in the almost equal rank case (\cite{P3}).
\begin{thm}\label{PA}
Suppose that $m \equiv n+1 \bmod 2$.
Let $\lam = (\phi, \eta)$ be as above, and set $\pi = \pi(\phi,\eta)$.
\begin{enumerate}
\item
Suppose that $\phi$ does not contain $\chi_V$.
Then there exist exactly two pairs $(r,s)$ such that 
$r+s = p+q+1$ and $\Theta_{r,s}(\pi) \not= 0$.
Moreover $(r,s)$ is given by 
\[
\left\{
\begin{aligned}
r &= \#\{i \in \{1, \dots, u\}\ |\ (-1)^{i-1}\eta(e_{V,2\alpha_i})\alpha_i > 0\}+v+1,\\
s &= \#\{i \in \{1, \dots, u\}\ |\ (-1)^{i-1}\eta(e_{V,2\alpha_i})\alpha_i < 0\}+v
\end{aligned}
\right.
\]
and
\[
\left\{
\begin{aligned}
r &= \#\{i \in \{1, \dots, u\}\ |\ (-1)^{i-1}\eta(e_{V,2\alpha_i})\alpha_i > 0\}+v,\\
s &= \#\{i \in \{1, \dots, u\}\ |\ (-1)^{i-1}\eta(e_{V,2\alpha_i})\alpha_i < 0\}+v+1.
\end{aligned}
\right.
\]

\item
Suppose that $\phi$ contains $\chi_V$ with odd multiplicity.
Then there exists a unique pair $(r,s)$ such that
$r+s = p+q-1$ and $\Theta_{r,s}(\pi) \not= 0$.
Moreover $(r,s)$ is given by 
\[
\left\{
\begin{aligned}
r &= \#\{i \in \{1, \dots, u\}\ |\ (-1)^{i-1}\eta(e_{V,2\alpha_i})\alpha_i > 0\}+v,\\
s &= \#\{i \in \{1, \dots, u\}\ |\ (-1)^{i-1}\eta(e_{V,2\alpha_i})\alpha_i < 0\}+v.
\end{aligned}
\right.
\]

\item
Suppose that $\phi$ contains $\chi_V$ with even multiplicity.
Then there exists a unique pair $(r,s)$ such that
$r+s = p+q-1$ and $\Theta_{r,s}(\pi) \not= 0$.
Moreover $(r,s)$ is given by 
\[
\left\{
\begin{aligned}
r &= \#\{i \in \{1, \dots, u\}\ |\ (-1)^{i-1}\eta(e_{V,2\alpha_i})\alpha_i > 0\}+v-1,\\
s &= \#\{i \in \{1, \dots, u\}\ |\ (-1)^{i-1}\eta(e_{V,2\alpha_i})\alpha_i < 0\}+v
\end{aligned}
\right.
\]
or
\[
\left\{
\begin{aligned}
r &= \#\{i \in \{1, \dots, u\}\ |\ (-1)^{i-1}\eta(e_{V,2\alpha_i})\alpha_i > 0\}+v,\\
s &= \#\{i \in \{1, \dots, u\}\ |\ (-1)^{i-1}\eta(e_{V,2\alpha_i})\alpha_i < 0\}+v-1.
\end{aligned}
\right.
\]

\end{enumerate}
\end{thm}

In fact, Paul (\cite[Theorem 3.4]{P3}) has determined $(r,s)$ in Theorem \ref{PA} (3) exactly
in terms of a system of positive roots.

\section{The definition and the main result}\label{result}
In this section, we state the main result and its corollary.

\subsection{Definition}
Before stating the main result, 
we define some notations.

\begin{defi}\label{def}
Fix $\kappa \in \{1,2\}$. 
Choose a character $\chi_V$ of $\C^\times$ such that $\chi_V|\R^\times = \sgn^{\kappa+n}$.
Let $\lam = (\phi, \eta)$ be a pair of $\phi \in \Phi_\temp(\U_n(\R))$ and $\eta \in \widehat{A_\phi}$.

\begin{enumerate}
\item
Consider the set $\TT$ containing $\kappa - 2$ 
and all integers $k > 0$ with $k \equiv \kappa \bmod 2$ satisfying the following conditions:
\begin{description}
\item[(chain condition)]
$\phi$ contains $\chi_V\chi_{k-1} + \chi_V\chi_{k-3} + \dots + \chi_V\chi_{-k+1}$;
\item[(odd-ness condition)] 
the multiplicity of $\chi_V\chi_{k+1-2i}$ in $\phi$ 
is odd for $i = 1, \dots, k$;
\item[(alternating condition)]
$\eta(e_{V,k+i-2i}) = -\eta(e_{V, k-1-2i})$ for $i=1, \dots, k-1$.
\end{description}
Here, $e_{V,2\alpha}$ is the element in $A_\phi$ corresponding to $\chi_V\chi_{2\alpha}$.
Set
\[
k_\lam = \max \, \TT.
\]

\item
Write
\begin{align*}
\phi\chi_V^{-1} =& 
\chi_{2\alpha_1} + \dots + \chi_{2\alpha_{u}} 
+ (\xi_1 + \dots + \xi_v) + ({}^c\xi_1^{-1} + \dots + {}^c\xi_v^{-1}), 
\end{align*}
where
\begin{itemize}
\item
$\alpha_i \in \half{1}\Z$ such that 
$2\alpha_i \equiv \kappa-1 \bmod 2$ and $\alpha_1 > \dots > \alpha_u$;
\item
$\xi_j$ is a unitary character of $\C^\times$ (which can be of the form $\chi_{2\alpha}$);
\item
$u+2v=n$.
\end{itemize}
Then $A_\phi \supset (\Z/2\Z)e_{V, 2\alpha_1} + \dots + (\Z/2\Z)e_{V,2\alpha_u}$.
Define $(r_\lam, s_\lam)$ by 
\[
\left\{
\begin{aligned}
r_\lam =& \#\{i \in \{1, \dots, u\}\ |\ 
|\alpha_i| \geq (k_\lam+1)/2,\ (-1)^{i-1}\eta(e_{V,\alpha_i})\alpha_i > 0 \} + v, \\
s_\lam =& \#\{i \in \{1, \dots, u\}\ |\ 
|\alpha_i| \geq (k_\lam+1)/2,\ (-1)^{i-1}\eta(e_{V,\alpha_i})\alpha_i < 0 \} + v.
\end{aligned}
\right.
\]

\item
Write
\begin{align*}
\phi\chi_V^{-1} =& 
m_1\chi_{2\alpha_1} + \dots + m_u\chi_{2\alpha_{u}} 
+ m'_1\chi_{2\alpha'_1} + \dots + m'_{u'}\chi_{2\alpha'_{u'}}
+ (\xi_1 + \dots + \xi_v) + ({}^c\xi_1^{-1} + \dots + {}^c\xi_v^{-1}), 
\end{align*}
where
\begin{itemize}
\item
$\alpha_i, \alpha'_{i'} \in \half{1}\Z$ such that 
$2\alpha_i \equiv 2\alpha'_{i'} \equiv \kappa-1 \bmod 2$, 
$\alpha_1 > \dots > \alpha_u$, $\alpha'_1 > \dots > \alpha'_{u'}$ and 
$\{\alpha_1, \dots, \alpha_u\} \cap \{\alpha'_1, \dots, \alpha'_{u'}\} = \emptyset$;
\item
$m_i \geq 1$ (\resp $m'_{i'} \geq 1$) is the multiplicity of 
$\chi_{2\alpha_1}$ (\resp $\chi_{2\alpha'_{i'}}$)
in $\phi\chi_V^{-1}$ such that $m_i$ is odd (\resp $m'_{i'}$ is even);
\item
$\xi_j$ is a unitary character of $\C^\times$ which is not of the form $\chi_{2\alpha}$
with $2\alpha \equiv \kappa-1 \bmod 2$;
\item
$(m_1 + \dots + m_u) + (m'_1 + \dots + m'_{u'}) + 2v=n$.
\end{itemize}
Define a subset $X_\lam$ of $\half{1}\Z \times \{\pm1\}$ by
\begin{align*}
X_\lam = 
&\{(\alpha_i, (-1)^{i-1}\eta(e_{V,2\alpha_i}))\ |\ i = 1, \dots, u\} 
\\&
\cup
\{(\alpha'_{i'}, +1), (\alpha'_{i'}, -1)\ |\ i' = 1, \dots, u',\ 
\eta(e_{V, \alpha'_{i'}}) \not= (-1)^{\#\{i \in \{1, \dots, u\}\ |\ \alpha_i > \alpha'_{i'}\}}
\}.
\end{align*}

\item
We define a sequence 
$X_\lam=X_\lam^{(0)} \supset X_\lam^{(1)} \supset \dots 
\supset X_\lam^{(n)} \supset \cdots$ as follows:
Let $\{\beta_1, \dots, \beta_{u_j}\}$ be the image of $X_\lam^{(j)}$ 
under the projection $\half{1}\Z \times \{\pm1\} \rightarrow \half{1}\Z$
such that $\beta_1 > \dots > \beta_{u_j}$.
Set $S_j$ to be the set of $i \in \{2, \dots, u_j\}$ such that
\begin{itemize}
\item
$(\beta_{i-1}, +1), (\beta_i, -1) \in X_\lam^{(j)}$; 
\item
$\min\{|\beta_{i-1}|, |\beta_i|\} \geq (k_\lam+1)/2$;
\item
$\beta_{i-1}\beta_i \geq 0$.
\end{itemize}
Then we define a subset $X_\lam^{(j+1)}$ of $X_\lam^{(j)}$ by
\[
X_\lam^{(j+1)} 
= X_\lam^{(j)} \setminus \left(\bigcup_{i \in S_j}\{(\beta_{i-1}, +1), (\beta_i, -1)\}\right).
\]
Finally, we set $X_\lam^{(\infty)} = X_\lam^{(n)} = X_\lam^{(n+1)}$.

\item
For an integer $T$ and $\epsilon \in \{\pm1\}$, we define a set $\CC^\epsilon_\lam(T)$ by 
\[
\CC^\epsilon_\lam(T) = \left\{(\alpha, \epsilon) \in X_\lam^{(\infty)} |\ 
0 \leq \epsilon\alpha + \half{k_\lam-1} < T
\right\}.
\]
In particular, if $T \leq 0$, then $\CC^\epsilon_\lam(T) = \emptyset$.

\end{enumerate}
\end{defi}

\subsection{Main result}
The main result is the following:
\begin{thm}\label{main}
Assume the local Gan--Gross--Prasad conjecture (Conjecture \ref{GGP}).
Let $\lam = (\phi, \eta)$ be a pair of $\phi \in \Phi_\temp(\U_n(\R))$ and $\eta \in \widehat{A_\phi}$.
Set $\pi = \pi(\phi,\eta)$.
Let $k=k_\lam$, $r=r_\lam$, $s=s_\lam$ and $\CC^\epsilon_\lam(T)$ 
be as in Definition \ref{def}.
\begin{enumerate}
\item
Suppose that $k = -1$.
Then for integers $l$ and $t \geq 1$, 
the theta lift $\Theta_{r+2t+l+1, s+l}(\pi)$ is nonzero if and only if
\[
l \geq 0
\quad\text{and} \quad
\#\CC^\epsilon_\lam(t+l) \leq l
\quad\text{for each $\epsilon \in \{\pm1\}$}.
\]
Moreover, for an integer $l$, 
the theta lift $\Theta_{r+l+1, s+l}(\pi)$ is nonzero if and only if
\[
\left\{
\begin{aligned}
&l \geq 0 \iif \text{$\phi$ does not contain $\chi_V$},\\
&l \geq -1 \iif \text{$\phi$ contains $\chi_V$ and $(0,\pm1) \not\in X_\lam$},\\
&l \geq 1 \iif \text{$\phi$ contains $\chi_V$ and $(0,\pm1) \in X_\lam$}.
\end{aligned}
\right.
\]

\item
Suppose that $k \geq 0$.
Then for integers $l$ and $t \geq 1$, 
the theta lift $\Theta_{r+2t+l, s+l}(\pi)$ is nonzero if and only if
\[
l \geq k
\quad\text{and} \quad
\#\CC^\epsilon_\lam(t+l) \leq l
\quad\text{for each $\epsilon \in \{\pm1\}$}.
\]
Moreover, we consider the following three conditions:
\begin{description}
\item[(chain condition 2)]
$\phi\chi_V^{-1}$ contains both $\chi_{k+1}$ and $\chi_{-(k+1)}$, 
so that 
\[
\phi \chi_V^{-1} \supset 
\underbrace{\chi_{k+1} + \chi_{k-1} + \dots + \chi_{-(k-1)} + \chi_{-(k+1)}}_{k+2};
\]
\item[(even-ness condition)]
at least one of $\chi_{k+1}$ and $\chi_{-(k+1)}$ is 
contained in $\phi\chi_V^{-1}$ with even multiplicity;
\item[(alternating condition 2)]
$\eta(e_{V, k+1-2i}) \not= \eta(e_{V, k-1-2i})$ for $i=0, \dots, k$. 
\end{description}
Then for an integer $l$, 
the theta lift $\Theta_{r+l, s+l}(\pi)$ is nonzero if and only if 
\[
\left\{
\begin{aligned}
&l \geq -1 \iif 
\text{$\lam = (\phi, \eta)$ satisfies the three conditions},\\
&l \geq 0 \other.
\end{aligned}
\right.
\]

\end{enumerate}
\end{thm}

\begin{rem}
\begin{enumerate}
\item
When $\phi \in \Phi_\disc(\U_n(\R))$, 
we need the local Gan--Gross--Prasad conjecture only for discrete series representations.
Since He \cite{He} has established the conjecture in this case, 
the statements in Theorem \ref{main} for discrete series representations hold unconditionally.

\item
If $\phi \in \Phi_\disc(\U_n(\R))$, 
then $\pi$ is a discrete series representation of some $\U(p,q)$.
By Theorem \ref{LLC} (4), 
we can translate Definition \ref{def} and Theorem \ref{main} 
in terms of Harish-Chandra parameters, 
and we obtain Definition \ref{intro.def} and Theorem \ref{intro.main}.

\item
When $\phi \in \Phi_\disc(\U_n(\R))$ and $t = 0, 1$, 
Theorem \ref{main} is a translation of results of 
Paul (\cite[Proposition 3.4, Theorem 3.14]{P2}).
\item
When $\pi$ is a representation of a compact unitary group and $l = 0$, 
Theorem \ref{main} is compatible with results of \cite{KV} and \cite{Li}
(see also \cite[Proposition 6.6]{A}).
\end{enumerate}
\end{rem}
\par

By Proposition \ref{vee} together with the following lemma, 
we can obtain the first occurrence index of any Witt tower of theta lifts of 
any irreducible tempered representation $\pi$ of $\U(p,q)$.

\begin{lem}\label{Xv}
Let $\lam = (\phi, \eta)$ as in Definition \ref{def}, and
set $\lam^\vee = (\phi^\vee \otimes \chi_V^2, \eta^\vee)$.
\begin{enumerate}
\item
We have $k_{\lam^\vee} = k_\lam$ and $(r_{\lam^\vee}, s_{\lam^\vee}) = (s_\lam,r_\lam)$.
\item
Suppose that $\chi_{2\alpha}$ is contained in $\phi\chi_V^{-1}$ with odd multiplicity.
Then $(\alpha, \epsilon) \in X_\lam$ if and only if $(-\alpha, \epsilon) \in X_{\lam^\vee}$.
\item
Suppose that $\chi_{2\alpha}$ is contained in $\phi\chi_V^{-1}$ with even multiplicity.
Then $(\alpha, \pm1) \in X_\lam$ if and only if $(-\alpha, \pm1) \not\in X_{\lam^\vee}$.
\item
In general, if $(\alpha, \epsilon) \in X_{\lam}$, then $(-\alpha, -\epsilon) \not \in X_{\lam^\vee}$.
\end{enumerate}
\end{lem}
\begin{proof}
Write
\[
\phi\chi_V^{-1} =
\chi_{2\alpha_1} + \dots + \chi_{2\alpha_{u}} 
+ (\xi_1 + \dots + \xi_v) + ({}^c\xi_1^{-1} + \dots + {}^c\xi_v^{-1})
\]
as in Definition \ref{def} (2).
Then 
\[
\phi^\vee\chi_V =
\chi_{-2\alpha_u} + \dots + \chi_{-2\alpha_{1}} 
+ (\xi_1^{-1} + \dots + \xi_v^{-1}) + ({}^c\xi_1 + \dots + {}^c\xi_v). 
\]
\par

Suppose that $\chi_{2\alpha}$ is contained in $\phi\chi_V^{-1}$ with odd multiplicity.
This means that $\alpha = \alpha_i$ for some $i$.
Then
\begin{align*}
(\alpha, \epsilon) \in X_\lam
&\iff
(-1)^{i-1}\eta(e_{V,2\alpha}) = \epsilon,\\
&\iff
(-1)^{u-i}\eta^\vee(e_{V,-2\alpha}) = \epsilon
\iff (-\alpha, \epsilon) \in X_{\lam^\vee}.
\end{align*}
Hence we have (2).
This easily implies that 
$k_{\lam^\vee} = k_\lam$ and $(r_{\lam^\vee}, s_{\lam^\vee}) = (s_\lam,r_\lam)$.
Hence we have (1).
\par

Now suppose that $\chi_{2\alpha}$ is contained in $\phi\chi_V^{-1}$ with even multiplicity.
Then 
\begin{align*}
(\alpha, \pm1) \in X_\lam
&\iff
\eta(e_{V,2\alpha}) = (-1)^{\#\{ i \in \{1, \dots, u\} \ |\ \alpha_i > \alpha\} + 1},\\
&\iff
\eta^\vee(e_{V,-2\alpha}) = (-1)^{\#\{ i \in \{1, \dots, u\} \ |\ -\alpha_i > -\alpha\}}
\iff (-\alpha, \pm1) \not\in X_{\lam^\vee}.
\end{align*}
Hence we have (3).
The assertion (4) follows from (2) and (3).
\end{proof}

\subsection{A corollary}
As a consequence of Theorem \ref{main}, 
we obtain a new relation between theta lifts and induced representations.

\begin{cor}\label{ind}
Assume the local Gan--Gross--Prasad conjecture (Conjecture \ref{GGP}).
Let $\pi \in \Irr_\temp(\U(p,q))$ and $\pi_0 \in \Irr_\temp(\U(p-1,q-1))$.
Suppose that there exists a unitary character $\xi$ such that 
\[
\pi \subset \Ind_P^{\U(p,q)}(\xi \otimes \pi_0), 
\]
where $P$ is a parabolic subgroup of $\U(p,q)$ 
with Levi subgroup $M_P = \C^\times \times \U(p-1,q-1)$.
For $(r,s)$, we have the following:
\begin{enumerate}
\item
Suppose that $\Theta_{r-1,s-1}(\pi_0) \not=0$. 
If $r+s \leq p+q$, then $\Theta_{r,s}(\pi) \not= 0$.
\item
Suppose that $\Theta_{r,s}(\pi) \not= 0$. 
If $r+s \geq p+q+1$, then $\Theta_{r-1,s-1}(\pi_0) \not=0$.
In general, $\Theta_{r,s}(\pi_0) \not=0$.
\end{enumerate}
\end{cor}
\begin{proof}
Let $\lam = (\phi,\eta)$ and $\lam_0 = (\phi_0,\eta_0)$ be the $L$-parameters 
of $\pi$ and $\pi_0$, respectively.
By Theorem \ref{LLC} (5), we have 
$\phi = \phi_0 + \xi + {}^c\xi^{-1}$ and $\eta|A_{\phi_0} = \eta_0$.
In particular, we see that $k_\lam = k_{\lam_0}$, 
$(r_{\lam}, s_{\lam}) = (r_{\lam_0}+1, s_{\lam_0}+1)$
and $X_{\lam_0} \subset X_\lam$.
\par

We show (1).
By Theorem \ref{main}, 
if $\Theta_{r-1,s-1}(\pi_0) \not=0$ and $r+s \leq p+q$,
then $|(r-s) - (r_{\lam_0}-s_{\lam_0})| \leq 1$.
Hence (1) follows from the last assertions of Theorem \ref{main} (1) and (2).
\par

We show (2).
To prove the first part, it suffices to show that
$\#\CC^\epsilon_{\lam_0}(T) \leq \#\CC^\epsilon_\lam(T)$ 
for any $T$ and $\epsilon \in \{\pm1\}$.
If $X_{\lam_0} = X_\lam$, 
it is clear that $\CC^\epsilon_{\lam_0}(T) = \CC^\epsilon_\lam(T)$ 
for any $T$ and $\epsilon \in \{\pm1\}$.
Hence we may assume that $\xi = \chi_V\chi_{2\alpha}$ 
and $(\alpha, \pm1) \in X_\lam \setminus X_{\lam_0}$.
First, we assume that $\alpha \geq (k_\lam+1)/2$.
Then $\CC^-_{\lam}(T) = \CC^-_{\lam_0}(T)$ for any $T$.
We note that
\[
\#\{(\alpha_0,+1) \in X_{\lam_0}^{(\infty)}\ |\ (\alpha_0,+1) \not\in X_\lam^{(\infty)}\} \leq 1. 
\]
If $X_{\lam_0}^{(\infty)} \subset X_\lam^{(\infty)}$, 
then $\CC^+_{\lam_0}(T) \subset \CC^+_\lam(T)$. 
Suppose that $(\alpha_0,+1) \in X_{\lam_0}^{(\infty)}$ 
but $(\alpha_0,+1) \not\in X_\lam^{(\infty)}$.
Then one of the following holds:
\begin{itemize}
\item
$X_\lam^{(\infty)} = (X_{\lam_0}^{(\infty)} \setminus \{(\alpha_0,+1)\}) \cup \{(\alpha', +1)\}$
for some $(\alpha', +1) \not\in X_{\lam_0}^{(\infty)}$ with $\alpha' \geq (k_\lam+1)/2$;
\item
$X_\lam^{(\infty)} = X_{\lam_0}^{(\infty)} \setminus \{(\alpha_0,+1), (\alpha', -1)\}$
for some $(\alpha', -1) \in X_{\lam_0}^{(\infty)}$ with $\alpha' \geq (k_\lam+1)/2$.
\end{itemize}
In both of two cases, we must have $\alpha' < \alpha_0$.
However, in the second case, $X_{\lam_0}^{(\infty)}$ contains both 
$(\alpha_0,+1)$ and $(\alpha', -1)$ with $\alpha_0 > \alpha' \geq (k_\lam+1)/2$.
This contradicts the definition of $X_{\lam_0}^{(\infty)}$ (see Definition \ref{def} (4)).
Hence we must have 
$X_\lam^{(\infty)} = (X_{\lam_0}^{(\infty)} \setminus \{(\alpha_0,+1)\}) \cup \{(\alpha', +1)\}$
for some $(\alpha', +1) \not\in X_{\lam_0}^{(\infty)}$ with 
$\alpha_0 > \alpha' \geq (k_\lam+1)/2$.
Then we have
\[
\#\CC^+_{\lam}(T) = 
\left\{
\begin{aligned}
&\#\CC^+_{\lam_0}(T)+1	\iif \alpha'-\half{k_\lam-1} < T \leq \alpha_0-\half{k_\lam-1},\\
&\#\CC^+_{\lam_0}(T)	\other.
\end{aligned}
\right.
\]
Therefore in any case, we have 
\[
\#\CC^+_{\lam_0}(T) \leq \#\CC^+_{\lam}(T).
\]
Similarly, if $\alpha \leq -(k_\lam+1)/2$, then 
$\CC^+_{\lam}(T) = \CC^+_{\lam_0}(T)$ and 
$\#\CC^-_{\lam_0}(T) \leq \#\CC^-_{\lam}(T)$ for any $T$.
Hence we have the first part of (2).
The last part follows from the last assertions of Theorem \ref{main} (1) and (2).
\end{proof}

\begin{rem}
A non-archimedean analogue also holds (see Theorem \ref{AGmain}).
In the non-archimedean case, the part (2) is a corollary of Kudla's filtration \cite{Ku1}, 
whereas, in the archimedean case, 
it would not follow from the induction principle (\cite[Theorem 4.5.5]{P1}).
\end{rem}

\section{Proof of Theorem \ref{main}}\label{proof}
In this section, we shall prove Theorem \ref{main}.
In \S \ref{nonvanish}, we show the sufficient conditions of the non-vanishing 
of theta lifts in Theorem \ref{main} (1), (2) when $t \geq 1$.
The proof is an induction using a seesaw identity (Proposition \ref{seesaw} (1)).
To use such a seesaw identity, 
for given $\pi \in \Irr_\temp(\U(p,q))$, 
we have to find a ``good'' representation $\pi'$ of $\U(p+1,q)$ 
such that $\Hom_{\U(p,q)}(\pi',\pi) \not= 0$. 
To do this, we use the Gan--Gross--Prasad conjecture (Conjecture \ref{GGP}) in \S \ref{find}.
The necessary conditions of the non-vanishing 
of theta lifts in Theorem \ref{main} (1), (2) when $t \geq 1$
are proven in \S \ref{vanish}.
For the proof, we use a seesaw identity (Proposition \ref{seesaw} (2))
and Proposition \ref{trivial}.
Finally, using the conservation relation (Theorem \ref{CR}), 
we show Theorem \ref{main} (1), (2) when $t=0$ in \S \ref{t=0}.

\subsection{Finding a GGP pair}\label{find}
A main tool in the proof of the main result is 
a seesaw identity (Proposition \ref{seesaw} (1)).
To use it, 
in this subsection, for given $\pi \in \Irr_\temp(\U(p,q))$, 
we find a ``good'' representation $\pi'$ of $\U(p+1,q)$ 
such that $\Hom_{\U(p,q)}(\pi',\pi) \not= 0$.
\par

For a pair $\lam = (\phi, \eta)$ of $\phi \in \Phi_\temp(\U_n(\R))$ and $\eta \in \widehat{A_\phi}$, 
let $k_\lam$, $r_\lam$, $s_\lam$ and $X_\lam$ be as in Definition \ref{def}.
Consider a representation $\phi''_0$ of $\C^\times$ defined so that
\[
\phi''_0\chi_V^{-1} = \bigoplus_{(\alpha, \epsilon) \in X_\lam} \chi_{2\alpha-\epsilon}.
\]
Note that $\dim(\phi_0'') \equiv n \bmod 2$.
For each $\beta \in \half{1}\Z$ with $2\beta \equiv \kappa \bmod 2$, 
the multiplicity of $\chi_{2\beta}$ in $\phi''_0\chi_V^{-1}$ is at most $2$.
Moreover, it is equal to $2$ if and only if 
both $(\beta+1/2, +1)$ and $(\beta-1/2, -1)$ are contained in $X_\lam$.
Define a representation $\phi'_0$ of $\C^\times$ so that
\[
\phi'_0\chi_V^{-1} = \phi''_0\chi_V^{-1} - 
\bigoplus\{2\chi_{2\beta}\ |\ (\beta+1/2, +1), (\beta-1/2, -1) \in X_\lam \}.
\]
Then $\phi'_0$ is multiplicity-free and $\dim(\phi_0') \equiv n \bmod 2$.
Put $v' = (n-\dim(\phi'_0))/2$.
\par

In this subsection, 
we choose arbitrary half integers $\beta_0, \beta_1, \ldots$ such that 
$2\beta_j \equiv \kappa \bmod 2$ and 
\[
\max\{\alpha_1+1,0\} < \beta_0 < \beta_1 < \cdots.
\] 
\par

To use a seesaw identity (Proposition \ref{seesaw} (1)), 
we need the following lemmas.

\begin{lem}\label{k=-1}
Assume the local Gan--Gross--Prasad conjecture (Conjecture \ref{GGP}).
Let $\lam = (\phi, \eta)$ be a pair of $\phi \in \Phi_\temp(\U_n(\R))$ 
and $\eta \in \widehat{A_\phi}$.
Suppose that $k_\lam=-1$.
Set 
\[
\phi' = \phi_0' + \chi_V(\chi_{2\beta_0} + \chi_{2\beta_1} + \dots + \chi_{2\beta_{2v'}})
\in \Phi_\disc(\U_{n+1}(\R)).
\]
We define $\eta' \in \widehat{A_{\phi'}}$ by setting 
$\eta'(e_{V,2\beta_j}) = 1$ for $j = 0, 1, \dots, 2v'$, and
\[
\eta'(e_{V,2\beta}) = -\eta(e_{V,2\alpha})
\]
when $\beta = \alpha-\epsilon/2$ with $(\alpha,\epsilon) \in X_\lam$.
Let $\lam' = (\phi',\eta')$ and $\pi'=\pi(\phi',\eta')$.
Then 
\begin{itemize}
\item
$\pi' \in \Irr_\disc(\U(p+1, q))$ and $\Hom_{\U(p,q)}(\pi',\pi) \not= 0$;
\item
$k_{\lam'} = 0$ and 
\[
(r_{\lam'}, s_{\lam'}) = 
\left\{
\begin{aligned}
&(r_\lam, s_\lam+1) \iif (0,\pm1) \in X_\lam,\\
&(r_\lam+1, s_\lam) \other;
\end{aligned}
\right.
\]
\item
for $T \leq \beta_0-1/2$ and $\epsilon \in \{\pm1\}$, 
the map $(\beta, \epsilon) \mapsto (\beta + \epsilon/2, \epsilon)$ gives an injection
\[
\CC^\epsilon_{\lam'}(T) \hookrightarrow \CC^\epsilon_{\lam}(T);
\]
\item
for fixed $l \geq 0$ and $1 \leq T \leq \beta_0-1/2$, 
if $\#\CC^\epsilon_\lam(T) \leq l$ for each $\epsilon \in \{\pm1\}$, 
then $\#\CC^\epsilon_{\lam'}(T) \leq l$ for each $\epsilon \in \{\pm1\}$.
\end{itemize}
\end{lem}

\begin{lem}\label{k=0}
Assume the local Gan--Gross--Prasad conjecture (Conjecture \ref{GGP}).
Let $\lam = (\phi, \eta)$ be a pair of $\phi \in \Phi_\temp(\U_n(\R))$ 
and $\eta \in \widehat{A_\phi}$.
Suppose that $k_\lam=0$.
Put $\alpha_+ = \min\{\alpha \ |\ (\alpha, \epsilon) \in X_\lam, \alpha>0\}$ and
$\alpha_- = \max\{\alpha \ |\ (\alpha, \epsilon) \in X_\lam, \alpha<0\}$.
We consider the following cases separately:
\begin{description}
\item[Case $1$]
Both $(1/2, +1)$ and $(-1/2, -1)$ are not contained in $X_\lam$.
\item[Case $2$]
Both $(1/2, +1)$ and $(-1/2, -1)$ are contained in $X_\lam$.
\item[Case $3$]
There is $\delta \in \{\pm1\}$ such that
\[
\alpha_{\delta} = \delta/2, \quad 
(\alpha_\delta, \delta) \in X_\lam
\quad\text{and}\quad 
\alpha_{-\delta} \not=-\delta/2, \quad
(\alpha_{-\delta}, -\delta) \in X_\lam.
\]
\item[Case $4$]
There is $\delta \in \{\pm1\}$ such that
\[
\alpha_\delta=\delta/2, \quad 
(\alpha_\delta,\delta) \in X_\lam
\quad\text{and}\quad
\alpha_{-\delta} \not=-\delta/2, \quad
(\alpha_{-\delta},-\delta) \not\in X_\lam.
\]
\item[Case $5$]
There is $\delta \in \{\pm1\}$ such that
\[
\alpha_\delta=\delta/2, \quad 
(\alpha_\delta,\delta) \in X_\lam
\quad\text{and}\quad
\alpha_{-\delta} =-\delta/2, \quad
(\alpha_{-\delta},-\delta) \not\in X_\lam.
\]
\end{description}
The case $5$ cannot occur if $\phi \in \Phi_\disc(\U_n(\R))$.
We set $\phi' \in \Phi_\disc(\U_{n+1}(\R))$ so that
\[
\phi'\chi_V^{-1}= \left\{
\begin{aligned}
&\phi_0'\chi_V^{-1} + (\chi_{2\beta_0} + \chi_{2\beta_1} + \dots + \chi_{2\beta_{2v'}})
\iif \text{$\lam$ is in case $1$, $2$ or $5$},\\
&\phi_0'\chi_V^{-1} - (\chi_{2\beta_{+}}+\chi_{2\beta_{-}})
+ (\chi_{2\beta_0} + \chi_{2\beta_1} + \dots + \chi_{2\beta_{2v'+2}})
\iif \text{$\lam$ is in case $3$},\\
&\phi_0'\chi_V^{-1} - \1 + \chi_{-2\delta}
+ (\chi_{2\beta_0} + \chi_{2\beta_1} + \dots + \chi_{2\beta_{2v'}})
\iif \text{$\lam$ is in case $4$}.
\end{aligned}
\right.
\]
Here, in the case $3$, we put $\beta_\pm = \alpha_{\pm} \mp 1/2$.
Also we define $\eta' \in \widehat{A_{\phi'}}$ by setting 
$\eta'(e_{V,2\beta_j}) = 1$ for $j= 0,1, \ldots$, and 
\[
\eta'(e_{V,2\beta}) = -\eta(e_{V,2\alpha})
\]
when $\beta = \alpha-\epsilon/2$ with $(\alpha,\epsilon) \in X_\lam$.
In the case $4$, we set $\eta'(e_{V, -2\delta}) = -\eta(e_{V, 2\alpha_\delta})$.
Let $\lam' = (\phi',\eta')$ and $\pi'=\pi(\phi',\eta')$.
Then 
\begin{itemize}
\item
$\pi' \in \Irr_\disc(\U(p+1, q))$ and $\Hom_{\U(p,q)}(\pi',\pi) \not= 0$;
\item
$k_{\lam'}$ and $(r_{\lam'}, s_{\lam'})$ are given by
\begin{align*}
k_{\lam'} &= \left\{
\begin{aligned}
&-1 \iif \text{$\lam$ is in case $1$, $2$, $3$ or $4$},\\
&1	\iif \text{$\lam$ is in case $5$},
\end{aligned}
\right.
\\
(r_{\lam'}, s_{\lam'}) &= \left\{
\begin{aligned}
&(r_{\lam}+1, s_{\lam}) \iif \text{$\lam$ is in case $1$},\\
&(r_{\lam}, s_{\lam}+1) \iif \text{$\lam$ is in case $2$, $3$ or $4$},\\
&(r_{\lam}, s_{\lam})	\iif \text{$\lam$ is in case $5$};
\end{aligned}
\right.
\end{align*}
\item
for $T \leq \beta_0$ and $\epsilon \in \{\pm1\}$, 
the map $(\beta, \epsilon) \mapsto (\beta + \epsilon/2, \epsilon)$ gives an injection
\[
 \left\{
\begin{aligned}
\CC^\epsilon_{\lam'}(T-1) &\hookrightarrow \CC^\epsilon_{\lam}(T) 
\iif \text{$\lam$ is in case $1$, $2$, $3$ or $4$},\\
\CC^\epsilon_{\lam'}(T) &\hookrightarrow \CC^\epsilon_{\lam}(T) 
\iif \text{$\lam$ is in case $5$};
\end{aligned}
\right.
\]
\item
for fixed $l \geq 0$ and $1 \leq T \leq \beta_0$,  
if $\#\CC^\epsilon_\lam(T) \leq l$ for each $\epsilon \in \{\pm1\}$, then 
\[
\left\{
\begin{aligned}
\#\CC^\epsilon_{\lam'}(T-1) &\leq l \iif \text{$\lam$ is in case $1$},\\
\#\CC^\epsilon_{\lam'}(T-1) &\leq l-1	\iif \text{$\lam$ is in case $2$, $3$ or $4$},\\
\#\CC^\epsilon_{\lam'}(T) &\leq l \iif \text{$\lam$ is in case $5$}
\end{aligned}
\right.
\]
for each $\epsilon \in \{\pm1\}$.
\end{itemize}
\end{lem}

\begin{lem}\label{k>0}
Assume the local Gan--Gross--Prasad conjecture (Conjecture \ref{GGP}).
Let $\lam = (\phi, \eta)$ be a pair of $\phi \in \Phi_\temp(\U_n(\R))$ 
and $\eta \in \widehat{A_\phi}$.
Suppose that $k=k_\lam>0$.
Put $\alpha_+ = \min\{\alpha \ |\ (\alpha, \epsilon) \in X_\lam, \alpha>(k-1)/2\}$ and
$\alpha_- = \max\{\alpha \ |\ (\alpha, \epsilon) \in X_\lam, \alpha<-(k-1)/2\}$.
There is unique $\delta \in \{\pm1\}$ such that
$((k-1)/2, \delta), \dots, ((-k+1)/2, \delta) \in X_\lam$.
We consider the following cases separately. 
\begin{description}
\item[Case $1$]
$\alpha_\delta \not= (k+1)\delta/2$ or $(\alpha_\delta, \delta) \not\in X_\lam$.
\item[Case $2$]
$\alpha_\delta = (k+1)\delta/2$, $(\alpha_\delta, \delta) \in X_\lam$
and $\alpha_{-\delta} = -(k+1)\delta/2$, $(\alpha_{-\delta},-\delta) \in X_\lam$.
\item[Case $3$]
$\alpha_\delta = (k+1)\delta/2$, $(\alpha_\delta, \delta) \in X_\lam$
and $\alpha_{-\delta} \not= -(k+1)\delta/2$, $(\alpha_{-\delta},-\delta) \in X_\lam$.
\item[Case $4$]
$\alpha_\delta = (k+1)\delta/2$, $(\alpha_\delta, \delta) \in X_\lam$
and $\alpha_{-\delta} \not= -(k+1)\delta/2$, $(\alpha_{-\delta},-\delta) \not\in X_\lam$.
\item[Case $5$]
$\alpha_\delta = (k+1)\delta/2$, $(\alpha_\delta, \delta) \in X_\lam$
and $\alpha_{-\delta} = -(k+1)\delta/2$, $(\alpha_{-\delta},-\delta) \not\in X_\lam$.
\end{description}
The case $5$ cannot occur if $\phi \in \Phi_\disc(\U_n(\R))$.
We set $\phi' \in \Phi_\disc(\U_{n+1}(\R))$ so that
\[
\phi'\chi_V^{-1}= \left\{
\begin{aligned}
&\phi_0'\chi_V^{-1} + (\chi_{2\beta_0} + \chi_{2\beta_1} + \dots + \chi_{2\beta_{2v'}})
\iif \text{$\lam$ is in case $1$, $2$ or $5$},\\
&\phi_0'\chi_V^{-1} - (\chi_{2\beta_{+}}+\chi_{2\beta_{-}})
+ (\chi_{2\beta_0} + \chi_{2\beta_1} + \dots + \chi_{2\beta_{2v'+2}})
\iif \text{$\lam$ is in case $3$},\\
&\phi_0'\chi_V^{-1} - \chi_{-k\delta}+\chi_V\chi_{-(k+2)\delta}
+ (\chi_{2\beta_0} + \chi_{2\beta_1} + \dots + \chi_{2\beta_{2v'}})
\iif \text{$\lam$ is in case $4$}.
\end{aligned}
\right.
\]
Here, in the case $3$, we put 
$\beta_{\delta} = -k\delta/2$ and $\beta_{-\delta} = \alpha_{-\delta}+\delta/2$.
Also we define $\eta' \in \widehat{A_{\phi'}}$ by setting 
$\eta'(e_{V,2\beta_j}) = 1$ for $j= 0,1, \ldots$, and 
\[
\eta'(e_{V,2\beta}) = -\eta(e_{V,2\alpha})
\]
when $\beta = \alpha-\epsilon/2$ with $(\alpha,\epsilon) \in X_\lam$.
In the case $4$, we set $\eta'(e_{V, -(k+2)\delta}) = -\eta(e_{V, -(k-1)\delta})$.
Let $\lam' = (\phi',\eta')$ and $\pi'=\pi(\phi',\eta')$.
Then 
\begin{itemize}
\item
$\pi' \in \Irr_\disc(\U(p+1, q))$ and $\Hom_{\U(p,q)}(\pi',\pi) \not= 0$;
\item
$k_{\lam'}$ and $(r_{\lam'}, s_{\lam'})$ are given by
\begin{align*}
k_{\lam'} &= \left\{
\begin{aligned}
&k-1 \iif \text{$\lam$ is in case $1$, $2$, $3$ or $4$},\\
&k+1	\iif \text{$\lam$ is in case $5$},
\end{aligned}
\right.
\\
(r_{\lam'}, s_{\lam'}) &= \left\{
\begin{aligned}
&(r_{\lam}+1, s_{\lam}+1) \iif \text{$\lam$ is in case $1$, $2$, $3$ or $4$},\\
&(r_{\lam}, s_{\lam})	\iif \text{$\lam$ is in case $5$};
\end{aligned}
\right.
\end{align*}
\item
for $T \leq \beta_0+k/2$ and $\epsilon \in \{\pm1\}$, 
the map $(\beta, \epsilon) \mapsto (\beta + \epsilon/2, \epsilon)$ gives an injection
\[
\left\{
\begin{aligned}
\CC^\epsilon_{\lam'}(T-1) &\hookrightarrow \CC^\epsilon_{\lam}(T)
\iif \text{$\lam$ is in case $1$, $2$, $3$ or $4$},\\
\CC^\epsilon_{\lam'}(T) &\hookrightarrow \CC^\epsilon_{\lam}(T)
\iif \text{$\lam$ is in case $5$};
\end{aligned}
\right.
\]
\item
for fixed $l \geq 0$ and $1 \leq T \leq \beta_0 + k/2$, 
if $\#\CC^\epsilon_\lam(T) \leq l$ for each $\epsilon \in \{\pm1\}$, then 
\[
\left\{
\begin{aligned}
\#\CC^\epsilon_{\lam'}(T-1) &\leq l-1	\iif \text{$\lam$ is in case $1$, $2$, $3$ or $4$},\\
\#\CC^\epsilon_{\lam'}(T) &\leq l \iif \text{$\lam$ is in case $5$}
\end{aligned}
\right.
\]
for each $\epsilon \in \{\pm1\}$.
\end{itemize}
\end{lem}

The proofs of Lemmas \ref{k=-1}, \ref{k=0} and \ref{k>0} are straightforward and similar to each other.
So we only prove Lemma \ref{k=-1}.

\begin{proof}[Proof of Lemma \ref{k=-1}]
We check the conditions in Proposition \ref{GGPprop}.
Write
\begin{align*}
\phi\chi_V^{-1} =& 
\chi_{2\alpha_1} + \dots + \chi_{2\alpha_{u}} 
+ (\xi_1 + \dots + \xi_v) + ({}^c\xi_1^{-1} + \dots + {}^c\xi_v^{-1}), 
\end{align*}
where
\begin{itemize}
\item
$\alpha_i \in \half{1}\Z$ such that 
$2\alpha_i \equiv \kappa-1 \bmod 2$ and $\alpha_1 > \dots > \alpha_{u}$;
\item
$\xi_i$ is a unitary character of $\C^\times$ (which can be of the form $\chi_{2\alpha}$);
\item
$u+ 2v=n$.
\end{itemize}
Fix $\alpha$ such that $\chi_V\chi_{2\alpha} \subset \phi$.
Then 
\begin{align*}
\#\{\beta\ |\ \chi_V\chi_{2\beta} \subset \phi',\ \beta < \alpha\}
&\equiv 
\left\{
\begin{aligned}
&\#\{i \in \{1, \dots, u\} \ |\ \alpha_i < \alpha\}+1 &\bmod 2 &\iif (\alpha, +1) \in X_\lam,\\
&\#\{i \in \{1, \dots, u\} \ |\ \alpha_i < \alpha\} &\bmod 2 &\other.
\end{aligned}
\right.
\end{align*}
This implies that 
\[
(-1)^{\#\{\beta\ |\ \chi_V\chi_{2\beta} \subset \phi',\ \beta < \alpha\}} 
= (-1)^n\eta(e_{V,2\alpha}).
\]
Similarly, for each $\beta$ such that $\beta \not= \beta_j$ and $\chi_V\chi_{2\beta} \subset \phi'$, 
we have
\begin{align*}
&\#\{i \in \{1,\dots, u\} \ |\ \alpha_i < \beta\}
\\&
\equiv
\left\{
\begin{aligned}
&\#\{i \in \{1,\dots, u\} \ |\ \alpha_i < \beta +1/2\}+1	&\bmod 2 
&\iif (\beta, -1) \in X_{\lam'}, (\beta-1/2, +1) \not\in X_\lam\\
&\#\{i \in \{1,\dots, u\} \ |\ \alpha_i < \beta -1/2\} &\bmod 2 &\other.
\end{aligned}
\right.
\end{align*}
so that
\[
(-1)^{\#\{i \in \{1,\dots, u\} \ |\ \alpha_i < \beta\}} = (-1)^n\eta'(e_{V,2\beta}). 
\]
This equation also holds when $\beta = \beta_j$ for $j = 0, 1, \dots, 2v'$.
By the local Gan--Gross--Prasad conjecture 
(Conjecture \ref{GGP} and Proposition \ref{GGPprop}), 
we conclude that $\pi' \in \Irr_\disc(\U(p+1,q))$ and $\Hom_{\U(p,q)}(\pi', \pi) \not= 0$.
\par

By the construction, 
if $X_{\lam'}$ contains $(1/2, -1)$ (\resp $(-1/2, +1)$), 
then $X_{\lam}$ must contain $(0,-1)$ (\resp $(0,+1)$).
Since $k_\lam = -1$, in this case, we must have $(0,\pm1) \in X_\lam$, 
so that $X_{\lam'}$ cannot contain $(-1/2, -1)$ (\resp $(1/2,+1)$).
Hence $k_{\lam'} = 0$.
Note that for $(\alpha, \epsilon) \in X_\lam$ with $\alpha \not=0$, 
\[
\epsilon \alpha > 0 \iff \epsilon (\alpha - \half{\epsilon}) >0.
\]
This implies that
when $(0,\pm1) \not\in X_\lam$, we have $(r_{\lam'}, s_{\lam'}) = (r_\lam+1, s_\lam)$.
If $(\beta, \epsilon) = (\pm1/2, \mp1)$, then $\epsilon \beta < 0$.
This implies that
when $(0,\pm1) \in X_\lam$, we have $(r_{\lam'}, s_{\lam'}) = (r_\lam, s_\lam+1)$.
\par

Finally, 
by definition, we have
\[
X_{\lam'} \subset \{(\beta_j, (-1)^j)\ |\ j = 0, \dots, 2v'\} 
\cup \{(\beta, \epsilon)\ |\ (\beta+\epsilon/2, \epsilon) \in X_\lam\}.
\]
Moreover, by construction of $\lam'$, 
we see that if $(\beta, \epsilon) \in X_{\lam'}^{(\infty)}$ with $\beta \not= \beta_0$, 
then $(\alpha, \epsilon) \in X_\lam^{(\infty)}$ with $\alpha = \beta + \epsilon/2$.
In this case, 
\begin{align*}
0 \leq \epsilon\beta-\half{1} < T
&\iff
0 \leq \epsilon\left(\alpha-\half{\epsilon}\right)-\half{1} < T
\\&\iff
0 \leq \epsilon\alpha-1 < T.
\end{align*}
Hence if $(\beta,\epsilon) \in \CC^\epsilon_{\lam'}(T)$ with $\beta \not= \beta_0$, 
then $(\alpha,\epsilon) \in \CC^\epsilon_\lam(T)$ with $\alpha = \beta + \epsilon/2$.
When $T < \beta_0$, we see that 
$(\beta_0,+1)$ is not contained in $\CC^\epsilon_{\lam'}(T)$, so that
we may consider the map 
\[
\CC^\epsilon_{\lam'}(T) \rightarrow \CC^\epsilon_\lam(T),\ 
(\beta, \epsilon) \mapsto \left(\beta+\half{\epsilon}, \epsilon\right).
\] 
This map is clearly injective.
Hence we have $\#\CC^\epsilon_{\lam'}(T) \leq \#\CC^\epsilon_\lam(T)$.
This completes the proof of Lemma \ref{k=-1}.
\end{proof}

\subsection{Non-vanishing}\label{nonvanish}
In this subsection, we prove sufficient conditions
of the non-vanishing of theta lifts in Theorem \ref{main}.
\par

Let $\lam = (\phi, \eta)$ be a pair of $\phi \in \Phi_\temp(\U_n(\R))$ and 
$\eta \in \widehat{A_\phi}$.
Set $k=k_\lam$, $r=r_\lam$, $s=s_\lam$, and $\pi = \pi(\phi,\eta)$.
For non-negative integers $t$ and $l$, 
consider the following statements:
\begin{description}
\item[$(S)_{-1,t,l}$]
Suppose that $k_\lam = -1$.
If $\#\CC^\epsilon_\lam(t+l) \leq l$ for each $\epsilon \in \{\pm1\}$, 
then 
\[
\left\{
\begin{aligned}
&\Theta_{r+2t+l+1,s+l}(\pi) \not=0 \iif (0,\pm1) \not\in X_\lam,\\
&\Theta_{r+2t+l,s+l+1}(\pi) \not=0 \iif (0,\pm1) \in X_\lam.
\end{aligned}
\right.
\]
\item[$(S)_{k,t,l}$]
Suppose that $k_\lam=k \geq 0$ and $l \geq k$.
If $\#\CC^\epsilon_\lam(t+l) \leq l$ for each $\epsilon \in \{\pm1\}$, 
then $\Theta_{r+2t+l,s+l}(\pi) \not=0$.
\end{description}
\par

First, we consider these statements for the discrete series representations.
We have implications.
\begin{prop}\label{red}
Consider the statement $(S)_{k,t,l}$ only when $\phi \in \Phi_\disc(\U_n(\R))$.
\begin{enumerate}
\item
For $t \geq 0$ and $l \geq 0$, 
we have $(S)_{0,t,l} \Rightarrow (S)_{-1,t,l}$.
\item
For $t \geq 1$ and $l \geq 0$, 
we have $(S)_{-1,t-1,l} + (S)_{-1,t,l-1} \Rightarrow (S)_{0,t,l}$.
\item
For $t \geq 0$ and $l \geq k >0$, we have $(S)_{k-1,t,l-1} \Rightarrow (S)_{k,t,l}$.
\end{enumerate}
Here, we interpret $(S)_{k,t,-1}$ to be empty. 
\end{prop}
\begin{proof}
Suppose that $\lam = (\phi, \eta)$ is a pair of 
$\phi \in \Phi_\disc(\U_n(\R))$ and $\eta \in \widehat{A_\phi}$.
Let $\lam' = (\phi',\eta')$ be as in Lemma \ref{k=-1}, \ref{k=0} or \ref{k>0}.
Here, we take $\beta_0$ so that $\beta_0+k_\lam/2 \geq t+l$.
Since $\phi \in \Phi_{\disc}(\U_n(\R))$ and $\phi' \in \Phi_{\disc}(\U_{n+1}(\R))$, 
the local Gan--Gross--Prasad conjecture for $\phi$ and $\phi'$
has been established by He \cite{He}.
So we have $\Hom_{\U(p,q)}(\pi',\pi) \not= 0$ unconditionally.
\par

We show (1).
Suppose that $k_{\lam} = -1$.
By Lamma \ref{k=-1}, we have 
$k_{\lam'} = 0$, $(r_{\lam'}, s_{\lam'}) = (r_\lam+1,s_\lam)$, 
and if $\#\CC_\lam^\epsilon(t+l) \leq l$, then $\#\CC^\epsilon_{\lam'}(t+l) \leq l$.
Hence we can apply $(S)_{0,t,l}$ to $\lam'$, and we obtain that
\[
\Theta_{(r_\lam+1)+2t+l, s_\lam+l}(\pi') \not=0.
\]
Since $\Hom_{\U(p,q)}(\pi', \pi) \not= 0$, 
by the seesaw
\[
\xymatrix{
\U(p+1, q) \ar@{-}[d] \ar@{-}[dr] & 
\U(r_\lam+1+2t+l, s_\lam+l) \times \U(r_\lam+1+2t+l, s_\lam+l) \ar@{-}[d]\\
\U(p,q) \times \U(1,0) \ar@{-}[ur] & \U(r_\lam+1+2t+l, s_\lam+l),
}
\]
we conclude that $\Theta_{r_\lam+1+2t+l,s_\lam+l}(\pi) \not= 0$. 
Therefore, we have $(S)_{0,t,l} \Rightarrow (S)_{-1,t,l}$.
\par

The proofs of (2) and (3) are similar. 
Note that the cases $5$ in Lemmas \ref{k=0} and \ref{k>0} cannot occur
since $\phi \in \Phi_\disc(\U_n(\R))$. 
We omit the detail.
\end{proof}

\begin{cor}\label{true}
The statement $(S)_{k,t,l}$ is true 
for $\lam = (\phi, \eta)$ such that $\phi \in \Phi_\disc(\U_n(\R))$.
\end{cor}
\begin{proof}
When $k > 0$, 
the statement $(S)_{k,t,l}$ is reduced to $(S)_{0,t,l-k}$ by Proposition \ref{red} (3).
We prove $(S)_{k,t,l}$ for $k \leq 0$ by induction on $t+l$.
For $k=-1, 0$ and $T \geq 0$, we consider the following statement:
\begin{description}
\item[$(S')_{k,T}$]
the statement $(S)_{k,t,l}$ is true for any $t,l \geq 0$ such that $t+l \leq T$.
\end{description}
Note that $(S)_{0,0,0}$ is Paul's result (Theorem \ref{PE}).
This implies $(S)_{-1,0,0}$ (Proposition \ref{red} (1)).
In particular, $(S')_{k,0}$ is true.
Also the tower property (Proposition \ref{tower}) implies
$(S)_{k,0,l}$ for $k=-1,0$ and $l \geq 0$.
By Proposition \ref{red} (1) and (2), we have
\[
(S')_{-1,T-1} \Rightarrow (S')_{0,T} \Rightarrow (S')_{-1,T}
\]
for any $T > 0$.
Hence by induction, we obtain $(S')_{k,T}$ for $k =-1, 0$ and $T \geq 0$.
\end{proof}
\par

Now we obtain the sufficient conditions of the non-vanishing of theta lifts.
\begin{cor}\label{suf}
Assume the local Gan--Gross--Prasad conjecture (Conjecture \ref{GGP}).
Then the statement $(S)_{k,t,l}$ is true in general.
\end{cor}
\begin{proof}
The statement $(S)_{0,0,l}$ follows from Paul's result (Theorem \ref{PE}) and 
the tower property (Proposition \ref{tower}).
By a similar argument to the proof of Proposition \ref{red}, 
the other assertions follow from Corollary \ref{true}
by using a seesaw identity (Proposition \ref{seesaw} (1))
and Lemmas \ref{k=-1}, \ref{k=0} and \ref{k>0}.
We omit the detail.
\end{proof}
\par

Let $\lam = (\phi, \eta)$ be a pair of 
$\phi \in \Phi_\temp(\U_n(\R))$ and $\eta \in \widehat{A_\phi}$.
Set $k=k_\lam$, $r=r_\lam$, $s=s_\lam$, and $\pi = \pi(\phi,\eta)$.
Recall as in Theorem \ref{PA} (3) that 
when $k = -1$ but $\phi$ contains $\chi_V$ (with even multiplicity), 
exactly one of $\Theta_{r-1,s}(\pi)$ and $\Theta_{r,s-1}(\pi)$ is nonzero.
Now we can determine which is nonzero in terms of $\lam$.

\begin{cor}\label{00}
Assume the local Gan--Gross--Prasad conjecture (Conjecture \ref{GGP}).
Let $\lam = (\phi, \eta)$ be a pair of 
$\phi \in \Phi_\temp(\U_n(\R))$ and $\eta \in \widehat{A_\phi}$.
Set $k=k_\lam$, $r=r_\lam$, $s=s_\lam$, and $\pi = \pi(\phi,\eta)$.
Suppose that $k = -1$ but $\phi$ contains $\chi_V$ (with even multiplicity).
Then $\Theta_{r,s-1}(\pi) \not= 0$ if and only if $(0,\pm1) \not\in X_\lam$.
\end{cor}
\begin{proof}
Note that $r+s=p+q$.
By a result of Paul (Theorem \ref{PA} (3)) and the conservation relation (Theorem \ref{CR}), 
we see that exactly one of $\Theta_{r+1,s}(\pi)$ and $\Theta_{r,s+1}(\pi)$ is nonzero.
Hence $\Theta_{r,s-1}(\pi) \not= 0$ if and only if $\Theta_{r+1,s}(\pi) \not= 0$.
By the statement $(S)_{-1,0,0}$, we have 
\[
\left\{
\begin{aligned}
&\Theta_{r+1, s}(\pi) \not= 0 \iif (0,\pm1) \not\in X_\lam, \\
&\Theta_{r, s+1}(\pi) \not= 0 \iif (0,\pm1) \in X_\lam. 
\end{aligned}
\right.
\]
This completes the proof.
\end{proof}

\subsection{Vanishing}\label{vanish}
In this subsection, 
we prove necessary conditions of non-vanishing of theta lifts in Theorem \ref{main}.
First, we prove the following:

\begin{prop}\label{zero}
Let $\lam = (\phi, \eta)$ be a pair of 
$\phi \in \Phi_\temp(\U_n(\R))$ and $\eta \in \widehat{A_\phi}$.
Set $\pi = \pi(\phi,\eta)$.
Let $k=k_\lam$, $r=r_\lam$, $s=s_\lam$, $X_\lam$ and $\CC^\epsilon_\lam(T)$ 
be as in Definition \ref{def}.
\begin{enumerate}
\item
Suppose that $k=-1$.
For $t \geq 1$ and $l \geq 0$, 
if $\CC^+_\lam(t+l) > l$ or $\CC^-_\lam(t+l) > l$, 
then $\Theta_{r+2t+l+1, s+l}(\pi) =0$.
\item
Suppose that $k \geq 0$.
For $t \geq 1$ and $l \geq k$, 
if $\CC^+_\lam(t+l) > l$ or $\CC^-_\lam(t+l) > l$, 
then $\Theta_{r+2t+l, s+l}(\pi) =0$.
\end{enumerate}
\end{prop}
\begin{proof}
The proof is similar to that of \cite[Theorem 3.14]{P2}.
Suppose for the sake of contradiction that 
for some $t \geq 1$ and $l \geq \max\{0,k\}$, 
\[
\left\{
\begin{aligned}
&\Theta_{r+2t+l, s+l}(\pi) \not= 0 \iif k \geq 0,\\
&\Theta_{r+1+2t+l, s+l}(\pi) \not= 0 \iif k=-1,
\end{aligned}
\right.
\]
but there exists $\epsilon \in \{\pm1\}$ such that $\CC^\epsilon_\lam(t+l) > l$.
Write 
\[
\CC^\epsilon_\lam(t+l) 
= \{(\alpha_1, \epsilon), (\alpha_2, \epsilon), \dots, (\alpha_{l+1}, \epsilon), \dots\}
\]
with $\epsilon\alpha_1 < \dots < \epsilon\alpha_{l+1} < \dots < t+l - (k-1)/2$.
We set
\[
a = \left\{
\begin{aligned}
&\alpha_{l+1}+\epsilon/2 \iif \text{$k$ is even},\\
&\alpha_{l+1}	\iif \text{$k$ is odd}.
\end{aligned}
\right.
\]
Then $a$ is an integer.
By the definition of $\CC^\epsilon_\lam(t+l)$ and an easy calculation, 
we have
\begin{align*}
&\#\{(\alpha, \epsilon) \in X_\lam\ |\ \epsilon\alpha > 0,\ \epsilon\alpha < \epsilon a\}
-
\#\{(\alpha, -\epsilon) \in X_\lam\ |\ -\epsilon\alpha < 0,\ -\epsilon\alpha > -\epsilon a\}
\\&
= \left\{
\begin{aligned}
&l+1-\half{k} \iif \text{$k$ is even},\\
&l-\half{k+1} \iif \text{$k$ is odd, $k>0$ and $(0, \epsilon) \in X_\lam$},\\ 
&l-\half{k-1} \iif \text{$k$ is odd, $k>0$ and $(0, -\epsilon) \in X_\lam$},\\
&l \iif \text{$k=-1$ and $(0,\pm1) \not\in X_\lam$},\\
&l+1 \iif \text{$k=-1$ and $(0,\pm1) \in X_\lam$}.
\end{aligned}
\right.
\end{align*}
Also, since $l \geq \max\{0,k\}$, we have $\epsilon a > 0$.
Hence
\begin{align*}
\{(\alpha, -\epsilon) \in X_\lam\ |\ -\epsilon\alpha > 0,\ -\epsilon\alpha < -\epsilon a\}
=
\{(\alpha, \epsilon) \in X_\lam\ |\ \epsilon\alpha < 0,\ \epsilon\alpha > \epsilon a\}
= \emptyset.
\end{align*}
\par

We set
\[
\delta = \left\{
\begin{aligned}
&1 \iif \text{$k$ is odd, $k>0$ and $(\alpha_{l+1}, -\epsilon) \in X_\lam$},\\
&0 \other,
\end{aligned}
\right.
\]
and
\begin{align*}
x &= \#\{(\alpha, \epsilon') \in X_\lam\ |\ \epsilon'\alpha > 0,\ \epsilon'\alpha > \epsilon' a\},\\
y &= \#\{(\alpha, \epsilon') \in X_\lam\ |\ \epsilon'\alpha > 0,\ \epsilon'\alpha < \epsilon' a\},\\
z &= \#\{(\alpha, \epsilon') \in X_\lam\ |\ \epsilon'\alpha < 0,\ \epsilon'\alpha > \epsilon' a\},\\
w &= \#\{(\alpha, \epsilon') \in X_\lam\ |\ \epsilon'\alpha < 0,\ \epsilon'\alpha < \epsilon' a\}, 
\end{align*}
and $v = (n-\#X_\lam)/2$.
Then by Definition \ref{def}, we have
\[
(x+y+v, z+w+v) = \left\{
\begin{aligned}
&(r + \half{k}, s + \half{k}) \iif \text{$k$ is even},\\
&(r + \half{k-1}-1, s + \half{k-1}-\delta) \iif \text{$k$ is odd and $k>0$},\\
&(r-1,s-\delta) \iif \text{$k=-1$ and $(0,\pm1) \not\in X_\lam$},\\
&(r-2,s-1-\delta) \iif \text{$k=-1$ and $(0,\pm1) \in X_\lam$}.
\end{aligned}
\right.
\]
On the other hand, by the above calculation, we have
\[
y-z = \left\{
\begin{aligned}
&l+1-\half{k} \iif \text{$k$ is even},\\
&l-\half{k+1} \iif \text{$k$ is odd, $k>0$ and $(0, \epsilon) \in X_\lam$},\\ 
&l-\half{k-1} \iif \text{$k$ is odd, $k>0$ and $(0, -\epsilon) \in X_\lam$},\\
&l \iif \text{$k=-1$ and $(0,\pm1) \not\in X_\lam$},\\
&l+1 \iif \text{$k=-1$ and $(0,\pm1) \in X_\lam$}.
\end{aligned}
\right.
\]
\par

Now we consider $\lam_a = (\phi \otimes \chi_{-2a}, \eta)$ 
and $\pi(\phi \otimes \chi_{-2a}, \eta) = \pi \otimes {\det}^{-a}$.
Note that $k_{\lam_a} \equiv k \bmod 2$.
There exists a bijection 
\[
X_\lam \rightarrow X_{\lam_a},\ 
(\alpha, \epsilon') \mapsto (\alpha-a, \epsilon').
\]
We set
\begin{align*}
r' &= \#\{(\alpha',\epsilon') \in X_{\lam_a}\ |\ \epsilon'\alpha' > 0\} + v +\delta,\\
s' &= \#\{(\alpha',\epsilon') \in X_{\lam_a}\ |\ \epsilon'\alpha' < 0\} + v +\delta.
\end{align*}
Note that
\[
r'+s' = \left\{
\begin{aligned}
&n \iif \text{$k$ is even},\\
&n-1+\delta \iif \text{$k$ is odd}.
\end{aligned}
\right.
\]
By the above bijection, we have 
\begin{align*}
r' &= \left\{
\begin{aligned}
&x+z+v \iif \text{$k$ is even},\\
&x+z+v+\delta \iif \text{$k>0$, $k$ is odd and $(0,\epsilon) \in X_\lam$}, \\
&x+z+v+\delta+1 \iif \text{$k>0$, $k$ is odd and $(0,-\epsilon) \in X_\lam$}, \\
&x+z+v+\delta \iif \text{$k=-1$ and $(0,\pm1) \not\in X_\lam$}, \\
&x+z+v+\delta+1 \iif \text{$k=-1$ and $(0,\pm1) \in X_\lam$}
\end{aligned}
\right.
\\
&= \left\{
\begin{aligned}
&r+k-l-1 \iif \text{$k$ is even},\\
&r+k-l-1+\delta \iif \text{$k>0$ and $k$ is odd}, \\
&r-l-1+\delta \iif \text{$k=-1$ and $(0,\pm1) \not\in X_\lam$}, \\
&r-l-2+\delta \iif \text{$k=-1$ and $(0,\pm1) \in X_\lam$}
\end{aligned}
\right.
\end{align*}
and
\begin{align*}
s' &= \left\{
\begin{aligned}
&y+w+v \iif \text{$k$ is even},\\
&y+w+v+\delta+1 \iif \text{$k>0$, $k$ is odd and $(0,\epsilon) \in X_\lam$}, \\
&y+w+v+\delta \iif \text{$k>0$, $k$ is odd and $(0,-\epsilon) \in X_\lam$}, \\
&y+w+v+\delta \iif \text{$k=-1$ and $(0,\pm1) \not\in X_\lam$}, \\
&y+w+v+\delta+1 \iif \text{$k=-1$ and $(0,\pm1) \in X_\lam$}
\end{aligned}
\right.
\\
&= \left\{
\begin{aligned}
&s+l+1 \iif \text{$k$ is even},\\
&s+l \iif \text{$k>0$ and $k$ is odd}, \\
&s+l \iif \text{$k=-1$ and $(0,\pm1) \not\in X_\lam$}, \\
&s+l+1 \iif \text{$k=-1$ and $(0,\pm1) \in X_\lam$}.
\end{aligned}
\right.
\end{align*}
By Theorems \ref{PE}, \ref{PA} (1) and Corollary \ref{00}, 
we see that 
\[
\left\{
\begin{aligned}
&\Theta_{r'-1,s'}(\pi \otimes {\det}^{-a}) \not= 0 \iif \delta = 1,\\
&\Theta_{r',s'}(\pi \otimes {\det}^{-a}) \not= 0 \other, 
\end{aligned}
\right.
\]
i.e., $\Theta_{r'-\delta,s'}(\pi \otimes {\det}^{-a}) \not= 0$.
Hence $\Theta_{s',r'-\delta}(\pi^\vee \otimes {\det}^{a} \otimes \chi_V^2) \not= 0$ 
by Proposition \ref{vee}.
By the seesaw (Proposition \ref{seesaw} (2))
\[
\xymatrix{
\U(p, q) \times \U(p,q) \ar@{-}[d] \ar@{-}[dr] & 
\U((r+2t+l)+s', (s+l)+(r'-\delta)) \ar@{-}[d]\\
\U(p,q) \ar@{-}[ur] & \U(r+2t+l,s+l) \times \U(s', r'-\delta), 
}
\]
we deduce that
\[
\left\{
\begin{aligned}
&\Theta_{(r+2t+l)+s', (s+l)+(r'-\delta)}({\det}^{a} \cdot \chi_{V_{0,0}}) \not= 0 \iif k\geq0,\\
&\Theta_{(r+1+2t+l)+s', (s+l)+(r'-\delta)}({\det}^{a} \cdot \chi_{V_{0,0}}) \not= 0 \iif k=-1.
\end{aligned}
\right.
\]
Here, $\chi_{V_{0,,0}} = \chi_{V_{(r+2t+l)+s', (s+l)+(r'-\delta)}}$ if $k \geq 0$, 
and  $\chi_{V_{0,,0}} = \chi_{V_{(r+1+2t+l)+s', (s+l)+(r'-\delta)}}$ if $k=-1$.
\par

On the other hand, since 
\[
r+s = 
\left\{
\begin{aligned}
&n-k \iif k \geq 0,\\
&n	\iif k=-1,
\end{aligned}
\right.
\]
we have
\begin{align*}
(r+2t+l)+s' &= \left\{
\begin{aligned}
&(n-1)+2t+2l+2-k \iif \text{$k$ is even},\\
&(n-1)+2t+2l-(k-1)  \iif \text{$k>0$ and $k$ is odd},\\
&(n-1)+2t+2l+1 \iif \text{$k=-1$ and $(0,\pm1) \not\in X_\lam$},\\
&(n-2)+2t+2l+3 \iif \text{$k=-1$ and $(0,\pm1) \in X_\lam$},\\
\end{aligned}
\right.
\\
(s+l)+(r'-\delta) &= \left\{
\begin{aligned}
&n-2 \iif \text{$k=-1$ and $(0,\pm1) \in X_\lam$},\\
&n-1 \other.
\end{aligned}
\right.
\end{align*}
In particular, $\min\{(r+2t+l)+s', (s+l)+(r'-\delta)\} = (s+l)+(r'-\delta) < n$ and
\[
((r+2t+l)+s')-((s+l)+(r'-\delta)) 
= \left\{
\begin{aligned}
&2t+2l+2-k \iif \text{$k$ is even},\\
&2t+2l-(k-1)  \iif \text{$k>0$ and $k$ is odd},\\
&2t+2l+1 \iif \text{$k=-1$ and $(0,\pm1) \not\in X_\lam$},\\
&2t+2l+3 \iif \text{$k=-1$ and $(0,\pm1) \in X_\lam$}.
\end{aligned}
\right.
\]
Moreover, since $\epsilon\alpha_{l+1}+(k-1)/2 < t+l$ and $\epsilon a> 0$, we have
\[
\left\{
\begin{aligned}
&0 < \epsilon a < t+l+1-\half{k} \iif \text{$k$ is even},\\
&0 < \epsilon a < t+l-\half{k-1}  \iif \text{$k$ is odd}.
\end{aligned}
\right.
\]
By Proposition \ref{trivial}, we must have
\[
\left\{
\begin{aligned}
&\Theta_{(r+2t+l)+s', (s+l)+(r'-\delta)}({\det}^{a} \cdot \chi_{V_{0,0}}) = 0 \iif k\geq0,\\
&\Theta_{(r+1+2t+l)+s', (s+l)+(r'-\delta)}({\det}^{a} \cdot \chi_{V_{0,0}}) = 0 \iif k=-1.
\end{aligned}
\right.
\]
We obtain a contradiction.
\end{proof}

By a similar argument, we obtain the following.
\begin{prop}\label{lk}
Let $\lam = (\phi, \eta)$ be a pair of 
$\phi \in \Phi_\temp(\U_n(\R))$ and $\eta \in \widehat{A_\phi}$.
Set $\pi = \pi(\phi,\eta)$.
Let $k=k_\lam$, $r=r_\lam$ and $s=s_\lam$ be as in Definition \ref{def}.
\begin{enumerate}
\item
Suppose that $k=-1$.
If $\Theta_{r+2t+l+1, s+l}(\pi) \not= 0$ for some $t \geq 1$, then $l \geq 0$.
\item
Suppose that $k \geq 0$.
If $\Theta_{r+2t+l, s+l}(\pi) \not= 0$ for some $t \geq 1$, then $l \geq k$. 
\end{enumerate}
\end{prop}
\begin{proof}
We give an outline of the proof.
Suppose that 
\[
\left\{
\begin{aligned}
&\Theta_{r+2t+l, s+l}(\pi) \not= 0 \iif k \geq 0,\\
&\Theta_{r+1+2t+l, s+l}(\pi) \not= 0 \iif k=-1
\end{aligned}
\right.
\]
for some $t \geq 1$.
Set 
\[
a = \left\{
\begin{aligned}
& \half{k}\epsilon \iif \text{$k$ is even, $k>0$ and $(\half{k-1}, \epsilon) \in X_\lam$}, \\
& \half{k-1}\epsilon \iif \text{$k$ is odd, $k>0$ and $(\half{k-1}, \epsilon) \in X_\lam$}, \\
& 0 \iif \text{$k=-1$ or $k=0$}.
\end{aligned}
\right.
\]
By Theorems \ref{PE}, \ref{PA} (1) and Corollary \ref{00}, 
we have
\[
\left\{
\begin{aligned}
&\Theta_{r, s+k}(\pi \otimes {\det}^{-a}) \not = 0 
\iif \text{$k$ is even}, \\
&\Theta_{r, s+k-1}(\pi \otimes {\det}^{-a}) \not = 0 
\iif \text{$k$ is odd and $k > 0$}, \\
&\Theta_{r, s-1}(\pi \otimes {\det}^{-a}) \not = 0 
\iif \text{$k=-1$ and $(0,\pm1) \not\in X_\lam$ but $\chi_V \subset \phi$}, \\
&\Theta_{r, s+1}(\pi \otimes {\det}^{-a}) \not = 0 
\iif \text{$k=-1$ and $(0,\pm1) \in X_\lam$ if $\chi_V \subset \phi$}.
\end{aligned}
\right.
\]
By a similar argument to the proof of Proposition \ref{zero}, 
a seesaw identity (Proposition \ref{seesaw} (2)) and Proposition \ref{trivial}
imply that $l \geq \max\{0,k\}$.
We omit the detail.
\end{proof}

By Propositions \ref{zero} and \ref{lk}, we have the necessary conditions.
\begin{cor}\label{nec}
Let $\lam = (\phi, \eta)$ be a pair of 
$\phi \in \Phi_\temp(\U_n(\R))$ and $\eta \in \widehat{A_\phi}$.
Set $\pi = \pi(\phi,\eta)$.
Let $k=k_\lam$, $r=r_\lam$, $s=s_\lam$, $X_\lam$ and $\CC^\epsilon_\lam(T)$ 
be as in Definition \ref{def}.
Let $t$ be a positive integer.
\begin{enumerate}
\item
Suppose that $k=-1$.
If $\Theta_{r+2t+l+1, s+l}(\pi) \not= 0$, then
$l \geq 0$ and $\CC^\epsilon_\lam(t+l) \leq l$ for each $\epsilon \in \{\pm1\}$.
\item
Suppose that $k \geq 0$.
If $\Theta_{r+2t+l, s+l}(\pi) \not= 0$, then
$l \geq k$ and $\CC^\epsilon_\lam(t+l) \leq l$ for each $\epsilon \in \{\pm1\}$.
\end{enumerate}
\end{cor}

\subsection{Going-down towers}\label{t=0}
By Corollaries \ref{suf} and \ref{nec}, 
we can determine the first occurrence indices $m_d(\pi)$ of
the $d$-th Witt tower of theta lifts of $\pi = \pi(\phi,\eta)$
when $d - (r_\lam-s_\lam) > 1$ with $\lam = (\phi,\eta)$.
By Lemma \ref{Xv}, we can also determine $m_d(\pi)$ when $d -(r_\lam-s_\lam) < -1$.
In particular, if $|d -(r_\lam-s_\lam)| > 1$, then
\[
m_d(\pi) \geq n+2.
\]
In this case, we call the $d$-th Witt tower a going-up tower with respect to $\pi$.
When $|d -(r_\lam-s_\lam)| \leq 1$, 
we call the $d$-th Witt tower a going-down tower with respect to $\pi$.
By the conservation relation, 
we can determine the first occurrence indices of the going-down Witt towers.

\begin{prop}\label{down}
Assume the local Gan--Gross--Prasad conjecture (Conjecture \ref{GGP}).
Let $\lam = (\phi, \eta)$ be a pair of 
$\phi \in \Phi_\temp(\U_n(\R))$ and $\eta \in \widehat{A_\phi}$.
Set $\pi = \pi(\phi,\eta)$.
Let $k=k_\lam$, $r=r_\lam$, $s=s_\lam$ and $X_\lam$ be as in Definition \ref{def}.
\begin{enumerate}
\item
Suppose that $k=-1$.
Then for an integer $l$, we have
\[
\Theta_{r+1+l,s+l}(\pi) \not= 0 \iff
\left\{
\begin{aligned}
&l \geq 0 \iif \text{$\phi$ does not contain $\chi_V$},\\
&l \geq -1 \iif \text{$\phi$ contains $\chi_V$ but $(0,\pm1) \not \in X_\lam$},\\
&l \geq 1 \iif \text{$\phi$ contains $\chi_V$ and $(0,\pm1) \in X_\lam$}, 
\end{aligned}
\right.
\]
and
\[
\Theta_{r+l,s+1+l}(\pi) \not= 0 \iff
\left\{
\begin{aligned}
&l \geq 0 \iif \text{$\phi$ does not contain $\chi_V$},\\
&l \geq 1 \iif \text{$\phi$ contains $\chi_V$ but $(0,\pm1) \not \in X_\lam$},\\
&l \geq -1 \iif \text{$\phi$ contains $\chi_V$ and $(0,\pm1) \in X_\lam$}.
\end{aligned}
\right.
\]

\item
Suppose that $k \geq 0$.
Consider the following three conditions on $\lam = (\phi, \eta)$:
\begin{description}
\item[(chain condition 2)]
$\phi\chi_V^{-1}$ contains both $\chi_{k+1}$ and $\chi_{-(k+1)}$, 
so that 
\[
\phi \chi_V^{-1} \supset 
\underbrace{\chi_{k+1} + \chi_{k-1} + \dots + \chi_{-(k-1)} + \chi_{-(k+1)}}_{k+2};
\]
\item[(even-ness condition)]
at least one of $\chi_{k+1}$ and $\chi_{-(k+1)}$ is 
contained in $\phi\chi_V^{-1}$ with even multiplicity;
\item[(alternating condition 2)]
$\eta(e_{V, k+1-2i}) \not= \eta(e_{V, k-1-2i})$ for $i=0, \dots, k$.
\end{description}
Then for an integer $l$, we have
\[
\Theta_{r+l,s+l}(\pi) \not= 0 \iff
\left\{
\begin{aligned}
&l \geq -1	\iif \text{$\lam$ satisfies these three conditions}, \\
&l \geq 0	\other.
\end{aligned}
\right.
\]
\end{enumerate}
\end{prop}

When $\phi \in \Phi_\disc(\U_n(\R))$, 
this proposition has been proven by Paul \cite[Proposition 3.4]{P2}
in terms of Harish-Chandra parameters
(not using the Gan--Gross--Prasad conjecture).

\begin{proof}[Proof of Proposition \ref{down}]
We show (1).
Suppose that $k=-1$.
If $\phi$ does not contain $\chi_V$, 
then by Theorem \ref{PA} (1), we have
$\Theta_{r+1,s}(\pi) \not= 0$ and $\Theta_{r,s+1}(\pi) \not= 0$.
By the tower property (Proposition \ref{tower}) and the conservation relation (Theorem \ref{CR}), 
we see that $\Theta_{r+1+l,s+l}(\pi) \not= 0 \iff l \geq 0$, and
$\Theta_{r+l,s+1+l}(\pi) \not= 0 \iff l \geq 0$.
\par

Now we assume that $\phi$ contains $\chi_V$ (with even multiplicity).
Then by Corollary \ref{00}, we have
\[
\left\{
\begin{aligned}
&\Theta_{r,s-1}(\pi) \not=0 \iif (0,\pm1) \not\in X_\lam,\\
&\Theta_{r-1,s}(\pi) \not=0 \iif (0,\pm1) \in X_\lam.
\end{aligned}
\right.
\]
When $(0,\pm1) \in X_\lam$, by Corollary \ref{suf}, 
$\Theta_{r+2, s+1}(\pi) \not=0$ 
if $\CC^\epsilon_\lam(1) = \emptyset$ for each $\epsilon \in \{\pm1\}$.
By Definition \ref{def} (5), 
$\CC^\epsilon_\lam(1) = \emptyset$ for each $\epsilon \in \{\pm1\}$ 
if and only if $(\pm1,\pm1) \not\in X_\lam^{(\infty)}$.
However, since $(0,\pm1) \in X_\lam$, by Definition \ref{def} (4), 
we see that $X_\lam^{(\infty)}$ cannot contain $(\pm1,\pm1)$.
Hence we deduce that $\Theta_{r+2, s+1}(\pi) \not=0$.
By the tower property (Proposition \ref{tower}) and the conservation relation (Theorem \ref{CR}), 
we see that $\Theta_{r+1+l,s+l}(\pi) \not= 0 \iff l \geq 1$, and
$\Theta_{r+l,s+1+l}(\pi) \not= 0 \iff l \geq -1$.
\par

Now suppose that $\phi$ contains $\chi_V$ (with even multiplicity) but $(0,\pm1) \not\in X_\lam$.
Set $\lam^\vee = (\phi^\vee \otimes \chi_V^2, \eta^\vee)$.
By Lemma \ref{Xv}, we have $(0,\pm1) \in X_{\lam^\vee}$ so that 
$\Theta_{s+2,r+1}(\pi^\vee \otimes \chi_V^2) \not= 0$ by the above case.
By Proposition \ref{vee}, we deduce that $\Theta_{r+1,s+2}(\pi) \not= 0$.
By the tower property (Proposition \ref{tower}) and the conservation relation (Theorem \ref{CR}), 
we see that $\Theta_{r+1+l,s+l}(\pi) \not= 0 \iff l \geq -1$, and
$\Theta_{r+l,s+1+l}(\pi) \not= 0 \iff l \geq 1$.
This completes the proof of (1).
\par

We show (2). 
Suppose that $k \geq 0$.
By Corollaries \ref{suf}, \ref{nec} and Proposition \ref{vee}, we see that
\[
m_d(\pi) \geq r+s+2k+2=n+k+2
\]
for any integer $d$ such that $d \not= r-s$.
Hence $\min\{m_+(\pi), m_-(\pi)\} = m_{r-s}(\pi)$.
Moreover, if $|d-(r-s)| > 2$, then $m_d(\pi) \geq n+k+4$.
We compute $m_{r-s+2}(\pi)$ and $m_{r-s-2}(\pi)$.
By Corollaries \ref{suf} and \ref{nec}, 
$\Theta_{r+2+k, s+k}(\pi) \not=0$ if and only if 
$\#\CC^\epsilon_\lam(1+k) \leq k$ for each $\epsilon \in \{\pm1\}$.
\par

First, we consider the case where $k=0$.
Then by Definition \ref{def} (5), 
$\CC^\epsilon_\lam(1) = \emptyset$ for each $\epsilon \in \{\pm1\}$ 
if and only if $(\pm1/2,\pm1) \not\in X_{\lam}^{(\infty)}$.
Also by Definition \ref{def} (4), 
$(\pm1/2,\pm1) \in X_{\lam}^{(\infty)}$ if and only if $(\pm1/2,\pm1) \in X_{\lam}$.
Hence
\[
\Theta_{r+2, s}(\pi) \not=0 \iff \left(\pm\half{1},\pm1\right) \not\in X_{\lam}.
\]
Similarly, by using Proposition \ref{vee}, we see that
\[
\Theta_{r, s+2}(\pi) \not=0 \iff \left(\pm\half{1},\pm1\right) \not\in X_{\lam^\vee}.
\]
Now we assume that both $\Theta_{r+2, s}(\pi)$ and $\Theta_{r, s+2}(\pi)$ are zero.
By Lemma \ref{Xv} (4), 
this condition is equivalent that 
$(1/2, +1) \in X_\lam \cap X_{\lam^\vee}$ or $(-1/2, -1) \in X_\lam \cap X_{\lam^\vee}$. 
We check (chain condition 2) and (even-ness condition).
For $\epsilon \in \{\pm1\}$, 
if $\chi_{-\epsilon}$ were not contained in $\phi\chi_V^{-1}$, 
then we must have $(\epsilon/2, \epsilon) \in X_\lam$ 
and $(-\epsilon/2,-\epsilon) \in X_{\lam^\vee}$.
This contradicts Lemma \ref{Xv} (4).
Hence $\phi\chi_V^{-1}$ contains both $\chi_{1}$ and $\chi_{-1}$.
If both $\chi_{1}$ and $\chi_{-1}$ were contained in $\phi\chi_V^{-1}$
with odd multiplicities, then by Lemma \ref{Xv} (2), 
there must be $\epsilon \in \{\pm1\}$ such that 
$(\epsilon/2, \epsilon) \in X_\lam \cap X_{\lam^\vee}$.
This implies that $(1/2, \epsilon), (-1/2, \epsilon) \in X_\lam$, 
which contradicts that $k_\lam=0$ (see Definition \ref{def} (1)).
\par

We claim that under assuming (chain condition 2) and (even-ness condition), 
both $\Theta_{r+2, s}(\pi)$ and $\Theta_{r, s+2}(\pi)$ are zero
if and only if 
$\lam$ satisfies (alternating condition 2), which is equivalent that
$\eta(e_{V, 1}) \not= \eta(e_{V, -1})$.
Replacing $\lam$ with $\lam^\vee$ if necessary, 
we may assume that $\chi_{1}$ appears in $\phi\chi_V^{-1}$ with even multiplicity.
Write
\[
\phi\chi_V^{-1} =
\chi_{2\alpha_1} + \dots + \chi_{2\alpha_{u}} 
+ (\xi_1 + \dots + \xi_v) + ({}^c\xi_1^{-1} + \dots + {}^c\xi_v^{-1}), 
\]
where
\begin{itemize}
\item
$\alpha_i \in \half{1}\Z$ such that 
$2\alpha_i \equiv \kappa-1 \bmod 2$ and $\alpha_1 > \dots > \alpha_{u}$;
\item
$\xi_i$ is a unitary character of $\C^\times$ (which can be of the form $\chi_{2\alpha}$);
\item
$u+ 2v=n$.
\end{itemize}
Then we see that
\begin{align*}
\left(\half{1}, \pm1\right) \in X_\lam 
&\iff
\eta(e_{V, 1}) = (-1)^{\#\{ i \in \{1, \dots, u\}\ |\ \alpha_i > 1/2\} + 1}, \\
\left(-\half{1}, \pm1\right) \in X_{\lam^\vee} 
&\iff
\eta(e_{V, 1}) = (-1)^{\#\{ i \in \{1, \dots, u\}\ |\ \alpha_i > 1/2\}}, \\
\left(-\half{1}, -1\right) \in X_\lam
&\iff
\eta(e_{V, -1}) = (-1)^{\#\{ i \in \{1, \dots, u\}\ |\ \alpha_i > -1/2\}+1}, \\
\left(\half{1}, +1\right) \in X_{\lam^\vee}
&\iff
\eta(e_{V, -1}) = (-1)^{\#\{ i \in \{1, \dots, u\}\ |\ \alpha_i > -1/2\}}. 
\end{align*}
In particular, both $\Theta_{r+2, s}(\pi)$ and $\Theta_{r, s+2}(\pi)$ are zero, i.e., 
$(1/2, +1) \in X_\lam \cap X_{\lam^\vee}$ or $(-1/2, -1) \in X_\lam \cap X_{\lam^\vee}$
if and only if
\[
\left\{
\begin{aligned}
&\eta(e_{V, 1}) = (-1)^{\#\{ i \in \{1, \dots, u\}\ |\ \alpha_i > 1/2\} + 1}, \\
&\eta(e_{V, -1}) = (-1)^{\#\{ i \in \{1, \dots, u\}\ |\ \alpha_i > -1/2\}}, 
\end{aligned}
\right.
\quad\text{or}\quad
\left\{
\begin{aligned}
&\eta(e_{V, 1}) = (-1)^{\#\{ i \in \{1, \dots, u\}\ |\ \alpha_i > 1/2\}}, \\
&\eta(e_{V, -1}) = (-1)^{\#\{ i \in \{1, \dots, u\}\ |\ \alpha_i > -1/2\} + 1}.
\end{aligned}
\right.
\]
Since $1/2 \not\in \{\alpha_1, \dots, \alpha_u\}$, 
we see that
\[
\#\{ i \in \{1, \dots, u\}\ |\ \alpha_i > 1/2\} = \#\{ i \in \{1, \dots, u\}\ |\ \alpha_i > -1/2\}.
\]
Hence we see that 
under assuming (chain condition 2) and (even-ness condition), 
both $\Theta_{r+2, s}(\pi)$ and $\Theta_{r, s+2}(\pi)$ are zero
if and only if $\eta(e_{V, 1}) \not= \eta(e_{V, -1})$.
\par

Similarly, when $k > 0$, we see that 
\[
\Theta_{r+2+k, s+k}(\pi) \not=0 \iff \left(\half{k+1}\epsilon, \epsilon\right) \not\in X_\lam
\]
and
\[
\Theta_{r+k, s+2+k}(\pi) \not=0 \iff \left(\half{k+1}\epsilon, \epsilon\right) \not\in X_{\lam^\vee}, 
\]
where $\epsilon$ is the unique element in $\{\pm1\}$ such that
$((k-1)/2, \epsilon), \dots, (-(k-1)/2, \epsilon) \in X_\lam$.
Also, it is easy to see that
if both $\Theta_{r+2+k, s+k}(\pi)$ and $\Theta_{r+k, s+2+k}(\pi)$ are zero, 
then $\lam$ satisfies (chain condition 2) and (even-ness condition).
Furthermore, when $\chi_{k+1}$ appears in $\phi\chi_V^{-1}$ with even multiplicity, 
we see that
\begin{align*}
\left(\half{k+1}\epsilon, \epsilon\right) \in X_\lam 
&\iff 
\eta(e_{V,(k+1)\epsilon}) = (-1)^{\#\{i \in \{1, \dots, u\}\ |\ \alpha_i > (k+1)\epsilon/2\} + 1}, \\
\left(\half{k+1}\epsilon, \epsilon\right) \in X_{\lam^\vee} 
&\iff 
\eta(e_{V,-(k+1)\epsilon}) = (-1)^{\#\{i \in \{1, \dots, u\}\ |\ \alpha_i > -(k+1)\epsilon/2\}}.
\end{align*}
Since $((k-1)/2, \epsilon), \dots, (-(k-1)/2, \epsilon) \in X_\lam$, we have
\[
\eta(e_{V, k+1-2j}) = \epsilon(-1)^{\#\{i \in \{1, \dots, u\}\ |\ \alpha_i > (k+1-2j)\epsilon/2\}}
\]
for any $j = 1, \dots, k$.
Since $\chi_{k+1-2j}$ appears in $\phi\chi_V^{-1}$ with odd multiplicity 
for any $j = 1, \dots, k$ (see (odd-ness condition) in Definition \ref{def} (1)), 
we see that
\[
\#\{i \in \{1, \dots, u\}\ |\ \alpha_i > (k+1-2j)\epsilon/2\}
-
\#\{i \in \{1, \dots, u\}\ |\ \alpha_i > (k-1-2j)\epsilon/2\} = -1
\]
for any $j = 1, \dots, k-1$.
Hence under assuming (chain condition 2) and (even-ness condition), 
both $\Theta_{r+2+k, s}(\pi)$ and $\Theta_{r, s+2+k}(\pi)$ are zero
if and only if
\[
\eta(e_{V,k+1}) \not= \eta(e_{V,k-1}) \not= \dots \not= \eta(e_{V,-(k-1)}) \not= \eta(e_{V,-(k+1)}), 
\]
which is (alternating condition 2).
\par

We have shown that for $k \geq 0$, 
both $\Theta_{r+2+k, s}(\pi)$ and $\Theta_{r, s+2+k}(\pi)$ are zero
if and only if $\lam$ satisfies 
(chain condition 2), (even-ness condition) and (alternating condition 2).
In this case, there exists $\epsilon \in \{\pm1\}$ such that 
$\chi_{(k+1)\epsilon}$ is contained in $\phi\chi_V^{-1}$ with even multiplicity.
Suppose that $((k+1)\epsilon/2, \pm1) \in X_\lam$
(so that $((k-1)/2, \epsilon), \dots, (-(k-1)/2, \epsilon) \in X_\lam$ when $k>0$).
Then by Definition \ref{def} (4) and (5), 
we see that $((k+3)\epsilon/2, \epsilon) \not\in X_\lam^{(\infty)}$, so that
\[
\CC^\epsilon_\lam(k+2) = 
\left\{ 
\left(\half{k+1}\epsilon, \epsilon\right), 
\left(\half{k-1}\epsilon, \epsilon\right), \dots, \left(-\half{k-1}\epsilon, \epsilon\right)
\right\}.
\]
Hence $\#\CC^\epsilon_\lam(k+2) = k+1$.
Moreover, since $((k+1)\epsilon/2, \epsilon) \in X_{\lam^\vee}$ so that 
$(-(k+1)\epsilon/2, -\epsilon) \not\in X_{\lam}$ by Lemma \ref{Xv} (4), 
we see that $\#\CC^{-\epsilon}_\lam(k+2) \leq k+1$.
Hence by Corollary \ref{suf}, we have $\Theta_{r+3+k,s+1+k}(\pi) \not=0$.
Similarly, if $((k+1)\epsilon/2, \pm1) \not\in X_\lam$
(so that $((k-1)/2, -\epsilon), \dots, (-(k-1)/2, -\epsilon) \in X_\lam$ when $k>0$), 
then $(-(k+1)\epsilon/2, \pm1) \in X_{\lam^\vee}$, 
so that we have $\Theta_{r+1+k,s+3+k}(\pi) \not=0$.
In any case, we have
\[
\min\{m_{r-s+2}(\pi), m_{r-s-2}(\pi)\} = r+s+4+2k = n+4+k.
\]
\par

By the conservation relation (Theorem \ref{CR}), we conclude that 
\[
m_{r-s}(\pi) = 
\left\{
\begin{aligned}
&n-2-k \iif \text{$\lam$ satisfies the three conditions},\\
&n-k	\other.
\end{aligned}
\right.
\]
By the tower property (Proposition \ref{tower}), we obtain (2).
\end{proof}

By Corollaries \ref{suf}, \ref{nec} and Proposition, \ref{down}, 
we obtain Theorem \ref{main}.

\appendix
\section{Explicit local Langlands correspondence for discrete series representations}
\label{sec.explicit}
In this appendix, we review the local Langlands correspondence established by 
Langlands himself \cite{L}, Vogan \cite{V3} and Shelstad \cite{S1}, \cite{S2}, \cite{S3}, and
explain the relation between the Harish-Chandra parameters and $L$-parameters
for discrete series representations of unitary groups.
\par

\subsection{Weil groups and representations}
Recall that the Weil group $W_\R$ of $\R$ (\resp $W_\C$ of $\C$)
is defined by 
\[
W_\C = \C^\times, \quad
W_\R = \C^\times \cup \C^\times j
\]
with
\[
j^2 = -1 \in \C^\times,\quad
jzj^{-1} = \overline{z}
\]
for $z \in \C^\times \subset W_\R$.
Then we have an exact sequence
\[
\begin{CD}
1 @>>> W_\C @>>> W_\R @>>> \Gal(\C/\R) @>>>1, 
\end{CD}
\]
where the last map is defined so that $j \mapsto \text{(the complex conjugate)} \in \Gal(\C/\R)$.
Also, the map
\begin{align*}
j \mapsto -1,\quad \C^\times \ni z\mapsto z\overline{z} 
\end{align*}
gives an isomorphism $W_\R^{\mathrm{ab}} \rightarrow \R^\times$.
\par

For $F = \R$ or $F=\C$, a representation of $W_F$ is 
a semisimple continuous homomorphism $\varphi \colon W_F \rightarrow \GL_n(\C)$.
Hence $\varphi$ decomposes into a direct sum of irreducible representations.
\par

For $2\alpha \in \Z$ (i.e., $\alpha \in \half{1}\Z$), 
we define a character $\chi_{2\alpha}$ of $W_\C = \C^\times$ by 
\[
\chi_{2\alpha}(z) = \overline{z}^{-2\alpha}(z\overline{z})^{\alpha}
\]
for $z \in \C^\times$.
A representation $\phi \colon W_\C \rightarrow \GL_n(\C)$ 
is called conjugate self-dual of sign $b \in \{\pm1\}$ if 
there exists a non-degenerate bilinear form $B \colon \C^n \times \C^n \rightarrow \C$
such that
\[
\left\{
\begin{aligned}
&B(\phi(z)x, \phi(\overline{z})y) = B(x,y),\\
&B(y, \phi(-1)x) = b \cdot B(x,y)
\end{aligned}
\right.
\]
for $x,y \in \C^n$ and $z \in W_\C = \C^\times$.
Such a representation $\phi$ is of the form
\[
\phi =
\chi_{2\alpha_1} + \dots + \chi_{2\alpha_{u}} 
+ (\xi_1 + \dots + \xi_v) + ({}^c\xi_1^{-1} + \dots + {}^c\xi_v^{-1}), 
\]
where
\begin{itemize}
\item
$\alpha_i \in \half{1}\Z$ such that $(-1)^{2\alpha_i} = b$;
\item
$\xi_i$ is a character of $\C^\times$ (which can be of the form $\chi_{2\alpha}$);
\item
$u+ 2v=n$.
\end{itemize}
For more precisions, see e.g., \cite[\S 3]{GGP1}.

\subsection{$L$-groups and local Langlands correspondence for unitary groups}
Let $G=\U_n$ be a unitary group of size $n$, 
which is regarded as a connected reductive algebraic group over $\R$.
Hence $G(\R) = \U(p,q)$ for some $(p,q)$ such that $p+q = n$.
Its dual group $\widehat{G}$ is isomorphic to $\GL_n(\C)$.
The Weil group $W_\R = \C^\times \cup \C^\times j$ acts on $\widehat{G} = \GL_n(\C)$
as follows: 
$\C^\times$ acts trivially, and $j$ acts by
\[
j \colon \GL_n(\C) \rightarrow \GL_n(\C),\ 
g \mapsto 
\begin{pmatrix}
&&1\\
&\iddots&\\
(-1)^{n-1}&&
\end{pmatrix}
{}^tg^{-1}
\begin{pmatrix}
&&1\\
&\iddots&\\
(-1)^{n-1}&&
\end{pmatrix}^{-1}.
\]
The $L$-group of $G$ is the semi-direct product
${}^LG = \widehat{G} \rtimes W_\R = \GL_n(\C) \rtimes W_\R$.
\par

An admissible homomorphism of $G(\R) = \U_n(\R)$
is a homomorphism $\varphi \colon W_\R \rightarrow {}^LG$
such that the composition 
\[
\pr \circ \varphi \colon W_\R \rightarrow \GL_n(\C) \rtimes W_\R \twoheadrightarrow W_\R
\]
is identity and the restriction of $\varphi$ to $W_\C = \C^\times$ is continuous.
Let $\Phi(\U_n(\R))$ be the set of ${\widehat{G}}$-conjugacy classes of 
admissible homomorphisms of $\U_n(\R)$.
For $\varphi \in \Phi(\U_n(\R))$, we define the component group $A_\varphi$ of $\varphi$ by 
\[
A_\varphi = \pi_0(\cent(\im(\varphi), \widehat{G})).
\]
This is an elementary two abelian group.
For $\varphi \in \Phi(\U_n(\R))$, the restriction of $\varphi$ gives a conjugate self-dual representation
$\phi = \varphi|\C^\times$ of $W_\C$ of dimension $n$ and sign $(-1)^{n-1}$.
Via the map $\varphi \mapsto \phi = \varphi | \C^\times$, 
we obtain an identification
\[
\Phi(\U_n(\R)) = 
\{\text{conjugate self-dual representations of $W_\C$ of dimension $n$ and sign $(-1)^{n-1}$}\}.
\]
When $\phi = \varphi|\C^\times$, we also put $A_\phi = A_\varphi$.
\par

We say that $\phi \in \Phi(\U_n(\R))$ is discrete (\resp tempered) if
$\phi$ is of the form $\phi = \chi_{2\alpha_1} \oplus \dots \oplus \chi_{2\alpha_n}$
with $2\alpha_i \equiv n-1 \bmod 2$ and $\alpha_1 > \dots > \alpha_n$
(\resp $\phi$ is a direct sum of unitary characters of $\C^\times$).
We denote the subset of $\Phi(\U_n(\R))$ consisting of
discrete elements (\resp tempered elements) 
by $\Phi_\disc(\U_n(\R))$ (\resp $\Phi_\temp(\U_n(\R))$).
\par

When $\phi = \varphi|\C^\times = \chi_{2\alpha_1} \oplus \dots \oplus \chi_{2\alpha_n} \in \Phi_\disc(\U_n(\R))$, 
there exists a unique semisimple element $s_{2\alpha_i} \in \cent(\im(\varphi), \GL_n(\C))$ 
such that $W_\C$ acts on the $(-1)$-eigenspace of $s_{2\alpha_i}$ by $\chi_{2\alpha_i}$.
Let $e_{2\alpha_i}$ be the image of $s_{2\alpha_i}$ in $A_\phi = \pi_0(\cent(\im(\varphi), \GL_n(\C)))$.
Then we have
\[
A_\phi = (\Z/2\Z)e_{2\alpha_1} \oplus \dots \oplus (\Z/2\Z)e_{2\alpha_n}.
\]
For more precisions, see \cite[\S 4]{GGP1}.
In particular, we have $|A_\phi| = 2^n$ for each $\phi \in \Phi_\disc(\U_n(\R))$.
\par

The local Langlands correspondence for unitary groups is as follows:
\begin{thm}\label{LLC2}
\begin{enumerate}
\item
There exists a canonical surjection
\[
\bigsqcup_{p+q=n}\Irr_\temp(\U(p,q)) \rightarrow \Phi_\temp(\U_n(\R)).
\]
For $\phi \in \Phi_\temp(\U_n(\R))$, 
we denote by $\Pi_\phi$ the inverse image of $\phi$ under this map, 
and call $\Pi_\phi$ the $L$-packet associated to $\phi$.
\item
$\#\Pi_\phi = \#A_\phi$;
\item
$\pi \in \Pi_{\phi}$ is discrete series if and only if $\phi$ is discrete.
\item
The map $\pi \mapsto \phi$ is compatible with 
parabolic inductions (c.f., Theorem \ref{LLC} (5)).
\item
If $\phi = \chi_{2\alpha_1} \oplus \dots \oplus \chi_{2\alpha_n} \in \Phi_\disc(\U_n(\R))$, 
then $\Pi_\phi$ is the set of all discrete series representations of various $\U(p,q)$ 
whose infinitesimal characters are equal to $(\alpha_1, \dots, \alpha_n)$
via the Harish-Chandra map.
\end{enumerate}
\end{thm}
Note that there exist exactly $(p+q)!/(p! \cdot q!)$ discrete series representations 
of $\U(p,q)$ with a given infinitesimal character.
Theorem \ref{LLC} (2) and (5) are compatible with the well-known equation
\[
\sum_{p+q = n} \frac{(p+q)!}{p! \cdot q!} = 2^n.
\]
\par

\subsection{Whittaker data and generic representations}
The unitary group $\U(p,q)$ with $p+q = n$ is quasi-split if and only if $|p-q| \leq 1$.
For such $(p,q)$, 
a Whittaker datum of $\U(p,q)$ is the conjugacy class of pairs $\w = (B,\mu)$, 
where $B=TU$ is an $\R$-rational Borel subgroup of $\U(p,q)$, and
$\mu \colon U(\R) \rightarrow \C^\times$ is a unitary generic character.
Here, $T$ is a maximal $\R$-torus of $B$ and $U$ is the unipotent radical of $B$, and
$T(\R)$ acts on $U(\R)$ by conjugation.
A unitary character $\mu$ of $U(\R)$ is called generic if 
the stabilizer of $\mu$ in $T(\R)$ is equal to the center $Z(\R)$ of $\U(p,q)$.
We say that $\pi \in \Irr_\temp(\U(p,q))$ is $\w$-generic 
if 
\[
\Hom_{\U(p,q)}(\pi, C^{\infty}_{\mg}(U(\R) \bs \U(p,q), \mu)) \not=0,
\] 
where $C^{\infty}_{\mg}(U(\R) \bs \U(p,q), \mu)$ is 
the set of moderate growth $C^{\infty}$-functions $W \colon \U(p,q) \rightarrow \C$ satisfying that
$W(ug) = \mu(u)W(g)$ for $u \in U(\R)$ and $g \in \U(p,q)$, 
and $\U(p,q)$ acts on $C^{\infty}_{\mg}(U(\R) \bs \U(p,q), \mu)$ by the right translation.
\par

For each $\phi \in \Phi_\temp(\U_n(\R))$, 
the $L$-packet $\Pi_\phi$ is parametrized the Pontryagin dual $\widehat{A_\phi}$ of $A_\phi$
if a quasi-split form $\U_n(\R)$ and its Whittaker datum are fixed.

\begin{thm}\label{LLC3}
Fix $(p,q)$ such that $p+q = n$ and $|p-q| \leq 1$, and 
a Whittaker datum $\w$ of $\U(p,q)$.
Then 
\begin{enumerate}
\item
for $\phi \in \Phi_\temp(\U_n(\R))$, 
there exists a bijection
\[
J_\w \colon \Pi_\phi \rightarrow \widehat{A_\phi}
\]
which satisfies certain character identities;
\item
for each $\phi \in \Phi_\temp(\U_n(\R))$, 
the $L$-packet $\Pi_\phi$ has a unique $\w$-generic representation $\pi_\w$; 
\item
in particular, the bijection $J_\w$ requires satisfying that 
$J_\w(\pi_\w)$ is the trivial character of $A_\phi$.
\end{enumerate}
\end{thm}
\par

In the next subsection, 
we will review the definition of $J_\w$ when $\phi \in \Phi_\disc(\U_n(\R))$, 
and give an explicit relation between $J_\w(\pi)$ and the Harish-Chandra parameter of $\pi$
for $\pi \in \Pi_\phi$.
To give a such relation, we need to specify which representation is $\w$-generic.
\par

Fix $(p,q)$ such that $p+q = n$ and $|p-q| \leq 1$.
By \cite[\S6, Theorem 6.2]{V1}, 
for $\pi \in \Irr_\disc(\U(p,q))$ with Harish-Chandra parameter $\lam$, 
the following are equivalent;
\begin{itemize}
\item
$\pi$ is large;
\item
$\pi$ is $\w$-generic for some Whittaker datum of $\U(p,q)$;
\item
all simple roots in $\Delta^+_\lam$ are non-compact, i.e., 
do not belong to $\Delta(K_{p,q},T_{p,q})$.
\end{itemize}
Here, 
\begin{itemize}
\item
$K_{p,q} \cong \U(p) \times \U(q)$ is the usual maximal compact subgroup of $\U(p,q)$; 
\item
$T_{p,q}$ is the usual maximal compact torus of $\U(p,q)$;
\item
$\Delta(\U(p,q), T_{p,q})$ (\resp $\Delta(K_{p,q}, T_{p,q})$) is the set of 
roots of $T_{p,q}$ in $\U(p,q)$ (\resp in $K_{p,q}$);
\item
$\Delta_{\lam}^+$ is the unique positive system of $\Delta(\U(p,q), T_{p,q})$
for which $\lam$ is dominant.
\end{itemize}
\par

When $n$ is odd, there exist exactly two quasi-split forms 
$\U((n+1)/2, (n-1)/2)$ and $\U((n-1)/2, (n+1)/2)$.
For $\epsilon \in \{\pm1\}$, we put 
$(p_\epsilon, q_\epsilon) = ((n+\epsilon)/2, (n-\epsilon)/2)$.
Then there exists a unique Whittaker datum $\w_\epsilon$ of $\U(p_\epsilon,q_\epsilon)$.
For integers $\alpha_1 > \dots > \alpha_n$, it is easy to see that
the Harish-Chandra parameter of
the unique large discrete series representation $\pi_{\w_\pm}$ of $\U(p_\pm, q_\pm)$
with infinitesimal character $(\alpha_1, \dots, \alpha_n)$ 
is given by
\[
\left\{
\begin{aligned}
&\HC(\pi_{\w_+}) = (\alpha_1, \alpha_3, \dots, \alpha_n; \alpha_2, \alpha_4, \dots, \alpha_{n-1}),\\
&\HC(\pi_{\w_-}) = (\alpha_2, \alpha_4, \dots, \alpha_{n-1}; \alpha_1, \alpha_3, \dots, \alpha_n).
\end{aligned}
\right.
\]
\par

When $n = 2m$ is even, $\U(m,m)$ is the unique quasi-split form of size $n$.
It has exactly two Whittaker data $\w_\pm$ constructed as follows.
First, we set
\begin{align*}
G_m &= \left\{g \in \GL_{2m}(\C)\ |\ 
{}^t\overline{g}
\begin{pmatrix}
\1_{m} & 0\\
0 & -\1_{m}
\end{pmatrix}
g
= 
\begin{pmatrix}
\1_{m} & 0\\
0 & -\1_{m}
\end{pmatrix}
\right\}, \\
G'_m &= \left\{g' \in \GL_{2m}(\C)\ |\ 
{}^t\overline{g'}
\begin{pmatrix}
&& 1\\
&\iddots&\\
1&&
\end{pmatrix}
g'
= 
\begin{pmatrix}
&& 1\\
&\iddots&\\
1&&
\end{pmatrix}
\right\}.
\end{align*}
Note that $G_m$ is the usual coordinate of $\U(m,m)$, and
these two groups are isomorphic to each other.
Put
\[
T_{m}=\frac{1}{\sqrt{2}}
\left(
\begin{array}{ccc|ccc}
1&& &&&-1\\
&\ddots& &&\iddots&\\
&&1 &-1&&\\ \hline
&&1 &1&&\\
&\iddots& &&\ddots&\\
1&& &&&1
\end{array}
\right)
\in \GL_{2m}(\C),
\]
so that
\[
T_{m}^{-1}=\frac{1}{\sqrt{2}}
\left(
\begin{array}{ccc|ccc}
1&& &&&1\\
&\ddots& &&\iddots&\\
&&1 &1&&\\ \hline
&&-1 &1&&\\
&\iddots& &&\ddots&\\
-1&& &&&1
\end{array}
\right) \in \GL_{2m}(\C).
\]
Then
\[
{}^t\overline{T_m}
\begin{pmatrix}
&& 1\\
&\iddots&\\
1&&
\end{pmatrix}
T_m
=
\begin{pmatrix}
\1_{m} & 0\\
0 & -\1_{m}
\end{pmatrix}, 
\]
so that the map
\[
f_{T_m} \colon G_m \rightarrow G'_m,\ g \mapsto g' \coloneqq T_m g T_m^{-1}
\]
gives an isomorphism.
Let $B'=T'U'$ be the Borel subgroup of $G'_m$ consisting of upper triangular matrices, 
where $T'$ is the maximal torus of $G'_m$ consisting of diagonal matrices and 
$U'$ is the unipotent radical of $B'$.
Define a generic character $\mu_\pm$ of $U'(\R)$ by
\[
\mu_\pm(u) = \exp(\mp\pi\I \tr_{\C/\R}(\I (u_{1,2} + \dots + u_{m,m+1}))).
\]
We set the Whittaker datum $\w_\pm$ of $G_m$ to be the conjugacy class of
\[
\w_\pm = (f_{T_m}^{-1}(B'), \mu_\pm \circ f_{T_m}).
\]
\par

Fix half-integers $\alpha_1 > \dots > \alpha_n$. 
Let $\pi$ and $\pi'$ be the discrete series representations with
Harish-Chandra parameters $\lam$ and $\lam'$ given by
\[
\left\{
\begin{aligned}
&\lam=(\alpha_1, \alpha_3, \dots, \alpha_{n-1}; \alpha_2, \alpha_4, \dots, \alpha_{n}),\\
&\lam'=(\alpha_2, \alpha_4, \dots, \alpha_{n}; \alpha_1, \alpha_3, \dots, \alpha_{n-1}),
\end{aligned}
\right.
\]
respectively. 
Then $\pi$ and $\pi'$ are the two large discrete series representations
with infinitesimal character $(\alpha_1, \dots, \alpha_n)$. 
Hence there exists $\epsilon \in \{\pm1\}$ such that 
$\pi$ (\resp $\pi'$) is $\w_{\epsilon}$-generic (\resp $\w_{-\epsilon}$-generic).
\par

To give an explicit description of the local Langlands correspondence for $\U_n(\R)$, 
we have to determine $\epsilon$.
The following proposition says that $\epsilon = +1$.
It seems to be well-known (c.f., \cite{M2}), 
but we give a proof for  the convenience of the reader.

\begin{prop}\label{generic}
Assume that $n = 2m$ is even.
Fixing half-integers $\alpha_1 > \dots > \alpha_n$, 
we let $\pi$ (\resp $\pi'$) be the large discrete series representation of $\U(m,m)$
whose Harish-Chandra parameter $\lam$ (\resp $\lam'$) is given as above.
Then $\pi$ is $\w_+$-generic (\resp $\pi'$ is $\w_-$-generic).
\end{prop}
\begin{proof}
We prove the proposition by induction on $m$.
First suppose that $m=1$ so that $G_1=\U(1,1)$.
Note that
\begin{align*}
\g = \Lie(G_1) 
&= \left\{
\begin{pmatrix}
a\I & b + c\I\\
b-c\I & d\I
\end{pmatrix}
\in \M_2(\C)\ |\ 
a,b,c,d \in \R
\right\}.
\end{align*}
Set 
\[
H=
\begin{pmatrix}
1 & 0 \\ 0 & -1
\end{pmatrix},
\quad
X_+ = 
\begin{pmatrix}
0&1\\0&0
\end{pmatrix},
\quad
X_- = 
\begin{pmatrix}
0&0\\1&0
\end{pmatrix}
\in \g_\C = \M_2(\C).
\]
Then $(H,X_+,X_-)$ is an $\mathfrak{sl}_2$-triple, i.e., 
\[
[X_+,X_-] = H,
\quad
[H,X_+] = 2X_+,
\quad
[H, X_-] = -2X_-.
\]
Set 
\[
n(z) = f_1^{-1}(
\begin{pmatrix}
1 & z \\ 0 & 1
\end{pmatrix}
)
= 
\begin{pmatrix}
1+z/2 & z/2 \\ -z/2 & 1-z/2
\end{pmatrix},
\quad
m(x) = f_1^{-1}(
\begin{pmatrix}
x & 0 \\ 0 & x^{-1}
\end{pmatrix}
)
=\half{1}
\begin{pmatrix}
x+x^{-1} & -x+x^{-1}\\
-x+x^{-1} & x+x^{-1}
\end{pmatrix}
\]
for $z \in \I \R$ and $x \in \R_{>0}$.
By the Iwasawa decomposition, any $g \in G$ has a unique decomposition
\[
g = n(z)m(x)
\begin{pmatrix}
t_1 & 0 \\ 0 & t_2
\end{pmatrix}
\]
for $t_1,t_2 \in \C^1$, $z \in \I \R$ and $x \in \R_{>0}$.
Define a moderate growth $C^\infty$-function $W$ on $G_1$ by
\[
W(n(z)m(x)
\begin{pmatrix}
t_1 & 0 \\ 0 & t_2
\end{pmatrix}
)
= \exp(-\pi\I\tr_{\C/\R}(\I z)) 
\cdot x^{\alpha_1-\alpha_2+1} e^{-2\pi x^2} 
\cdot t_1^{\alpha_1+1/2}t_2^{\alpha_2-1/2}
\]
for $t_1,t_2 \in \C^1$, $z \in \I \R$ and $x \in \R_{>0}$.
Then it is easy to see that
\[
X_- \cdot W = 0, 
\]
so that $W$ generates the discrete series representation of $G_1$ 
whose Harish-Chandra parameter is $\lam = (\alpha_1; \alpha_2)$.
We conclude that $\pi$ is $\w_+$-generic, so that $\pi'$ is $\w_-$-generic.
\par

Next, we assume that $m > 1$.
Let $P = f_{T_m}^{-1}(P') = MN$ be the parabolic subgroup of $G_m=\U(m,m)$
given by
\[
P' = \left\{
\begin{pmatrix}
a & * & * \\
0 & g_0' & * \\
0 & 0 & \overline{a}^{-1}
\end{pmatrix}
\ |\ 
a \in \C^\times,\ g_0' \in G'_{m-1}
\right\}.
\]
Here, $M = f_{T_m}^{-1}(M')$ is the Levi subgroup of $P$ such that
$M'$ consists of the block diagonal matrices of $P'$.
Let $\pi_0$ be the discrete series representation of $\U(m-1,m-1)$
whose Harish-Chandra parameter $\lam_0$ is given by
\[
\lam_0 = (\alpha_3, \alpha_5, \dots, \alpha_{n-1}; \alpha_4, \alpha_6, \dots, \alpha_n), 
\]
and $\chi$ be the character of $\C^\times$ given by
\[
\chi(ae^{\theta\I}) = a^{\alpha_1-\alpha_2} e^{(\alpha_1+\alpha_2)\theta\I}
\]
for $a>0$ and $\theta \in \R/2\pi\Z$.
Consider the normalized induced representation $I(\pi_0) = \Ind_P^G(\chi \boxtimes \pi_0)$.
We denote the Harish-Chandra characters of $\pi$ and $\Ind_P^G(\chi \boxtimes \pi_0)$ by
\[
\Theta_\lam
\quad \text{and} \quad
\Theta_{I(\pi_0)}, 
\]
respectively.
Now we use Schmid's character identity (\cite[(9.4) Theorem]{Sc1}, \cite[Theorem (b)]{Sc2}).
This focuses on representations of semisimple groups, but 
since $\U(m,m)$ is generated by its center and a semisimple group $\mathrm{SU}(m,m)$, 
we can apply Schmid's character identity to $\pi$.
It asserts that
\[
\Theta_{\lam} + \Theta'_{\lam} = \Theta_{I(\pi_0)}, 
\]
where $\Z^{2m} \ni \mu \mapsto \Theta'_{\mu}$ is 
the coherent continuation of Harish-Chandra characters satisfying that:
\begin{itemize}
\item
$\Theta'_{\mu}$ is a virtual character corresponding to 
the infinitesimal character $\mu$ via the Harish-Chandra map;
\item
if $\mu = (\mu_1, \dots, \mu_m; \mu_{m+1}, \dots, \mu_{2m}) \in \Z^{2m}$ 
satisfies that
\[
\mu_{m+1} > \mu_1 > \mu_2, 
\quad\text{and}\quad
\mu_2 > \mu_{m+2} > \mu_3 > \mu_{m+3} > \dots > \mu_{m} > \mu_{2m}, 
\]
then $\Theta'_{\mu}$ is the character of the discrete series representation 
whose Harish-Chandra parameter is $\mu$.
\end{itemize}
By the theory of the wall-crossing of coherent families (see e.g., \cite[Corollary 7.3.9]{V2}), 
we see that $\Theta'_{\lam}$ is the character of a representation of $G$.
This implies that $\pi$ is a subquotient of $I(\pi_0)$.
By the additivity of Whittaker model (see e.g., \cite[\S3]{M1}) 
and Hashizume's result (\cite[Theorem 1]{Ha}, \cite[Theorem 3.4.1]{M2}), 
we see that if $\pi$ is $\w_\epsilon$-generic, then so is $\pi_0$.
By the induction hypothesis, we must have $\epsilon = +1$, as desired.
\end{proof}

\subsection{Explicit description of discrete $L$-packets}
In this subsection, we recall the definition of $J_{\w_+}$ in Theorem \ref{LLC3}
when $\phi \in \Phi_\disc(\U_n(\R))$ (see e.g., \cite[\S 5.6]{Ka}),
and explain Theorem \ref{LLC} (4).
\par

Let $G=\U(p,q)$. 
We set $G^* = \U(m,m)$ and $(p_0,q_0) = (m,m)$ if $p+q=2m$, 
and $G^* = \U(m+1,m)$ and $(p_0,q_0) = (m+1,m)$ if $p+q=2m+1$.
We choose an isomorphism $\psi_{p,q} \colon G^* \rightarrow G$ over $\C$
and a $1$-cocycle $z_{p,q} \in Z^1(\R, G^*)$
by 
\[
\left\{
\begin{aligned}
&\psi_{p,q} = \Int
\begin{pmatrix}
\1_p &&\\
&\I\1_{p_0-p}&\\
&&\1_{q_0}
\end{pmatrix},
&z_{p,q}(c) = 
\begin{pmatrix}
\1_p &&\\
&-\1_{p_0-p}&\\
&&\1_{q_0}
\end{pmatrix}
\iif p \leq p_0,\\
&\psi_{p,q} = \Int
\begin{pmatrix}
\1_{p_0} &&\\
&\I\1_{q_0-q}&\\
&&\1_q
\end{pmatrix},
&z_{p,q}(c) = 
\begin{pmatrix}
\1_{p_0} &&\\
&-\1_{q_0-q}&\\
&&\1_q
\end{pmatrix}
\iif q \leq q_0.
\end{aligned}
\right.
\]
Here, we denote by $c$ the complex conjugate in $\Gal(\C/\R)$.
Then we have
\[
\psi_{p,q}^{-1} \cdot c(\psi_{p,q}) = \Int(z_{p,q}(c)).
\]
Hence $(G, \psi_{p,q}, z_{p,q})$ is a pure inner twist of $G^*$.
In particular, $G^*$ is a quasi-split pure inner form of $G$.
\par

Let $\phi = \chi_{2\alpha_1} \oplus \dots \oplus \chi_{2\alpha_n} \in \Phi_\disc(\U_n(\R))$ 
with $n=p+q$.
We can realize $\phi$ as
\[
\phi(ae^{\I\theta}) = 
\begin{pmatrix}
e^{2\alpha_1\I\theta} &&\\
&\ddots& \\
&&e^{2\alpha_{n}\I\theta} \\
\end{pmatrix}
\]
for $a \in \R_{>0}$ and $\theta \in \R/2\pi\Z$.
Note that
if $\phi = \varphi|W_\C$ for an admissible homomorphism 
$\varphi \colon W_\R \rightarrow {}^L\U_n$, 
then $\varphi(j)$ acts on $\phi(ae^{\I\theta})$ by the inverse.
\par

We denote the canonical right action of $S_n$ on $(\C^1)^n$ by
\[
(t_1, \dots, t_n)^\sigma = (t_{\sigma(1)}, \dots, t_{\sigma(n)})
\]
for $\sigma \in S_{n}$ and $(t_1, \dots, t_n) \in (\C^1)^n$.
For $\sigma \in S_{n}$, we define an embedding
\[
\eta_{p,q}^\sigma \colon (\C^1)^n \hookrightarrow G= \U(p,q)
\]
by
\[
\eta_{p,q}^\sigma(t_1, \dots, t_n) = 
\begin{pmatrix}
t_{\sigma(1)}  &&     &          &&\\
&\ddots&                  &         &&\\
&&t_{\sigma(p)}  &         &&\\
&&                            & t_{\sigma(p+1)}&&\\
&&                            &        &\ddots&\\
&&                             &       &&t_{\sigma(n)}
\end{pmatrix}
\in G= \U(p,q).
\]
for $t_1, \dots, t_n \in \C^1$.
We call $\eta_{p,q}^\sigma$ an admissible embedding of $(\C^1)^n$ (see \cite[\S 5.6]{Ka}).
The image of $\eta_{p,q}^\sigma$ is independent of $\sigma$, and is denoted by $T_{p,q}$.
Note that $\eta_{p,q}^\sigma$ and $\eta_{p,q}^{\sigma'}$ are $\U(p,q)$-conjugate if and only if
\[
\sigma^{-1}\sigma'  \in S_{p} \times S_{q}.
\]
Hence when we consider the $\U(p,q)$-conjugacy class of $\eta_{p,q}^\sigma$, 
we may assume that
\[
\sigma(1) < \dots < \sigma(p),\quad
\sigma(p+1) < \dots < \sigma(n).
\]
For such $\eta_{p,q}^\sigma$, 
we denote by $\pi_{p,q}^\sigma$ the irreducible discrete series representation of $\U(p,q)$
with Harish-Chandra parameter
\[
\lam^\sigma = 
(\alpha_{\sigma(1)}, \dots, \alpha_{\sigma(p)}; 
\alpha_{\sigma(p+1)}, \dots, \alpha_{\sigma(n)}).
\]
Let $\Theta_{\pi_{p,q}^\sigma}$ be the Harish-Chandra character of $\pi_{p,q}^\sigma$, 
which is a real analytic function on the regular set $\U(p,q)^\reg$ of $\U(p,q)$.
We put $K_{p,q} = \U(p) \times \U(q)$ to be 
the usual maximal compact subgroup of $\U(p,q)$,
which contains $T_{p,q}$ as a maximal torus.
On $T_{p,q,\reg} \coloneqq \U(p,q)^\reg \cap T_{p,q}$, 
we have
\[
\Theta_{\pi_{p,q}^\sigma}(t) = (-1)^{\half{1}\dim(\U(p,q)/K_{p,q})}
\frac{\sum_{w \in W_{K_{p,q}}}\sgn(w)t^{w(\lam^\sigma)}}
{\prod_{\alpha \in \Delta_{\lam^\sigma}^+}(t^{\alpha/2}-t^{-\alpha/2})}, 
\]
where 
$\Delta_{\lam^\sigma}^+$ is the unique positive system of $\Delta(\U(p,q), T_{p,q})$
such that $\pair{\lam^\sigma, \alpha^\vee} > 0$ for any $\alpha \in \Delta_{\lam^\sigma}^+$.
Note that $\dim(\U(p,q)/K_{p,q}) = 2pq$ so that $(-1)^{\half{1}\dim(\U(p,q)/K_{p,q})} = (-1)^{pq}$.
We put
\[
\Pi_\phi^{\U(p,q)} = \{\pi_{p,q}^\sigma\ |\ \sigma \in S_n/(S_{p} \times S_{q})\}
\subset \Irr_\disc(\U(p,q)).
\]
Then the $L$-packet $\Pi_\phi$ associated to $\phi$ is defined by
\[
\Pi_\phi = \bigsqcup_{p+q = n} \Pi_\phi^{\U(p,q)}.
\]
\par

Recall that 
\[
A_\phi = (\Z/2\Z) e_{2\alpha_1} \oplus \dots \oplus (\Z/2\Z) e_{2\alpha_n}.
\]
Let $\pi_{p,q}^\sigma \in \Pi_{\phi}^{\U(p,q)}$ with $\sigma \in S_{n}$ satisfying 
$\sigma(1) < \dots < \sigma(p)$ and $\sigma(p+1) < \dots < \sigma(n)$.
We define $g_{p,q}^\sigma \in \U(p_0,q_0, \C)$ as follows.
There is an element $h_{p,q}^\sigma \in \U(n,0)$ such that
$\Int(h_{p,q}^\sigma)$ is equal to 
\[
\begin{pmatrix}
t_1 &&\\
&\ddots&\\
&&t_n
\end{pmatrix}
\mapsto 
\begin{pmatrix}
t_{\sigma_+^{-1}\sigma(1)} &&\\
&\ddots&\\
&&t_{\sigma_+^{-1}\sigma(n)}
\end{pmatrix}
\]
on $T_{n,0}$.
Here, $\sigma_+ \in S_n$ is defined by
\[
\sigma_+(i) = \left\{
\begin{aligned}
&2i-1 \iif i \leq p_0,\\
&2(i-p_0) \iif i > p_0
\end{aligned}
\right.
\]
so that
\[
(\alpha_{\sigma_+(1)}, \dots, \alpha_{\sigma_+(p_0)}; 
\alpha_{\sigma_+(p_0+1)}, \dots, \alpha_{\sigma_+(n)})
=
(\alpha_{1}, \alpha_3, \dots ; \alpha_{2}, \alpha_4, \dots).
\]
We put 
\[
g_{p,q}^\sigma = 
\begin{pmatrix}
\1_{p_0}&0\\0&\I\1_{q_0}
\end{pmatrix}^{-1}
h_{p,q}^\sigma
\begin{pmatrix}
\1_{p_0}&0\\0&\I\1_{q_0}
\end{pmatrix} \in \U(p_0, q_0, \C).
\]
Then we have
\[
\eta_{p,q}^\sigma = 
\psi_{p,q} \circ \Int(g_{p,q}^\sigma) \circ \eta_{p_0,q_0}^{\sigma_+}.
\]
We set $\inv(\pi_{p_0,q_0}^{\sigma_+}, \pi_{p,q}^\sigma) \in H^1(\R, T_{p_0,q_0})$ 
to be the class of
\[
c \mapsto (g_{p,q}^\sigma)^{-1} \cdot z_{p,q}(c) \cdot c(g_{p,q}^\sigma).
\]
We can compute $(g_{p,q}^\sigma)^{-1} \cdot z_{p,q}(c) \cdot c(g_{p,q}^\sigma)$ explicitly.
First, we have
\[
\begin{pmatrix}
\1_{p_0}&0\\0&\I\1_{q_0}
\end{pmatrix}
z_{p,q}(c)
\begin{pmatrix}
\1_{p_0}&0\\0&\I\1_{q_0}
\end{pmatrix}
=
\begin{pmatrix}
\1_p & 0 \\ 0 & -\1_q
\end{pmatrix}.
\]
For $i=1,\dots, n$, 
we define $\epsilon_{i} \in \{\pm1\}$ by
\[
\epsilon_i = \left\{
\begin{aligned}
&+1 \iif i \leq p,\\
&-1 \iif i > p.
\end{aligned}
\right.
\]
Then $(g_{p,q}^\sigma)^{-1} \cdot z_{p,q}(c) \cdot c(g_{p,q}^\sigma)$ is equal to
\[
\begin{pmatrix}
\epsilon_{\sigma^{-1}\sigma_+(1)}&&&&&\\
&\ddots &&&&\\
&& \epsilon_{\sigma^{-1}\sigma_+(p_0)}&&&\\
&&&-\epsilon_{\sigma^{-1}\sigma_+(p_0+1)}&&\\
&&&&\ddots&\\
&&&&&-\epsilon_{\sigma^{-1}\sigma_+(n)}
\end{pmatrix}.
\]
\par

Note that $H^1(\R, \C^1) \cong \{\pm1\}$.
The map $\eta_{p_0,q_0}^{\sigma_+} \colon (\C^1)^n \rightarrow T_{p_0,q_0}$
gives an isomorphism 
\[
(\eta_{p_0,q_0}^{\sigma_+})^* \colon H^1(\R, (\C^1)^n) \rightarrow H^1(\R, T_{p_0,q_0}).
\]
For each $e_{2\alpha_i} \in A_\phi$, 
we denote by $\pair{e_{2\alpha_i}, \cdot}$ the $i$-th projection
\[
H^1(\R, (\C^1)^n) \cong H^1(\R, \C^1)^n \cong \{\pm1\}^n \rightarrow \{\pm1\}.
\]
Then $J_{\w_+}(\pi_{p,q}^\sigma)(e_{2\alpha_i}) \in \{\pm1\}$ is defined by
\[
J_{\w_+}(\pi_{p,q}^\sigma)(e_{2\alpha_i}) = \pair{e_{2\alpha_i}, 
{(\eta_{p_0,q_0}^{\sigma_+})^*}^{-1}(\inv(\pi_{p_0,q_0}^{\sigma_+}, 
\pi_{p,q}^\sigma))}^{-1}.
\]
Hence
\begin{align*}
J_{\w_+}(\pi_{p,q}^\sigma)(e_{2\alpha_{\sigma_+(i)}}) 
&= 
\left\{
\begin{aligned}
&\epsilon_{\sigma^{-1}\sigma_+(i)} \iif i \leq p_0,\\
&-\epsilon_{\sigma^{-1}\sigma_+(i)} \iif i>p_0.
\end{aligned}
\right.
\end{align*}
If we put $j= \sigma_+(i)$, we see that $1 \leq i \leq p_0$ if and only if $j$ is odd.
Therefore, we conclude that
\[
J_{\w_+}(\pi_{p,q}^\sigma)(e_{2\alpha_{j}}) 
=\left\{
\begin{aligned}
&(-1)^{j-1} \iif \sigma^{-1}(j) \leq p,\\
&(-1)^j \iif \sigma^{-1}(j) > p.
\end{aligned}
\right.
\]
\par

The other bijection $J_{\w_-} \colon \Pi_\phi \rightarrow \widehat{A_\phi}$
is defined by using $\sigma_- \in S_n$ such that
\[
(\alpha_{\sigma_-(1)}, \dots, \alpha_{\sigma_-(p_0)}; 
\alpha_{\sigma_-(p_0+1)}, \dots, \alpha_{\sigma_-(n)})
=
(\alpha_{2}, \alpha_4, \dots ; \alpha_{1}, \alpha_3, \dots)
\]
in place of $\sigma_+$.
By a similar calculation, we have
\[
J_{\w_-}(\pi_{p,q}^\sigma)(e_{2\alpha_{j}}) 
=\left\{
\begin{aligned}
&(-1)^{j} \iif \sigma^{-1}(j) \leq p,\\
&(-1)^{j-1} \iif \sigma^{-1}(j) > p.
\end{aligned}
\right.
\]
In particular, the isomorphism 
$J_{\w_-} \circ (J_{\w_+})^{-1} \colon \widehat{A_\phi} \rightarrow \widehat{A_\phi}$
is given by
\[
J_{\w_-} \circ (J_{\w_+})^{-1} = \cdot \otimes \eta_{-1}, 
\]
where the character $\eta_{-1} \colon A_{\phi} \rightarrow \{\pm1\}$ is defined by
\[
\eta_{-1}(e_{2\alpha_j}) = -1
\]
for any $1 \leq j \leq n$.
\par

Hence we have the following theorem.
\begin{thm}
Let 
\[
\phi = \chi_{2\alpha_1} \oplus \dots \oplus \chi_{2\alpha_n} \in \Phi_\disc(\U_n(\R))
\]
with $2\alpha_i \in \Z$ satisfying $\alpha_1 > \dots > \alpha_n$ 
and $2\alpha_i \equiv n-1 \bmod 2$.
Then the $L$-packet $\Pi_\phi$ consists of 
all irreducible discrete series representations of 
$\U(p,q)$ with $p+q = n$
whose infinitesimal characters are equal to $(\alpha_1, \dots, \alpha_n)$.
Moreover, for each Whittaker datum $\w_\pm$ of $\U(p_0,q_0)$ with $|p_0-q_0| \leq 1$, 
there is a bijection
\[
J_{\w_\pm} \colon \Pi_\phi \rightarrow \widehat{A_\phi}
\]
such that 
the Harish-Chandra parameter of $\pi^\pm(\phi,\eta) \coloneqq (J_{\w_\pm})^{-1}(\eta)$
is given by $(\lam_1, \dots, \lam_{p}; \lam'_1, \dots, \lam'_{q})$, 
where
\begin{align*}
\{\lam_1, \dots, \lam_{p}\} &= \{\alpha_i\ |\ \eta(e_{2\alpha_i}) = \pm(-1)^{i-1}\},\\
\{\lam'_1, \dots, \lam'_{q}\} &= \{\alpha_i\ |\ \eta(e_{2\alpha_i}) = \pm(-1)^{i}\}.
\end{align*}
In particular, if we put $p = \#\{i\ |\ \eta(e_{2\alpha_i}) = \pm(-1)^{i-1}\}$ and
$q = \#\{i\ |\ \eta(e_{2\alpha_i}) = \pm(-1)^{i}\}$, 
then 
$\pi^\pm(\phi,\eta)$ is a representation of $\U(p,q)$.
There is a unique $\w_\pm$-generic representation in $\Pi_\phi$
which corresponds to
the trivial character of $A_\phi$ via $J_{\w_\pm}$.
The bijections $J_{\w_+}$ and $J_{\w_-}$ are related by
\[
J_{\w_-}(\pi) = J_{\w_+}(\pi) \otimes \eta_{-1}
\]
for any $\pi \in \Pi_\phi$, 
where $\eta_{-1} \in \widehat{A_\phi}$ is defined by
$\eta_{-1}(e_{2\alpha_i}) = -1$ for any $e_{2\alpha_i} \in A_\phi$.
\end{thm}
In this paper, we always use $J=J_{\w_+}$.
By this theorem, we obtain Theorem \ref{LLC} (4).


\end{document}